\documentclass[12pt]{amsart}
\usepackage{amssymb,latexsym,amscd}

\numberwithin{equation}{section}

\newtheorem{theorem}{Theorem}[section]
\newtheorem{proposition}[theorem]{Proposition}
\newtheorem{conjecture}[theorem]{Conjecture}
\newtheorem{corollary}[theorem]{Corollary}
\newtheorem{lemma}[theorem]{Lemma}

\newtheorem{maintheorem}[theorem]{Main Theorem}

\theoremstyle{definition}

\newtheorem{remark}[theorem]{Remark}
\newtheorem{example}[theorem]{Example}
\newtheorem{definition}[theorem]{Definition}

\newtheorem{problem}[theorem]{Problem}

\renewcommand{\eqref}[1]{{\rm (\ref{#1})}}

\def\proof{\smallskip\noindent {\bf Proof.\ }}
\def\endproof{\hfill$\square$\medskip}

\def\ZZ{\mathbb{Z}}
\def\CC{\mathbb{C}}

\def\RR{\mathcal{R}}
\def\QQ{\mathbb{Q}}

\def\FF{\mathcal{F}}

\def\UU{\mathcal{U}}
\def\VV{\mathcal{V}}
\def\OO{\mathcal{O}}
\def\BB{\mathcal{B}}

\def\nn{\mathfrak{n}}

\def\gg{\mathfrak{g}}
\def\hh{\mathfrak{h}}
\def\mm{\mathfrak{m}}

\def\uu{\mathfrak{u}}

\addtocounter{section}{-1}

\begin{document}
\title{ Braided Symmetric and Exterior Algebras}

\author{Arkady Berenstein}
\address{Department of Mathematics, University of Oregon,
Eugene, OR 97403, USA} \email{arkadiy@math.uoregon.edu}

\author{Sebastian Zwicknagl}
\address{\noindent Department of Mathematics, University of California, Riverside, California 92521, USA}
\email{zwick@math.ucr.edu}

\date{March 31, 2005, revised in May 2006}

\thanks{Research supported in part
by NSF grants \#DMS-0102382 and \#DMS-0501103}

\maketitle

\tableofcontents

\section{Introduction}
\label{sec:introduction} The goal of the paper is to introduce and
study  symmetric  and exterior algebras in certain braided monoidal
categories such as the category $\OO$ for quantum groups. We relate our {\it
braided  symmetric algebras} and {\it braided exterior algebras}
with their classical counterparts.

The original motivation of this work comes from the following classical problem.

\begin{problem}
\label{prob:multiplicity} Let $V$ be a simple finite-dimensional
module over a complex semisimple Lie algebra $\gg$, find the
decomposition of symmetric and exterior powers of $V$ into simple
$\gg$-modules.

\end{problem}

This problem is open and, to some extent, is not settled even for
$\gg=sl_2(\CC)$. 
A modern approach to the problem would consist of
constructing a basis (which can be referred to as a crystal basis or
a canonical basis)  in the symmetric algebra $S(V)$ (or in the
exterior algebra $\Lambda(V)$) which is somehow compatible with the
irreducible submodules.  Typically, the crystal (or canonical) bases
require some kind of $q$-deformation of the involved modules and
algebras.  Following  Lusztig's original idea of \cite{L0}, in order
to implement this program, one should first $q$-deform $S(V)$ (or
$\Lambda(V)$) and then try to construct the canonical (or crystal)
basis  for $S_q(V)$ (or $\Lambda_q(V)$) using a ``bar''-involution
and some standard PBW-like basis.

Therefore, the following problem naturally arises.
 \begin{problem}
 \label{Pr:qsymmetric} Let $U_q(\gg)$ be the quantized enveloping
algebra of $\gg$. For each finite dimensional $U_q(\gg)$-module $V$
construct (in the the category of $U_q(\gg)$-module algebras) an
analogue $S_q(V)$ of the symmetric algebra of $V$  and an analogue
$\Lambda_q(V)$ of the exterior algebra of $V$.
 \end{problem}

Apparently, the first work in which {\bf Problem}
\ref{Pr:qsymmetric} was mentioned is the paper \cite{R-D}. Its main
result asserts that if $V$ is the $4$-dimensional simple
$U_q(sl_2(\CC))$-module, then the algebra $S_q(V)$ cannot be a flat
deformation of $S(V)$. Unfortunately, the results of \cite{R-D}
neither suggested a solution to the problem nor generated a
follow-up research.  For the case when $V$ is the {\it quantum
adjoint  module}  over $U_q(sl_n(\CC))$,  a version of $S_q(V)$ was
constructed by J.~Donin in \cite{Donin}. Another, yet unrelated
approach, based on Woronowicz's {\it quantum differential calculus}
associates to $V$ an algebra $\BB(V)$ (which is sometimes referred
to as ''quantum symmetric algebra'' or Nichols-Woronowicz algebra,
see e.g., \cite{AS,DurdOzi,M, Wor}). However, $\BB(V)$ is typically
much larger than the ordinary symmetric algebra $S(V)$ .

We propose the following  construction of $S_q(V)$. First, recall
that the category $\OO_f$ of finite-dimensional $U_q(\gg)$-modules
is braided monoidal, where the braiding $\RR_{U,V}:U\otimes V\to
V\otimes U$ is the permutation of factors composed with  the
universal $R$-matrix.

For each object $V$ of $\OO_f$ define the {\it quantum exterior
square} $\Lambda_q^2 V\subset V\otimes V$ to be the linear span of
all {\it negative eigenvectors} of  $\RR_{V,V}:V\otimes V\to
V\otimes V$, that is, the eigenvectors that correspond to {\em
negative} (i.e., having the form $-q^r$, $r\in \ZZ$) eigenvalues of
$\RR$. Clearly, $\Lambda_q^2 V$ is a well-defined flat deformation
of $\Lambda^2 V$. Note that this definition makes sense because the
square $(\RR_{V,V})^2$ is a diagonalizable linear map $V\otimes V\to
V\otimes V$ which  eigenvalues are powers of $q$.

We  define the {\it braided symmetric algebra} $S_q(V)$ of  $V$ as the
quotient of the tensor algebra $T(V)$ by the ideal generated by
$\Lambda_q^2 V$. By the very definition the correspondence $V\mapsto S_q(V)$ is a functor
from the category $\OO_f$ to the category of graded $U_q(\gg)$-module algebras.

We prove (see Theorem \ref{th:flat poisson closure}) that the
braided symmetric and exterior powers are ''less or equal'' than
their classical counterparts, which motivates the following
definition. We say that a module $V$ is {\em flat} if $S_q(V)$ is a
flat deformation of the ordinary symmetric algebra $S(V)$.
\begin{problem}
Classify all flat modules in $\OO_f$.
\end{problem}

In section \ref{subsect:flat modules} (Proposition \ref{pr:flat
symmetric cube criterion}) we provide a sufficient criterion for
flatness, and numerous examples convince us that this criterion is
also necessary. Among the examples of flat modules is the standard
$n$-dimensional $U_q(sl_n(\CC))$-module $V$, and its braided
symmetric algebra is isomorphic to the algebra of $q$-polynomials in
$n$ variables.   In this paper we completely classify flat simple
$U_q(sl_2(\CC))$-modules and prove that a simple $n$-dimensional
$U_q(sl_2(\CC))$-module  is flat if and only if $n=1,2,3$. The fact
that the ''adjoint'' (i.e., $3$-dimensional simple) module $V$ is
flat agrees with the results of \cite{a}.

A general expectation is that for each $\gg$ there is a finite set
of dominant weights $\lambda$ such that the irreducible
$U_q(\gg)$-module $V_\lambda$ is flat. Informally speaking, we
expect that each ``small enough'' irreducible $U_q(\gg)$-module  is
flat. A complete classification of flat modules has recently been obtained in \cite{ZW}.

Let us reiterate that our immediate motivation for introducing  $S_q(V)$ and $\Lambda_q(V)$
is to prepare the ground for a general notion of the {\it canonical
basis}. 
It turns out that for each flat $U_q$-module $V$ from the list compiled in \cite{ZW}, the
algebra $S_q(V)$ is isomorphic to the associated graded of the quantized enveloping algebra $U_q(\uu)$ (where $\uu$ is a certain nilpotent Lie algebra) and, therefore, $S_q(V)$ carries a canonical basis which can be
defined along the lines of \cite{L0}. In its turn, this defines a
crystal basis for the ordinary symmetric algebra $S(\overline V)$,
where $\overline V$ is the corresponding $\gg$-module and,
therefore, solves {\bf Problem} \ref{prob:multiplicity}. A
surprising application of this crystal basis (and, in fact, another
motivation of the project) is a construction of a {\it geometric
crystal} on $\overline V$ (see \cite{bk1,bk2,bk3,bk4}).

A more ambitious goal is to study braided symmetric algebras of
non-flat modules. Already for $\gg=sl_2(\CC)$ we obtained an
unexpected result for each simple module $V$: the braided symmetric
cube $S^3_q V$ is multiplicity-free as a $U_q(sl_2(\CC))$-module
(Theorem \ref{th:exterior powers and symmetric cube}).  Numerous
examples convince us that this phenomenon should take place for any
braided symmetric power $S^n_q  V_\ell$ (Conjecture
\ref{conj:multiplicity-free}). Our general expectation is that the
braided symmetric powers of  finite dimensional $U_q(\gg)$-modules
have much nicer decompositions into the irreducibles than their
classical counterparts and, therefore, the ''braided'' version of
the above  {\bf Problem} \ref{prob:multiplicity} is much easier to
solve.

In a similar manner we define the {\it braided exterior algebra}
$\Lambda_q(V)$, to be the quotient of the tensor algebra $T(V)$ by
the ideal  generated by the {\it braided symmetric square}
$S^{\,2}_q V$. Similarly to the ordinary symmetric and exterior
algebras, one has a  {\it quadratic duality} between the quadratic
algebras $S_q(V)$ and $\Lambda_q(V^*)$. Again, similarly to the
classical case,  for each finite-dimensional module $V$ its braided
exterior algebra is also finite dimensional.

The structure of braided exterior algebras for
$U_q(sl_2(\CC))$-modules is more transparent: we prove that
$\Lambda^4_q V =0$ for any simple  $U_q(sl_2(\CC))$-module $V$ and,
moreover, if  $n=\dim V$ is even, then  $\Lambda^3_q V =0$.

Needless to say that this paper is only a starting point in
exploring the new braided algebras, and we expect many more
interesting results in their study.

{\bf Acknowledgments}. We express our  gratitude to Joel Kamnitzer
and Alexander Polishchuk for very stimulating discussions  during
the work on the paper.

\section{Definitions and notation}
\label{sect:definitions}

We start with the definition of the quantized enveloping algebra
associated with a complex reductive Lie algebra $\gg$ (our standard
reference here will be \cite{brown-goodearl}). Let $\hh\subset \gg$
be the Cartan subalgebra, and let $A=(a_{ij})$ be the Cartan matrix
for $\gg$.
\begin{definition}
\label{def:realization}
A \emph{realization} of~$A$ is a triple
$(P, \Pi, (\cdot \,|\, \cdot ))$, where $P$ is a lattice (of full dimension) in $\hh^*$,
$\Pi = \{\alpha_1, \dots, \alpha_r\} \subset P$,
is a subset,  and $(\cdot \,|\, \cdot )$ is a $\ZZ$-valued inner product on $P$ satisfying the following conditions:
\begin{itemize}
\item $\Pi$ is  linearly independent.

\item $\frac{2(\lambda  \,|\, \alpha_i)}{(\alpha_i \,|\, \alpha_i)}\in \ZZ$ and $\frac{(\alpha_i \,|\,
\alpha_i)}{2}\in \ZZ_{\ge 0}$ for all $\lambda\in P$,  $i\in [1,r]$.

\item $\frac{2(\alpha_i \,|\, \alpha_j)}{(\alpha_i \,|\, \alpha_i)} = a_{ij}$ for all $i, j$.

\end{itemize}
\end{definition}

The {\it quantized enveloping algebra} $U$ is a $\CC(q)$-algebra generated
by the elements $E_i$ and $F_i$ for $i \in [1,r]$, and $K_\lambda $ for $\lambda \in P$,
subject to the following relations:
%(see e.g.,\cite{Lusztig} or \cite{brown-goodearl}):
$K_\lambda  K_\mu = K_{\lambda +\mu}, \,\, K_0 = 1$
%(K_\ell )^{-1}=K_{-\ell }$$
for $\lambda , \mu \in P$; $K_\lambda E_i =q^{(\alpha_i\,|\,\lambda)} E_i K_\lambda,
\,\, K_\lambda F_i =q^{-(\alpha_i\,|\,\lambda)} F_i K_\lambda
$
for $i \in [1,r]$ and $\lambda\in P$;
\begin{equation}
\label{eq:upper lower relations}
E_i,F_j-F_jE_i=\delta_{ij}\frac{K_{\alpha_i}- K_{-\alpha_i}}{q^{d_i}-q^{-d_i}}
\end{equation}
for $i,j \in [1,r]$, where  $d_i=\frac{(\alpha_i\,|\,\alpha_i)}{2}$;
%where $K_i=K_{\alpha_i}$ and
%$q_i=\frac{(\alpha_i\,\,|\,\,\alpha_i)}{2}$;
and the {\it quantum Serre relations}
\begin{equation}
\label{eq:quantum Serre relations}
\sum_{p=0}^{1-a_{ij}} (-1)^p %\binom{1-a_{ij}}{p}_{q^{d_i}}
E_i^{(1-a_{ij}-p)} E_j E_i^{(p)} = 0,~\sum_{p=0}^{1-a_{ij}} (-1)^p %\binom{1-a_{ij}}{p}_{q^{d_i}}
F_i^{(1-a_{ij}-p)} F_j F_i^{(p)} = 0
\end{equation}
for $i \neq j$, where
%$\binom{1-a_{ij}}{p}_{q^{d_i}}$ is given by
%\eqref{eq:centered-binomial-coeff}.
the notation $X_i^{(p)}$ stands for the \emph{divided power}
\begin{equation}
\label{eq:divided-power}
X_i^{(p)} = \frac{X^p}{(1)_i \cdots (p)_i}, \quad
(k)_i = \frac{q^{kd_i}-q^{-kd_i}}{q^{d_i}-q^{-d_i}} \ .
\end{equation}

%It will be convenient to use
%
%$$(p)_i!={(1)_i \cdots (p)_i}\ ,  \quad  \binom{n}{m}_{q^{d_i}}=\frac{(n)_i!}{m_i! (n-m)_i!}\ ,   $$

The algebra $U$ is a $q$-deformation of the universal enveloping algebra of
the reductive Lie algebra~$\gg$ associated to~$A$, so it is commonly
denoted by $U = U_q(\gg)$.
It has a natural structure of a bialgebra with the co-multiplication $\Delta:U\to U\otimes U$
and the counit homomorphism  $\varepsilon:U\to \QQ(q)$
given by
\begin{equation}
\label{eq:coproduct}
\Delta(E_i)=E_i\otimes 1+K_{\alpha_i}\otimes E_i, \,
\Delta(F_i)=F_i\otimes K_{-\alpha_i}+ 1\otimes F_i, \, \Delta(K_\lambda)=
K_\lambda\otimes K_\lambda \ ,
\end{equation}
\begin{equation}
\label{eq:counit}
\varepsilon(E_i)=\varepsilon(F_i)=0, \quad \varepsilon(K_\lambda)=1\ .
\end{equation}
In fact, $U$ is a Hopf algebra with the antipode anti-homomorphism
$S: U \to U$ given by $S(E_i) = -K_{-\alpha_i} E_i, \,\, S(F_i) = -F_i K_{\alpha_i}, \,\,
S(K_\lambda) = K_{-\lambda}$.

Let $U^-$ (resp.~$U^0$; $U^+$) be the $\QQ(q)$-subalgebra of~$U$ generated by
$F_1, \dots, F_r$ (resp. by~$K_\lambda \, (\lambda\in P)$; by $E_1, \dots, E_r$).
It is well-known that $U=U^-\cdot U^0\cdot U^+$ (more precisely,
the multiplication map induces an isomorphism $U^-\otimes U^0\otimes U^+ \to U$).

%The algebra $U$ is graded by the root lattice $Q$:
%\begin{equation}
%\label{eq:U-grading}
%\text{$U=\bigoplus_{\alpha\in Q} U_\alpha, \quad
%U_\alpha = \{u \in U: K_\lambda u K_{-\lambda} = q^{(\lambda\,\,|\,\,\alpha)}\cdot u$
%for $\lambda\in P\}$.}
%\end{equation}
%Thus, we have
%$${\rm deg} E_i = \alpha_i, \quad {\rm deg} F_i = - \alpha_i, \quad
%{\rm deg} K_\lambda = 0 \ .$$

%%%%%%%%%%%%%%%%%%%%%%%%%%%%%%%%%%%%%%%%%%%%%
%%%%%%%%%%%%%%%%%%%%%%%%%%%%%%%%%%%%%%%%%%%%

We will consider the full sub-category  $\OO_{f}$ of the category
$U_q(\gg)-Mod$. The objects of $\OO_{f}$ are finite-dimensional
$U_{q}(\gg)$-modules $V$ having a weight decomposition
$$V=\oplus_{\mu\in P} V(\mu)\ ,$$
where each $K_\lambda$ acts on each {\it weight space} $V(\mu)$ by
the multiplication with $q^{(\lambda\,|\,\mu)}$  (see e.g.,
\cite{brown-goodearl}[I.6.12]). The category $\OO_{f}$ is semisimple
and the irreducible objects $V_\lambda$ are  generated by highest
weight spaces $V_\lambda(\lambda)=\CC(q)\cdot v_\lambda$, where
$\lambda$ is a {\it dominant weight}, i.e., $\lambda$ belongs to
$P^+=\{\lambda\in P:(\lambda\,|\,\alpha_i)\ge 0~ \forall ~i\in
[1,r]\}$, the monoid of dominant weights.

 Let $R \in U_{q}(\gg)\widehat \otimes U_{q}(\gg)$ be the universal $R$-matrix. By definition,
\begin{equation}
\label{eq:Jordan}
R=R_0R_1=R_1R_0
\end{equation}
 where $R_0$ is ''the diagonal part'' of $R$, and  $R_1$ is unipotent, i.e., $R_1$ is a formal power series
\begin{equation}
\label{eq:R1}
R_1=1\otimes 1+(q-1)x_1+ (q-1)^2x_2+\cdots \ ,
\end{equation}
where all $x_k\in {U'}^-_k\otimes_{\CC[q,q^{-1}]} {U'}^+_k$, where
${U'}^-$ (resp. ${U'}^+$) is the integral form of $U^+$, i.e.,
${U'}^-$ is a $\CC[q,q^{-1}]$-subalgebra of $U_q(\gg)$ generated by
all $F_i$ (resp. by all $E_i$) and ${U'}^-_k$ (resp. ${U'}^+_k$) is
the $k$-th graded component under the grading $deg(F_i)=1$ (resp.
$deg(E_i)=1$).

By definition, for any $U,V$ in $\OO_f$ and any highest weights elements $v_\lambda\in U(\lambda)$,  $v_\mu\in V(\mu)$ we have
$R_0(v_\lambda\otimes v_\mu)=q^{(\lambda\,|\,\mu)}v_\lambda\otimes v_\mu$.

Let $R^{op}$ be the opposite element of $R$, i.e., $R^{op}=\tau (R)$, where  $\tau:U_{q}(\gg)\widehat \otimes U_{q}(\gg)$
is the permutation of factors. Clearly, $R^{op}=R_0R_1^{op}=R_1^{op}R_0$.

Following \cite[Section 3]{DR1},  define $D\in  U_{q}(\gg)\widehat \otimes U_{q}(\gg)$ by
\begin{equation}
\label{eq:D}
D:=R_0\sqrt{R_1^{op}R_1}=\sqrt{R_1^{op}R_1}R_0 \ .
\end{equation}

Clearly, $D$ is well-defined because $R_1^{op}R_1$ is also unipotent as well as its square root.  By definition,
$D^{2}=R^{op}R$, $D^{op}R=RD$.

Furthermore, define
\begin{equation}
\label{eq:hat R}
\widehat R:=RD^{-1}=(D^{op})^{-1}R=R_1\left(\sqrt{R_1^{op}R_1}\right)^{-1}
\end{equation}
It is easy to see that
\begin{equation}
\label{eq:unitary}
 \widehat R^{op}=\widehat R^{\,\,-1}
\end{equation}
According to \cite[Proposition 3.3]{DR1}, the pair $(U_{q}(\gg), \widehat R)$ is a {\it coboundary} Hopf algebra.

The braiding in the category $\OO_{f}$ is defined by $\RR_{U,V}:U\otimes V\to V\otimes U$, where
$$\RR_{U,V}(u\otimes v)=\tau R(u\otimes v)$$
for any $u\in U$, $v\in V$, where $\tau:U\otimes V\to V\otimes U$ is the ordinary permutation of factors.

Denote by  $C\in Z(\widehat{U_{q}(\gg)})$   the {\it quantum
Casimir} element which acts on any irreducible $U_q(\gg)$-module
$V_\lambda$ in $\OO_f$ by the scalar multiple
$q^{(\lambda\,|\,\lambda+2\rho)}$, where $2\rho$ is the sum of
positive roots.

The following fact is  well-known.
\begin{lemma}
\label{le:braiding casimir} One has $\RR^2=\Delta(C^{-1})\circ
(C\otimes C)$. In particular, for each $\lambda,\mu,\nu\in  P_+$ the
restriction of $\RR^2$ to the $\nu$-th isotypic component
$I^\nu_{\lambda,\mu}$ of the tensor product $V_\lambda\otimes V_\mu$
is scalar multiplication by
$q^{(\lambda\,|\,\lambda)+(\mu\,|\,\mu)-(\nu\,|\,\nu))+(2\rho\,|\,\lambda+\mu-\nu)}$.

\end{lemma}

This allows to define the diagonalizable $\CC(q)$-linear  map
$D_{U,V}:U\otimes V\to U\otimes V$ by $D_{U,V}(u\otimes v)=
D(u\otimes v)$ for any objects $U$ and $V$ of $\OO_f$. It is easy to
see that the operator $D_{V_\lambda,V_\mu}:V_\lambda\otimes V_\mu\to
V_\lambda\otimes V_\mu$ acts on the $\nu$-th isotypic component
$I^\nu_{\lambda,\mu}$ in $V_\lambda\otimes V_\mu$ by the scalar
multiplication with
$q^{\frac{1}{2}(\,(\lambda\,|\,\lambda)+(\mu\,|\,\mu)-(\nu\,|\,\nu)\,)+(\rho\,|\,\lambda+\mu-\nu)}$.

For any  $U$ and $V$ in $\OO_f$ define the {\it normalized braiding} $\sigma_{U,V}$ by
\begin{equation}
\label{eq:sigma}
\sigma_{U,V}(u\otimes v)=\tau \widehat R (u\otimes v) \ ,
\end{equation}

Therefore, we have by (\ref{eq:hat R}):
\begin{equation}
\label{eq:formula for sigma}
\sigma_{U,V}=D_{V,U}^{-1} \RR_{U,V}=\RR_{U,V}D_{U,V}^{-1} \ .
\end{equation}
%where $D_{U,V}:U\otimes V\to U\otimes V$ is an invertible linear map defined by
%$$D_{U,V}(u\otimes v)=D(u\otimes v) \ .$$

We will will sometimes write $\sigma_{U,V}$ in a more explicit way:
\begin{equation}
\label{eq:sigma root}
\sigma_{U,V}=\sqrt{\RR_{V,U}^{-1}\RR_{U,V}^{-1}} \RR_{U,V}=\RR_{U,V}\sqrt{\RR_{U,V}^{-1}\RR_{V,U}^{-1}}
\end{equation}

%Moreover, when it is clear what $U$ and $V$ are, we will write:
%\begin{equation}
%\label{eq:sigma root 12}
%\sigma_{12}=\sqrt{\RR_{21}^{-1}\RR_{12}^{-1}} \RR_{12}=\RR_{12}\sqrt{\RR_{12}^{-1}\RR_{21}^{-1}}
%\end{equation}

The following fact is an obvious corollary of (\ref{eq:unitary}).
\begin{lemma}
\label{le:symm comm constraint} $\sigma_{V,U}\circ
\sigma_{U,V}=id_{U\otimes V}$ for any  $U,V$ in $\OO_f$. That is,
$\sigma$ is a symmetric commutativity constraint.
\end{lemma}

We also have the following coboundary relation (even though we will not use it).

\begin{lemma} \cite[section 3]{DR1}
Let $A,B,C$ be objects of $\OO_{f})$. Then, the following diagram commutes:
\begin{equation}
\label{eq:comm square}
\begin{CD}
A\otimes B \otimes C@>\sigma_{12,3}>>C\otimes A \otimes B\\
@V\sigma_{1,23}VV@VV\sigma_{23}V\\
B\otimes C \otimes A@>\sigma_{12}>>C\otimes B \otimes A
\end{CD}
\end{equation}
where we abbreviated
$$\sigma_{12,3}:=\sigma_{A\otimes B,C}:(A\otimes B) \otimes C\to  C\otimes (A \otimes B),$$
$$\sigma_{1,23}:=\sigma_{A, B\otimes C}:A \otimes (B\otimes C) \to  (B\otimes C)\otimes A \ .$$

\end{lemma}

\begin{remark} If one replaces the braiding  $\RR$ of $\OO_f$  by it inverse $\RR^{-1}$,
the symmetric commutativity constraint $\sigma$ will not change.

\end{remark}
\section{Main results}

\label{sec:main results}

\subsection{Braided symmetric and exterior powers}
In this section we will use the notation and conventions of Section \ref{sect:definitions}.

For any morphism $f:V\otimes V\to V\otimes V$ in $\OO_f$ and $n>1$
we denote by $f^{i,i+1}$, $i=1,2,\ldots,n-1$ the morphism
$V^{\otimes n}\to V^{\otimes n}$ which acts as $f$ on the $i$-th and
the $i+1$st factors. Note that $\sigma_{V,V}^{i,i+1}$ is always an
involution on $V^{\otimes n}$.

\begin{definition}
\label{def:symmetric power} For an object $V$ in $\OO_f$ and $n\ge
0$ define the {\it braided symmetric power} $S_\sigma^nV \subset
V^{\otimes n}$ and the {\it braided exterior power}
$\Lambda_\sigma^nV\subset V^{\otimes n}$  by:
$$S_\sigma^nV=\bigcap_{1\le i\le n-1}  (Ker~ \sigma_{i,i+1}-id)=\bigcap_{1\le i\le n-1} (Im~\sigma_{i,i+1}+id)\ ,$$
$$\Lambda_\sigma^nV=\bigcap_{1\le i\le n-1} (Ker~ \sigma_{i,i+1}+id)=\bigcap_{1\le i\le n-1} (Im~\sigma_{i,i+1}-id) ,$$
where we abbreviated $\sigma_{i,i+1}=\sigma_{V,V}^{i,i+1}$.
\end{definition}

\begin{remark} Clearly, $-\RR$ is also a braiding on $\OO_f$ and $-\sigma$ is the corresponding normalized braiding.
Therefore, $\Lambda_\sigma^nV=S_{-\sigma}^nV$ and
$S_\sigma^nV=\Lambda_{-\sigma}^nV$. That is, informally speaking,
the symmetric and exterior powers are mutually ''interchangeable''.

\end{remark}

\begin{remark} Another way to introduce the symmetric and exterior squares involves the well-known fact that the braiding $\RR_{V,V}$
is a semisimple operator $V\otimes V\to V\otimes V$, and all the
eigenvalues of $\RR_{V,V}$ are of the form $\pm q^r$, where $r\in
\ZZ$. Then {\it positive} eigenvectors of $\RR_{V,V}$ span
$S_\sigma^2V$ and {\it negative} eigenvectors  of $\RR_{V,V}$ span
$\Lambda_\sigma^2V$.

\end{remark}

Clearly, $S_\sigma^{0}V=\CC(q)\ ,\ S_\sigma^{1}V=V\ , \Lambda_\sigma^{0}V=\CC(q)\ ,\ \Lambda_\sigma^{1}V=V$, and
$$S_\sigma^2V=\{v\in V\otimes V\,|\, \sigma_{V,V}(v)=v\},~\Lambda_\sigma^2V=\{v\in V\otimes V\,|\, \sigma_{V,V}(v)=-v\} \ .$$

We can provide a uniform characterization of braided symmetric and exterior powers as follows.

\begin{definition}
\label{def:braided intersections} Let $V$ be an object of $\OO_f$
and $I$ be a sub-object of $V\otimes V$ in $\OO_f$. For each $n\ge
2$ define the {\it braided power} $P(I)^n\subset V^{\otimes n}$ of
$I$ by:
\begin{equation}
\label{eq:braided intersections}
P^n(I)=\bigcap_{1\le i\le n-1} V^{\otimes ~i-1}\otimes  I\otimes V^{\otimes~n-1-i} \ .
\end{equation}
(and set $P^0(I)=\CC(q)$, $P^1(I)=V$).
\end{definition}
The following fact is obvious.

\begin{lemma} For each $n\ge 0$ we have $S_\sigma^nV= P^n(S_\sigma^2 V)$ and $\Lambda_\sigma^nV= P^n(\Lambda_\sigma^2)$.

\end{lemma}

The following fact is obvious.

\begin{proposition}
\label{pr:injections} For each $n\ge 0$ the association $V\mapsto
S^n_\sigma V$ is a functor from $\OO_f$ to $\OO_f$ and the
association $V\mapsto \Lambda_\sigma^n V$ is a functor from $\OO_f$
to $\OO_f$. In particular, an embedding $U\hookrightarrow V$ in the
category $\OO_f$ induces injective morphisms
$$S_\sigma^n U \hookrightarrow S_\sigma^n V, ~ \Lambda_\sigma^n U\hookrightarrow \Lambda_\sigma^n V \ .$$

\end{proposition}

\begin{definition}
\label{def:symmetric algebra}
For any  $V\in Ob(\OO)$   define the {\it braided symmetric algebra} $S_\sigma(V)$ and the {\it braided  exterior algebra} $\Lambda_\sigma(V)$ by:
\begin{equation}
\label{eq:injections}
S_\sigma(V)=T(V)/\left< \Lambda_\sigma^{2}V\right>,~\Lambda_\sigma(V)=T(V)/ \left<S_\sigma^{2}V\right> \ ,
\end{equation}
where $T(V)$ is the tensor algebra of $V$ and $\left<I\right>$ stands for the two-sided  ideal in $T(V)$ generated by a subset $I\subset T(V)$.
\end{definition}

Note that the algebras $S_\sigma(V)$ and $\Lambda_\sigma(V)$ carry a natural $\ZZ_{\ge 0}$-grading:
$$S_\sigma(V)=\bigoplus_{n\ge 0} S_\sigma(V)_n, ~~~\Lambda_\sigma(V)=\bigoplus_{n\ge 0} \Lambda_\sigma(V)_n\ ,$$
since the respective ideals in $T(V)$ are homogeneous.

\smallskip

Denote by  $\OO_{gr,f}$ the sub-category of $U_q(\gg)-Mod$ whose objects are $\ZZ_{\ge 0}$-graded:
$$V=\bigoplus_{n\in \ZZ_{\ge 0}}  V_n$$
where each $V_n$ is an object of $\OO_f$; and morphisms are those homomorphisms of  $U_q(\gg)$-modules which preserve the $\ZZ_{\ge 0}$-grading.

Clearly,  $\OO_{gr,f}$ is a tensor category under the natural
extension of the tensor structure of $\OO_f$. Therefore, we can
speak of algebras and coalgebras in $\OO_{gr,f}$.

%Clearly, $\OO_{gr,f}$ is a braided monoidal category because the {\it braided} tensor product $A\otimes _\RR B$ %of two algebras $A,B\in Ob(\OO_{gr,f})$ is well-defined:
%$$(a\otimes b)(a'\otimes b')=a\RR_{B,A}(b\otimes a')b'$$
%for any $a,a'\in A$ and $b,b'\in B$ (and the braiding $\RR_{A,B}:A\otimes _\RR B\to B\otimes %_\RR A$ is inherited from $\OO_f$).
%
%
%\begin{proposition} ???? For any $U,V\in Ob(\OO_f)$ there exist surjective  morphisms in $\OO_{gr,f}$:
%\begin{equation}
%\label{eq:sub-exponential functors}
%\varphi_{U,V}:S_\sigma(U)\otimes _\RR S_\sigma(V)\twoheadrightarrow S_\sigma(U\oplus V), %~\psi_{U,V}:\Lambda_\sigma(U)\otimes _{-\RR} \Lambda_\sigma(V)\twoheadrightarrow \Lambda_\sigma(U\oplus V)\ .
%\end{equation}
%That is, the functors $S_\sigma$ and $\Lambda_\sigma$ are sub-exponential.
%
%\end{proposition}

By the very definition,  $S_\sigma(V)$ and $\Lambda_\sigma(V)$ are algebras in $\OO_{gr,f}$.

\begin{proposition} The assignments $V\mapsto S_\sigma(V)$ and $V\mapsto \Lambda_\sigma(V)$ define functors from $\OO_f$
to the category of algebras in $\OO_{gr,f}$.
\end{proposition}

Let $V$ be an object of $\OO_f$ and $I^2$ be a sub-object of
$V\otimes V$ in $\OO_f$ as in Definition \ref{def:braided
intersections}. Using the identity $P^m(I) \otimes  V^{\otimes n}
\cap   V^{\otimes m}\otimes P^n(I)=P^m(I)\otimes P^n(I)$ and natural
inclusions $P^{m+n}(I)\subset P^m(I) \otimes  V^{\otimes n}, ~
P^{m+n}(I)\subset V^{\otimes m}\otimes P^n(I)$, we obtain natural
embeddings $\Delta_{m,n}^I:P^{m+n}(I)\hookrightarrow P^m(I) \otimes
P^n(I)$ in the category $\OO_f$.

Let $P(I):=\bigoplus\limits_{n\ge 0} P^n(I)$; define the morphism
$\Delta_I:P(I)\to P(I)\otimes P(I)$ by
$\Delta^I=\bigoplus\limits_{n,m\ge 0} \Delta_{m,n}^I$ and let
$\varepsilon_I$ be the natural projection $P(I)\to P^0(I)= \CC(q)$.

\begin{lemma} The triple $(P(I),\Delta^I,\varepsilon_I)$ is a coalgebra in $\OO_{gr,f}$.
\end{lemma}

\proof Clearly, $\Delta_I$ is co-associative. The compatibility between the co-multiplication and the co-unit also follows.
\endproof

\begin{corollary} For each object $V$ of $\OO_f$ the braided symmetric space $S_\sigma V=P(S^2_\sigma V)$
and the braided exterior space $\Lambda_\sigma V=P(\Lambda^2_\sigma
V)$ are naturally co-algebras in $\OO_f$. The assignments $V\mapsto
S_\sigma V$ and $V\mapsto \Lambda_\sigma V$ define functors from
$\OO_f$ to the category of coalgebras in $\OO_{gr,f}$.

\end{corollary}

Note that both $\OO_f$ and $\OO_{gr,f}$ are categories with the duality ${}^*$.

\begin{proposition}
\label{pr:algebra coalgebra duality}
For any $V$ in $\OO_f$ we have the following isomorphisms in $\OO_{gr,f}$:

\noindent (a) Algebra isomorphisms: $(S_\sigma V)^*\cong S_\sigma(V^*)$, $(\Lambda_\sigma V)^*\cong \Lambda_\sigma(V^*)$.

\noindent (b) Coalgebra isomorphisms: $(S_\sigma (V))^*\cong S_\sigma V^*$, $(\Lambda_\sigma (V))^*\cong \Lambda_\sigma V^*$.

\noindent (c) Algebra isomorphisms: $\Lambda_\sigma(V)^{!}\cong
S_\sigma(V^*)$,  $S_\sigma(V)^{!} \cong \Lambda_\sigma(V^*)$, where
$A^!$ stands for the quadratic dual of a quadratic algebra $A$ in
$\OO_{gr, f}$.

\end{proposition}

\proof Follows from the fact that $T(V)^*\cong T(V^*)$ in
$\OO_{gr,f}$, $(\Lambda_\sigma^2 V)^\perp=  S_\sigma^2 V^*$,
$(S_\sigma^2 V)^\perp= \Lambda_\sigma^2 V^*$ in $\OO_f$ and the
general fact that for any $I,J\subset V\otimes V$ such that
$I+J=V\otimes V$ and $I\cap J=0$ one has:
$$P(I)^\perp=\langle J^\perp\rangle,~ P(J)^\perp=\langle I^\perp\rangle \, $$
in $\OO_{gr,f}$, where  $\langle X\rangle$ is the ideal in $T(V)^*$ generated by $X\subset (V\otimes V)^*=V^*\otimes V^*$.
%To prove part (d) fix  $I,J\subset V\otimes V$ such that $V\otimes V=I\oplus J$, $n\ge 1$ and abbreviate $U:=V^\otimes n$, $A_i:=V^{\otimes i-1} %\otimes I\otimes V^{\otimes n-i}$, $B_i:=V^{\otimes i-1} \otimes J\otimes V^{\otimes n-i}$ for $i=1,\ldots,n$. In this notation $P^{n+1}(I)=\cap_i %A_i$, the $n+1$-st component $\langle I\rangle_{n+1}$ of the ideal  $\langle I\rangle$ is given by $\langle I\rangle_{n+1}=\sum_i B_i$.
%
%
%in $T(V)$ generated by $I$ satisfies $\langle I\rangle^*=\langle I^*\rangle$. Taking into account that $T(V)^*=T(V^*)$,  and applying $*$ to the %isomorphism $T(V)\cong \langle I\rangle \oplus  T(V)/ \langle I\rangle $ in $\OO_{gr, f}$, we obtain an isomorphism $T(V^*)\cong  \langle I\rangle^* %\oplus (T(V)/\langle I\rangle)^* $ and $(T(V)/\langle I\rangle)^*\cong T(V^*)/ \langle I\rangle^*$. Finally, using the isomorphisms  $(\Lambda_\sigma^2 %V)^*= \Lambda_\sigma^2 V^*$,  $(S_\sigma^2 V)^*= S_\sigma^2 V^*$ in $\OO_f$, we finish the proof of (d). The proposition is proved.
\endproof

In particular, Proposition \ref{pr:algebra coalgebra duality}(c) asserts that the relationship between the braided symmetric (resp. exterior) powers and the corresponding components of braided symmetric (resp. exterior) algebras is given by the following canonical isomorphisms (for all $n\ge 2$).
\begin{equation}
\label{eq:double duals}
(S_\sigma^n V^*)^*\cong S_\sigma(V)_n, ~(\Lambda_\sigma^n V^*)^*\cong \Lambda_\sigma(V)_n  \ .
\end{equation} 

To discuss a possibility of a direct isomorphism, we propose the following problem.

\begin{problem}
\label{pr:algebra coalgebra isomorphism} Given $V$ in $\OO_f$, find
all the pairs of sub-objects $I,J\subset V\otimes V$ such that
$I+J=V\otimes V$, $I\cap J=0$, and the composition of the natural
inclusion $P(I) \hookrightarrow T(V)$ with the structure
homomorphism $T(V)\twoheadrightarrow T(V)/\langle J \rangle$ is an
isomorphism  $P(I)\cong T(V)/\langle J \rangle$ in $\OO_{gr, f}$.

\end{problem}

We would expect that  the pairs $(I,J)=(S^2_\sigma
V,\Lambda^2_\sigma V)$ and $(I,J)=(\Lambda^2_\sigma V,S^2_\sigma V)$
are solutions to {\bf Problem} \ref{pr:algebra coalgebra
isomorphism} for any $V\in \OO_f$.  This expectation is certainly
true whenever  when $\gg$ is abelian, i.e., when $\sigma=\tau$, and
$S_\sigma(V)=S(V)$, $\Lambda_\sigma(V)=\Lambda(V)$.

\begin{proposition}
\label{pr:n-th symm-power} For any  $V$ in $\OO_{f}$ each embedding
$V_{\lambda}\hookrightarrow V$ defines embeddings
$V_{n\lambda}\hookrightarrow S_\sigma^n V $ for all  $n\ge 2$. In
particular, the  algebra $S_\sigma(V)$ is infinite-dimensional.
\end{proposition}

\proof By definition, an embedding $V_{\lambda}\hookrightarrow V$ is
equivalent to the existence of a highest weight vector $v_\lambda$
in the weight space $V(\lambda)$. Then $V_\lambda\cong
U_q(\gg)(v_\lambda)$.

Denote  by $v_{n\lambda}=v_{\lambda}\otimes v_{\lambda}\otimes  \cdots \otimes v_{\lambda}\in V^{\otimes n}$.
First, show that   $v_{2\lambda}\in S_\sigma^2 V$. Indeed,
$$\RR_{V\otimes V}(v_{2\lambda})=\RR_{V\otimes V}(v_{\lambda}\otimes v_{\lambda})=
q^{(\lambda|\lambda)}v_{\lambda}\otimes v_{\lambda}=D(v_{\lambda}\otimes v_{\lambda})\ . $$
Therefore, $\sigma_{V\otimes V}(v_{2\lambda})=v_{2\lambda}$ and $v_{2\lambda}\in  S_\sigma^2 V$.
This implies that $\sigma_{i,i+1}(v_{n\lambda})=v_{n\lambda}$
for $i=1,\ldots,n-1$.
Hence $v_{n\lambda}\in S_\sigma^n V$ and $V_{n\lambda}\hookrightarrow  S_\sigma^n V$.
Proposition \ref{pr:n-th symm-power} is proved.
\endproof

\subsection{Braided algebras as deformations of Poisson structures}
\label{subsect:semiclassical limit}

In this section we study ''classical'' (at $q=1$) analogues of the
braided symmetric and exterior algebras and present the braided
objects as $q$-deformations of the classical ones. More precisely, using the fact that the ''classical limit'' of each braided algebra is a  commutative $\CC$-(super)algebra, we will investigate the (super-)Poisson structures emerging as ''classical limits'' of the (super-)commutators in the braided algebras (Theorem \ref{th:flat poisson closure} below).

\begin{definition}
\label{def:commutative superalgebra} Given a $\ZZ_2$-graded algebra
$A=A_{\bar 0}\oplus A_{\bar 1}$, (i.e., $A_{\varepsilon}\cdot
A_{\delta}\subset A_{\varepsilon+ \delta}$ for any
$\varepsilon,\delta\in \ZZ_2=\ZZ/(2)=\{\overline 0,\overline 1\}$),
we say that $A$ is a {\it commutative   superalgebra}  if
$$ba=(-1)^{\varepsilon\delta}ab$$ for any
$a\in A_\varepsilon$ and $b\in A_\delta$, $\varepsilon,\delta\in \ZZ_2=\{\bar 0,\bar 1\}$.

\end{definition}

\begin{definition}
\label{def:bracketed superalgebra}
 A {\it bracketed superalgebra} is a pair $(A,\{\cdot,\cdot \})$ where $A$ is a commutative  superalgebra  and a  bilinear map
$\{\cdot,\cdot \}:A\times A\to A$
such that $\{A_\varepsilon,A_\delta\}\subset A_{\varepsilon+ \delta}$ for any $\varepsilon,\delta\in \ZZ_2$ and the following identities hold:

\noindent (i) {\it super-anti-commutativity}:
$$\{ a,b\}+(-1)^{\varepsilon\delta}\{b,a\}=0$$
for any $a\in A_\varepsilon, b\in A_\delta$,

\noindent (ii) the {\it super-Leibniz rule}
\begin{equation}
\label{eq: super Leibniz rule}
\{a,bc \}= \{a,b \}c+(-1)^{\varepsilon\delta}b\{a,c\}
\end{equation}
for any $a\in A_\varepsilon, b\in A_\delta, c\in A_\gamma$.

If,  in addition, the bracket satisfies

\noindent (iii)  the {\it super-Jacobi identity}:
$$(-1)^{\varepsilon\gamma} \{a,\{b,c\}\}+ (-1)^{\gamma\delta} \{c,\{a,b\}\}+(-1)^{\delta\varepsilon} \{b,\{c,a\}\}=0$$
for any $a\in A_\varepsilon, b\in A_\delta, c\in A_\gamma$,
then we will refer to $A$ as a {\it Poisson superalgebra} and will refer to $\{\cdot,\cdot \}$ as {\it  super-Poisson bracket}.

\end{definition}

We define the tensor product $A\otimes B$ of two  bracketed
superalgebras $A$ and $B$ to be the usual tensor product in the
category of superalgebras with bracket $\{\cdot,\cdot\}$:

$$ \{a\otimes b, a'\otimes b'\}=(-1)^{\varepsilon'\delta}aa'\otimes \{b,b'\}+(-1)^{\delta\varepsilon'+\delta\delta'}\{a,a'\} \otimes bb', $$
for all $a\in A, b'\in B, a'\in A_{\varepsilon'}, b\in B_{\delta}$.

%???????Our two main examples of homogeneous bracketed superalgebra are $S_{\varphi}(V)$ and $\Lambda_\phi(V)$ where $U$ and $V$ are $\bf %k$-vectorspaces and we consider $S(U)$ to be generated in the degree $2$; i.e.,
%$$S_{\varphi_+}(U)=\oplus_{i\in \ZZ_{+}} S_{\varphi_+}(U)_{2i}\ ,~  S_{\varphi_+}(U)_{2i}=S^i_{\phi_+}(U)\ .$$

For each reductive    Lie algebra $\gg$ denote by $\overline
\OO_f=\overline \OO_f(\gg)$ the category of finite-dimensional
$\gg$-modules $\overline V$ such that the center $Z(\gg)$ acts
semi-simply on $\overline V$.

Recall that the classical $r$-matrix of $\gg$ is a $\gg$-invariant
element $r\in \gg\otimes \gg$ satisfying the classical Yang-Baxter
equation $[[r,r]]=0$ and $r+\tau (r)=c$, where $\tau$ is the
permutation of factors and $c\in S^2(\gg)$ is the Casimir element.
Denote by $r^-\in\gg \wedge \gg$ the anti-symmetrization of $r$,
i.e.,
\begin{equation}
\label{eq:rminus}
r^-=\frac{1}{2}(r-\tau(r)) \ .
\end{equation}
Note that $r^-=\sum_{\alpha} E_{\alpha} \otimes F_\alpha-F_\alpha\otimes F_\alpha $,
where the summation is over all positive roots of $\gg=\nn_-\oplus \hh\oplus \nn_+$,  $\{E_\alpha\}$ is a basis for $\nn_+$, and $\{F_\alpha\}$ is the basis in $\nn_+$ dual to the former one with  respect to the Killing form.

\begin{example}
\label{ex:gl2 rmatrix} Recall that $gl_2(\CC)$ has a standard basis of
matrix units $E,F,H_1,H_2$ and $r=2E\otimes F+H_1\otimes
H_1+H_2\otimes H_2$. Therefore, $r^-=\frac{1}{2}(r-\tau(r))=E\otimes
F-F\otimes E=E\wedge F$.

\end{example}

The following fact is obvious.

\begin{lemma}
\label{le:superbrackets}
 (a) Let $\overline V$ be any finite-dimensional $\gg$-module. Then the element $r^-$ considered as $\gg$-equivariant operator
 $r^-:\overline V\otimes \overline V\to \overline V\otimes \overline V$ satisfies
$$r^-(\Lambda^2(\overline V))\subset S^2(\overline V), ~r^-(S^2(\overline V))\subset \Lambda^2(\overline V) \ .$$

\noindent (b) The symmetric algebra $S(\overline V)$ considered as
an even superalgebra (i.e., $S(\overline V)_{\overline
0}=S(\overline V)$) has a unique homogeneous bracket $\{\cdot
,\cdot\}_+:S(\overline V)\times S(\overline V)\to S(\overline V)$
such that
$$ \{u ,v\}_+ =r^-(u\cdot v)$$
for any $u,v\in \overline V$.

\noindent  (c) The symmetric algebra $\Lambda(\overline V)$
considered as a superalgebra with odd $\overline V$ has a unique
homogeneous super-bracket $\{\cdot ,\cdot\}_-:\Lambda(\overline
V)\times \Lambda(\overline V)\to \Lambda(\overline V)$ such that
$$ \{u ,v\}_- =r^-(u\wedge v)$$
for any $u,v\in  \overline V$.

\end{lemma}

For any bracketed superalgebra define the {\it super-Jacobian} map $J:A\times A\times A\to A$ by:
\begin{equation}
\label{eq:Jacobian}
J(a,b,c)=(-1)^{\varepsilon\gamma} \{a,\{b,c\}\}+ (-1)^{\gamma\delta} \{c,\{a,b\}\}+(-1)^{\delta\varepsilon} \{b,\{c,a\}\}
\end{equation}
for any $a\in A_\varepsilon, b\in A_\delta, c\in A_\gamma$.

Define the ideals $I_+(\overline V)$ and $I_-( \overline V)$ in $S(
\overline V)$ and $\Lambda( \overline V)$ respectively by:
$$I_+(\overline V):=J_+(\overline V,\overline V,\overline V)\cdot S(\overline V),
~I_-(\overline V):=J_-(\overline V,\overline V,\overline V)\wedge \Lambda(\overline V) \ ,$$
where $J_+$ and $J_-$ are respective Jacobian maps as in (\ref{eq:Jacobian}). Then define  algebras
\begin{equation}
\label{eq:barSV barLambdaV}
\overline {S(\overline V)}:=S(\overline V)/I_+(\overline V),~\overline {\Lambda (\overline V)}:=\Lambda (\overline V)/I_-(\overline V)\ .
\end{equation}

\begin{theorem}
\label{th:symmetric exterior almost poisson} For each $\overline V$
in $\overline \OO_f$ we have (using the Sweedler-like notation
$r^-=r^-_{(1)}\otimes r^-_{(2)}\ $):

\noindent (a)  The brackets $\{\cdot,\cdot\}_+$ and $\{\cdot,\cdot\}_-$ are  given by:
\begin{equation}
\label{eq:almost poisson symmetric1}
\{a,b\}_+=r^-_{(1)}(a)\cdot r^-_{(2)}(b)
\end{equation}
for each $a,b\in S(\overline V)$.
\begin{equation}
\label{eq:almost poisson exterior1}
\{a,b\}_-=r^-_{(1)}(a)\wedge r^-_{(2)}(b)
\end{equation}
for each $a,b\in \Lambda(\overline V)$.

%In particular, $\{S(\overline V),I_+(\overline V)\}_+\subset I_+(\overline V)$ and $\{\Lambda(\overline V),I_-(\overline V)\}_-\subset I_-(\overline V)$.

\noindent (b) The homogeneous bracketed (super)algebras $\overline
{S(\overline V)}$ and $\overline {\Lambda(\overline V)}$ are Poisson
with the (super)Poisson brackets given by the above formulas
(\ref{eq:almost poisson symmetric1}) and (\ref{eq:almost poisson
exterior1}).
\end{theorem}

% \begin{lemma} $[[r^-,r^-]]=\frac{1}{2}[c_{13},r_{21}-r_{32}]$, where $c=\frac{1}{2}(r+\tau(r))$ is the Casimir element of $\gg$.
%\end{lemma}
%
%\proof Recall that
%$$[[\Phi,\Phi]]=[\Phi_{12},\Phi_{13}]+[\Phi_{12},\Phi_{23}]+[\Phi_{13},\Phi_{23}]$$
%for any $\Phi\in U(\gg)\otimes U(\gg)$. Conjugating with $\tau$, we obtain several equations:
%$$[[r,r]]=[r_{12},r_{13}]+[r_{12},r_{23}]+[r_{13},r_{23}]$$
%$$\tau_{12}([[r,r]])=[r_{21},r_{23}]+[r_{21},r_{13}]+[r_{23},r_{13}]$$
%$$\tau_{23}\tau_{12}([[r,r]]))=[r_{31},r_{32}]+[r_{31},r_{12}]+[r_{32},r_{12}]$$
%$$\tau_{13}([[r,r]])=[r_{32},r_{31}]+[r_{32},r_{21}]+[r_{31},r_{21}]$$
%$$\tau_{12}\tau_{23}([[r,r]])=[r_{23},r_{21}]+[r_{23},r_{31}]+[r_{21},r_{31}]$$
%Let us compute:
%$$4[[r^-,r^-]]=[r_{12}-r_{21},r_{13}-r_{31}]+[r_{12}-r_{21},r_{23}-r_{32}]+[r_{13}-r_{31},r_{23}-r_{32}]$$
%$$=(1-\tau_{13})[[r,r]]+[r_{13},r_{21}]+[r_{31},r_{12}]+[r_{23},r_{21}]+[r_{32},r_{12}]+[r_{23},r_{31}]+[r_{32},r_{13}]$$
%$$=(1-\tau_{13})[[r,r]]+[r_{13},r_{21}]+([r_{31},r_{12}]+[r_{32},r_{12}])+([r_{23},r_{21}]+[r_{23},r_{31}])+[r_{32},r_{13}]$$
%$$=(1-\tau_{13}+\tau_{23}\tau_{12}+\tau_{12}\tau_{23})[[r,r]]+[r_{13},r_{21}]-[r_{31},r_{32}]-[r_{21},r_{31}]+[r_{32},r_{13}]$$
%$$=(1-\tau_{13}+\tau_{23}\tau_{12}+\tau_{12}\tau_{23})[[r,r]]+[c_{13},r_{21}]+[r_{32},c_{13}]$$
%Since, $r$ satisfies the classical Yang-Baxter equation $[[r,r]]=0$, we obtain the desirable result.
%
%The lemma is proved. \endproof

%%%%%%%%%%%%%%%%%%

We will prove Theorem \ref{th:symmetric exterior almost poisson} in Section \ref{sect:proofs of main results}.
\begin{example}
\label{ex:r-gl2} Recall from Example \ref{ex:gl2 rmatrix} that for
$\gg=gl_2(\CC)$ we have $r^-=E\wedge F=E\otimes F-F\otimes E$. Theorem
\ref{th:symmetric exterior almost poisson} guarantees that for any
$\overline V$ in $\overline \OO_f$ one has:
$$\{a,b\}_+=E(a) F(b)-F(a)E(b)$$
for all $a,b\in S(\overline V)$ and
$$\{a,b\}_-=E(a)\wedge F(b)-F(a)\wedge E(b)$$
for all $a,b\in\Lambda(\overline V)$.

\end{example}

\begin{corollary} (from the proof  of Theorem \ref{th:symmetric exterior almost poisson})
\label{cor:poisson closures} The Poisson (super)algebras
$\overline{S(\overline V)}=S(\overline V)/I_+(\overline V)$ and
$\overline{\Lambda(\overline V)}=\Lambda (\overline V)/I_-(\overline
V)$  are  Poisson closures  of $S(\overline V)$ and $\Lambda
(\overline V)$ respectively.

\end{corollary}

The following is the second main result of this section.

\begin{maintheorem}
\label{th:flat poisson closure} Let $V$ be an object of $\OO_f$ and
let $\overline V$ in $\overline \OO_f$ be the classical limit of
$V$. Then:

\noindent (a) The braided symmetric algebra $S_\sigma(V)$ is a flat
$q$-deformation of a certain Poisson quotient of the Poisson algebra
$\overline{S(\overline V)}$. In particular, $\dim_{\CC(q)}
S_\sigma(V)_n\le \dim_{\CC} \overline{S(\overline V)}_n$.

\noindent (b) The braided exterior algebra $\Lambda_\sigma(V)$ is a
flat $q$-deformation of a certain Poisson quotient of the Poisson
superalgebra $\overline{\Lambda(\overline V)}$. E.g., $\dim_{\CC(q)}
\Lambda_\sigma(V)_n\le \dim_{\CC} \overline{\Lambda(\overline
V)}_n$.

\end{maintheorem}

We will prove Theorem \ref{th:flat poisson closure} in Section \ref{sect:proofs of main results}.

Now we can post two natural problems.

\begin{problem} Describe all those $U_q(\gg)$-modules in $\overline V$ for which $S_\sigma(V)$ is  a flat deformation of $\overline{S(\overline V)}$.

\end{problem}
Even though we have only a  little evidence,  we would expect that   each irreducible $V=V_\lambda$ in $\OO_f$ solves the problem.

Since $\overline {S_\sigma (\overline V)}$ is a quotient algebra of
the symmetric algebra $S(\overline V)$, we can define an affine
scheme $X(\overline V)=spec\,\, \overline {S_\sigma (\overline V)}$
in $\overline V^{\,*}$.

\begin{problem}
\label{pr:schemes} For each  $V\in Ob(\OO_f)$ study the scheme
$X(\overline V)$. In particular, describe all objects $V$ for which
the scheme $X(\overline V)$ is a variety.
\end{problem}

Clearly, if $V$  is flat (see Section \ref{subsect:flat modules}),
then $X(\overline V)=\overline V^{\,*}$ is a variety. It would be
interesting to see what happens with $X(\overline V)\ne \overline V$, i.e., $S_\sigma(V)$ is not a flat deformation of $S(\overline V)$ (see Sections \ref{subsect:flat modules} and \ref{subsect:gl_2} below).

\medskip

To study the scheme $X(\overline V)$ it would be helpful to  estimate the growth of braided symmetric algebras. We will conclude the section with a few corollaries from Theorem \ref{th:flat poisson closure} and one conjecture concerning the Hilbert series of braided algebras.

%%%%%%%%%%%%%%%%%%%%%%%%%%%%%%%%%%%%%

\begin{definition}
\label{def:Hilbert series}
For each $U=\oplus_{n\ge 0} U_n$ in the category $\OO_{gr,f}$ define the {\it Hilbert series} $h(U,t)\in \ZZ[[t]]$  by:
\begin{equation}
\label{eq:Hilbert series}
h(U,t)=\sum_{i=0}^{\infty} \dim(U_i)t^i.
\end{equation}
\end{definition}

According to  Theorem \ref{th:flat poisson closure}(b),
$\Lambda_\sigma(V)$ is bounded from above (as a vector space) by the
ordinary exterior algebra $\Lambda(V)$. Therefore, the Hilbert
series $h(\Lambda_\sigma(V),t)$ is a polynomial in $t$ of degree
at most $\dim(V)$.

\begin{corollary}
\label{cor:rationality}
For any $V\in Ob(\OO_f)$ the Hilbert series
$h(S_\sigma(V)),t)$ is a rational function in $t$ of the form $\frac{p(t)}{(1-t)^d}$, where $p(t)$ is a polynomial and $d\le \dim V$.
\end{corollary}

\proof  Theorem \ref{th:flat poisson closure} asserts that
$S_\sigma(V)$ is a flat deformation of a certain quotient algebra of
the symmetric algebra $S(\overline V)$, where $\overline V$ is ''the
classical limit'' of $V$. Therefore,  the Hilbert series
$h(S_\sigma(V),t)$ equals to the Hilbert series of a   quotient of
algebra of  $S(\overline V)$ by a certain homogeneous ideal
$\overline J$. It is well-known that $h(S(\overline
V))=\frac{1}{(1-t)^{\dim V}}$ and   well-known that the Hilbert
series of $S(\overline V)/\overline J$ can be expressed as a
rational function in $t$ with the denominator $(1-t)^d$. This proves
the corollary.
\endproof

%(where $P$ the weight lattice of $\gg$).  Assume that there exists a sub-algebra $\hh'\in \hh$ such that all the %$\hh'$-weights of  $V$ are on one side of some hyperplane.  Then it is clear that all $\hh'$-weight spaces of $T(V)$ %are finite-dimensional. And for any $A\in Ob(\OO_{gr,f})$ the $\hh'$-character is well-defined:
%$$ch_{{\hh}'}(A)=\sum_{\beta'\in P'} \dim A(\beta')e^{\beta'} \ ,$$
%where $P'$ is the weight lattice of $\hh'$.
%In particular, if $z\in Z(\gg)$ such that $z|_V=id_V$ and $z$ has weight $\beta_{0}$ then for $\hh'=\CC\cdot z$ one has
%$$h(A,e^{\beta_0})=ch_{\hh'}(A) \ .$$
%We will, in this case, use the notation $ch_{\hh'}(A) =ch_{z}(A)$.

We conclude the section with  a conjectural
observation about ''numerical Koszul duality'' between the braided symmetric and exterior algebras.

\begin{conjecture} For each $\lambda\in P_+$ we have
$$h(S_\sigma(V_\lambda),t)\cdot h(\Lambda_\sigma(V_\lambda^*),-t)=1+O(t^4) \ .$$
Or, equivalently, $\dim S_\sigma ^3 V_\lambda-\dim \Lambda_\sigma ^3 V_\lambda^*=(\dim V_\lambda)^2$.

\end{conjecture}

Indeed, if $S_\sigma(V_\lambda)$ is  Koszul,  then the conjecture obviously holds without the term $O(t^4)$. For non-Koszul $S_\sigma(V_\lambda)$, we have verified the conjecture numerically in several cases (see e.g., Corollary \ref{cor:symmetric exterior cube dimension} below). Moreover, in non-Koszul case we expect that the term $O(t^4)$ will never vanish.

\subsection{Flat  modules}
\label{subsect:flat modules}

We view $S_\sigma(V)$ and $\Lambda_\sigma(V)$ as deformations of the quadratic algebras  $S(V)$ and $\Lambda(V)$ respectively.
 Theorem \ref{th:flat poisson closure} (taken in conjunction with Proposition \ref{pr:algebra coalgebra duality}(b)) implies that
$$\dim  ~S_\sigma^n V=\dim  ~S_\sigma(V)_n \le \binom{\dim V+n-1} n $$
for all $n$.

Therefore, we  come to the following definition.

 \begin{definition}
\label{definition:flatness}
A finite dimensional $U_q(\gg)$-module is flat, if and only if
$$\dim ~S^n_\sigma V=\binom{\dim V +n-1}{n}  $$
for all $n\ge 0$; i.e., the braided symmetric power $S^n_\sigma V$ is isomorphic (as a vector space) to the ordinary symmetric  power $S^n V$.

\end{definition}

By definition,  $\dim ~S^n_\sigma V=\binom{\dim V+1}{n}  $ for
$n=0,1,2$. To determine whether $V$ is flat it is sufficient to
determine the flatness of  $S^3_\sigma V$ due to the following
result.

\begin{proposition}  
\label{pr:flatness}
A module $V$ in $\OO_f$  is flat if and only if
$$\dim ~S^3_\sigma V=\dim ~S^3 V=\binom{\dim V+2}{3} \ .$$
Moreover, if $S_\sigma(V)$ is flat, then it is Koszul.

\end{proposition}

\begin{proof} Follows from Drinfeld's result (\cite[Theorem 1]{DR2}) which asserts that if $A =T(V)/\langle I_2 \rangle$ is a quadratic algebra over $\CC(h)$ such that $I_2\subset V\otimes V$ is a flat $h$-deformation of the exterior square $\Lambda^2V$ (i.e., $I_2$ has a basis of the form $e_i\otimes e_j- e_j\otimes e_i-h\sum_{k,l} c_{ij}^{kl}e_k\otimes e_l$ for $1\le i<j\le \dim V$, where all $c_{ij}^{k,l}\in \CC[[h]]$) and  $\dim A_3=\dim S^3V$, then $\dim A_n=\dim S^n V$ for all $n$. The Koszulity of such an $A$ follows because $A$ has a PBW-basis. The proposition is proved. \end{proof}

The above result is a first step toward solving the following problem.

\begin{problem}
 Classify all flat  $U_q(\gg)$-modules $V$ in $\OO_f$.
\end{problem}

We expect that for each semisimple Lie algebras $\gg$ there are only finitely many flat simple modules $V_\lambda$. 
Next, we compute a lower bound for dimension of each $S^3_\sigma V_\lambda$. We will use the notation
\begin{equation}
\label{eq:notation highest vectors}
V^\mu
\end{equation}
for the space of highest weight vectors of weight $\mu$ in a module
$V$ in $\OO_f$. For $\lambda,\mu,\nu\in P_+$ denote
$c_{\lambda,\mu}^\nu=\dim (V_\lambda\otimes V_\mu)^\nu$, i.e.,
$c_{\lambda,\mu}^\nu$ is the tensor product multiplicity.  And for
any $\lambda,\mu\in P_+$ denote $c^+_{\lambda;\mu}=\dim (S^2_\sigma V_\lambda)^\mu$ and $c^-_{\lambda;\mu}=\dim (\Lambda^2_\sigma V_\lambda)^\mu$, so that $c^-_{\lambda;\mu}+c^+_{\lambda;\mu}=c_{\lambda,\lambda}^\mu$. Ultimately, define: 
$$d_\lambda^\mu:=\sum_{\nu\in P_+} (c_{\lambda;\nu}^+-c_{\lambda;\nu}^-) c_{\nu,\lambda}^\mu \ .$$

\begin{lemma}
\label{le:lower bound symmetric cube}
For each $\lambda\in P_+$ we have
\begin{equation}
\label{eq:lower bound symmetric cube}
\dim S^3_\sigma V_\lambda\ge \sum_{\mu\in P_+} \max(d_\lambda^\mu,0)\cdot  \dim V_\mu \ .
\end{equation}
\end{lemma}

\proof   By definition of $S^3_\sigma V_\lambda$, one has:
$(S^3_\sigma V_\lambda)^\mu=(S^2_\sigma V_\lambda\otimes
V_\lambda)^\mu\cap (V_\lambda\otimes S^2_\sigma V_\lambda)^\mu$.
Therefore,  we obtain the inequality:
$$\dim (S^3_\sigma V_\lambda)^\mu \ge \dim (S^2_\sigma V_\lambda\otimes V_\lambda)^\mu
+\dim (V_\lambda\otimes S^2_\sigma V_\lambda)^\mu-\dim (V_\lambda\otimes V_\lambda\otimes V_\lambda)^\mu$$
$$=\dim (S^2_\sigma V_\lambda\otimes V_\lambda)^\mu-\dim (V_\lambda\otimes \Lambda^2_\sigma V_\lambda)^\mu
=d_\lambda^\mu \ . $$
Taking into account that $\dim V=\sum_{\mu \in P_+} \dim V^\mu\cdot \dim V_\mu$, we obtain (\ref{eq:lower bound symmetric cube}).
\endproof

Using  this result, we obtain the following  sufficient criterion of flatness.
\begin{proposition}
\label{pr:flat symmetric cube criterion}
Assume that
\begin{equation}
\label{eq:flat symmetric cube equation}
\sum_{\mu\in P_+} \max(d_\lambda^\mu,0)\cdot  \dim V_\mu=\binom{\dim V_\lambda+2}{3}
\end{equation}
for some $\lambda\in P_+$. Then  the simple $U_q(\gg)$-module $V_\lambda$ is flat.
\end{proposition}

\proof Follows from (\ref{eq:lower bound symmetric cube}) and Proposition \ref{pr:flatness}.
\endproof

Based on the examples below (and on Corollary \ref{cor:symmetric
exterior cube dimension} below), we expect that this sufficient
criterion of flatness is also necessary. In what follows we will
provide some examples of flat modules $V_\lambda$,  each of which
satisfies (\ref{eq:flat symmetric cube equation}). Along these lines
we expect the ''orthogonality'' of the symmetric and exterior cubes:
if $(S^3_\sigma V_\lambda)^\mu\ne 0$ (in the notation of
(\ref{eq:notation highest vectors})) then $(\Lambda^3_\sigma
V_\lambda)^\mu=0$ and vice versa.

\subsection{Examples of flat modules} Recall that for  $U_{q}(gl_d(\CC))$-modules the dominant weights are
non-increasing $d$-tuples of integers  $\lambda=(\lambda_1\ge
\lambda_2\ge \cdots \lambda_d)$. The $d$-tuple
$\omega_i=(1,1,\ldots,1,0,\ldots,0)$ having $i$ ones and $d-i$ zeros
is referred to as the  {\it $i$-th  fundamental weight}. In particular, $V_{\omega_1}\cong (\CC(q))^d$ is the standard $U_{q}(gl_d(\CC))$-module.

\begin{lemma}
\label{le:q-polynomials}
The standard  $U_{q}(gl_d(\CC))$-module  $V_{\omega_1}$  is flat. More precisely,

\noindent (a) Each graded component $S_\sigma(V_{\omega_1})_n$ is irreducible and isomorphic to $V_{n\omega_1}$.

\noindent (b) As an algebra  $S_\sigma( V_{\omega_1})$ is isomorphic
to $\CC_q[x_1,\ldots,x_d]$, {\it $q$-polynomial algebra} generated
by $x_1,\ldots,x_n$ subject to the relations $x_jx_i=qx_ix_j$ for
$1\leq i<j\le d$.

\end{lemma}

\proof Prove (a) It suffices to show that $S_\sigma^n V_{\omega_1}\cong V_{n\omega_1}$. Indeed,
by Proposition \ref{pr:n-th symm-power} with $\lambda=\omega_1$ and $V=V_{\omega_1}$, one has an injective homomorphism
$V_{n\omega_1}\hookrightarrow S^n_\sigma V_{\omega_1}$.   And and by the dimension count, this homomorphism is also surjective.

Prove (b) now. Let us choose a standard weight basis
$x_1,\ldots,x_d$ in $V_{\omega_1}$. That is,
$E_i(x_j)=\delta_{i,j-1}x_{j-1}$ and $F_i(x_j)=\delta_{i,j}v_{j+1}$,
$K_\lambda(x_j)=q^{\lambda_j}x_j$. The formula (\ref{eq:braiding
gl_d})  implies that
\begin{equation}
\label{eq:symmetric exterior squares irreducible}
\Lambda_\sigma^2 V_{\omega_1}=\langle x_j\otimes x_i-qx_i\otimes x_j|i<j
\rangle,~S_\sigma^2 V_{\omega_1}=\langle x_i\otimes x_j+qx_j\otimes x_i| i\le j\rangle \ ,
\end{equation}
where $\langle\cdot \rangle$ denotes the $\CC(q)$-linear span.
Therefore, the quotient algebra
$S_\sigma(V_{\omega_1})=T(V_{\omega_1})/\langle \Lambda^2_\sigma
V_{\omega_1}\rangle$ is naturally isomorphic to
$\CC_q[x_1,\ldots,x_d]$. This proves (b).  The proposition is
proved.
\endproof

A more general example is the space $M_{d\times k}$ of all $d\times
k$-matrices over $\CC(q)$  regarded as a $U_q(gl_d(\CC)\times
gl_k(\CC))$-module. Recall that   $\CC_q[M_{d\times k}]$ is the
algebra of quantum regular functions on the space  $d\times
k$-matrices (for readers' convenience, the definition and main
properties of  $\CC_q[M_{d\times k}]$ are given in Appendix).

\begin{proposition}
\label{pr:q-matrices} For any $d,k\ge 1$ one has an isomorphism of
algebras in $\OO_{gr, f}$: $S_\sigma( M_{d\times k})\cong
\CC_q[M_{d\times k}]$. In particular,  the $U_q(gl_d(\CC)\times
gl_k(\CC))$-module $M_{d,k}$ is flat.
\end{proposition}
\proof Let us identify the weight lattice $P^+$ for $U_q(gl_d(\CC)\times gl_k(\CC))$ with $\ZZ^d\times \ZZ^k$. Let
 $V_1:=V_{(\omega_1,0)}$, $V_2:=V_{(0,\omega_1)}$ so that $M_{d\times k}\cong V_{(\omega_1,\omega_1)}\cong V_1\otimes V_2$
 as a $U_q(gl_d(\CC)\times gl_k(\CC))$-module.  By definition, both braiding operators
 $\RR_{V_1,V_2}:V_1\otimes V_2\to V_2\otimes V_1$  and $\RR_{V_2,V_1}:V_2\otimes V_1\to V_1\otimes V_2$ are merely permutations of factors. Therefore,
$$\RR_{M_{d\times k},M_{d\times k}}=\tau_{23}\circ (\RR_{V_1,V_1}\otimes \RR_{V_2,V_2})\circ \tau_{23}$$
where $\tau_{23}$ is the permutation of two middle factors in the tensor product of four modules. In its turn, it implies that
$$\Lambda^2_\sigma  M_{d\times k}=\tau_{23}(S^2_\sigma V_1\otimes \Lambda_\sigma^2 V_2\oplus \Lambda_\sigma^2 V_1\otimes S_\sigma^2 V_2) \ .$$
Then using the standard basis $x_1,x_2,\ldots,x_d$ for $V_1$ and the
standard $x'_1,x'_2,\ldots,x'_k$ for $V_2$, the formula
(\ref{eq:symmetric exterior squares irreducible}), and the notation
$x_{ij}:=x_i\otimes x'_j$, we see that $\Lambda^2_\sigma  M_{d\times
k}$ is freely spanned by  the elements of the form
$$x_{i,j'}\otimes x_{i,j}-qx_{i,j}\otimes x_{i,j'}, ~~x_{i,j'}\otimes x_{i',j}
+qx_{i',j'}\otimes x_{i,j}-qx_{i,j}\otimes x_{i',j'}-q^2x_{i',j}\otimes x_{i,j'}\ ,  $$
$$x_{i',j}\otimes x_{i,j}-qx_{i,j}\otimes x_{i',j}, ~~x_{i',j'}\otimes x_{i,j}
+qx_{i',j'}\otimes x_{i,j}-qx_{i,j}\otimes
x_{i',j'}-q^2x_{i,j'}\otimes x_{i',j}$$ for $1\le i< i'\le d$, $1\le
j<j'\le d$. Therefore, the quotient algebra $S_\sigma(M_{d\times
k})=T(M_{d\times k})/\langle \Lambda^2_\sigma  M_{d\times k}\rangle$
is naturally isomorphic to $\CC_q[M_{d\times k}]$.

The proposition is proved.
\endproof

\subsection{Braided symmetric and exterior powers of simple $U_q(gl_2(\CC))$-modules}
\label{subsect:gl_2} Recall that irreducible
$U_{q}(gl_2(\CC))$-modules $V_\lambda$ in $\OO_f$ are labeled by the
pairs $\lambda=(\lambda_1\ge \lambda_2)$ of integers and $\dim
V_\lambda=\lambda_1-\lambda_2+1$.

 In what follows  we will  abbreviate  $V_\ell:=V_{(\ell,0)}$ for $\ell\in \ZZ_{\ge 0}$.
 Note that in terms of braided symmetric powers, $V_\ell\cong S_\sigma^\ell V_1$ for each $\ell\ge 0$.

It is easy to see that  $(V_{(\lambda_1,\lambda_2)})^{\otimes n}\cong V_{(n\lambda_2,n\lambda_2)}\otimes (V_{\lambda_1-\lambda_2})^{\otimes n}$ and
$$ \Lambda_\sigma^n V_{(\lambda_1,\lambda_2)}\cong V_{(n\lambda_2,n\lambda_2)}\otimes \Lambda_\sigma^n
V_{\lambda_1-\lambda_2},~S_\sigma^n V_{(\lambda_1,\lambda_2)}\cong V_{(n\lambda_2,n\lambda_2)}\otimes S_\sigma^n V_{\lambda_1-\lambda_2}$$
for any $\lambda_1\ge \lambda_2$ and $n\ge 0$.
%One also has an isomorphisms $V_{(\lambda_1,\lambda_2)}\cong V_{1,1}^{\otimes \lambda_2}\bigotimes V_{(\lambda_1-\lambda_2,0)}$ and $V_{\lambda}^*\cong %V_{\lambda^*}$, where $(\lambda_1,\lambda_2)^*=(-\lambda_2,-\lambda_1)$.

The following is our first main result for $U_q(gl_2(\CC))$-modules.

\begin{maintheorem}
\label{th:exterior fourth power}

For each $\ell\ge 3$ and   $n\ge 4$ we have $\Lambda^n V_\ell=0$.

\end{maintheorem}

We will proof Theorem \ref{th:exterior fourth power} in Section
\ref{sect:proofs of main results}.

%In order to prove  Theorem  \ref{theorem:flat-sl2},  it suffices to determine the dimensions of  the braided exterior and symmetric cubes of finite dimensional irreducible representations $V_{\ell _{1},\ell _{2}}$ ( Proposition \ref{proposition:flatness}).  As $R\in U_{q}(\sl_{2})\widehat \otimes U_{q}(\sl_{2})$ it suffices to prove the theorem for $V_{\ell,0}=V_\ell$.

%%%%%%%%%%%%%%%%
%Since $\widehat R$ lies in the semisimple part of $U_{q}(gl_2(\CC))\widehat\otimes U_{q}(gl_2(\CC))$, it suffices (due to Remark \ref{rem:sl2-reduction}) to prove Theorem \ref{th:flat-sl2} for $V_{\ell,0}=V_\ell$.

Our second main result gives a complete description of  the  braided exterior and the braided symmetric cube of each $V_\ell$.

\begin{maintheorem}
\label{th:exterior powers and symmetric cube}

For each $\ell\ge 0$ and each $n\ge 3$ one has:
\begin{equation}
\label{eq:symmetric cube}
S^3_\sigma V_\ell\cong \begin{cases}
\bigoplus\limits_{0\le i\le \frac{\ell-1}{2}} V_{(3\ell -2i,2i)} &  \text {if $\ell$ is odd} \\
\bigoplus\limits_{0\le i\le \frac{3\ell}{4}} V_{(3\ell -2i,2i)} & \text {if $\ell$ is even}
\end{cases} .
\end{equation}
\begin{equation}
\label{eq:exterior powers}
\Lambda_\sigma^3 V_\ell \cong \begin{cases}
0 & \text {if $\ell$ is odd} \\
\bigoplus\limits_{\frac{\ell}{2}\le i\le \frac{3\ell-2}{4}} V_{(3\ell -2i-1,2i+1)} & \text{if $\ell$ is even} \\
\end{cases} .
\end{equation}

\end{maintheorem}

We will proof Theorem \ref{th:exterior powers and symmetric cube} in section \ref{sect:proofs of main results}.
The following is an immediate corollary of Theorem \ref{th:exterior powers and symmetric cube}.

\begin{corollary}
\label{cor:symmetric exterior cube dimension}
For each $\ell\ge 0$ we have:
%For any $\lambda=(\lambda_1\ge \lambda_2)$ one has an isomorphism in $\OO_f$:
%
%$$\Lambda_\sigma^n V_\lambda  \cong \Lambda_\sigma (V_\lambda)_n, S_\sigma^k V_\lambda  \cong S_\sigma (V_\lambda)_k$$
%for all $n\ge 0$ and all $k\le 3$.
%Also
\begin{equation}
\label{eq:exterior cube}
\dim \Lambda_\sigma^3 V_\ell =\delta_\ell\cdot \binom{\frac{\ell}{2}+1}{2} , ~~\dim S^3_\sigma V_\ell=
(\ell+1)^2 +\delta_\ell\cdot \binom{\frac{\ell}{2}+1}{2} \ ,
\end{equation}
where $\delta_\ell=0$ if $\ell$ is odd and $\delta_\ell=1$ if $\ell$ is even. In particular, $V_\ell$ is flat if and only if $\ell\in \{0,1,2\}$.

\end{corollary}

%By definition, $\Lambda^n V_\ell$ is flat for $n\le 2$,  more precisely,
%$\Lambda^0_\sigma V_\ell\cong V_0$, $\Lambda^1_\sigma V_\ell=V_\ell$, and
%$$\Lambda^2 V_\ell\cong \bigoplus\limits_{0\le i\le (\ell -1) /2}V_{(2\ell-2i-1, 2i+1)} \ .$$

Based on Theorem \ref{th:exterior powers and symmetric cube}, we
propose the following conjectural description of the higher braided
symmetric powers of  $V_\lambda$ for any $\lambda$.
\begin{conjecture}
\label{conj:multiplicity-free}
For any  $\ell\ge 0$ and any $n\ge 4$ one has:
\begin{equation}
\label{eq:multiplicity-free}
S^n_\sigma V_\ell \cong \begin{cases}
\bigoplus\limits_{0\le i\le \frac{\ell-1}{2}} V_{(n\ell -2i,2i)} &  \text {if $\ell$ is odd} \\
\bigoplus\limits_{0\le i\le \frac{n\ell }{4}} V_{(n\ell -2i,2i)} & \text {if $\ell$ is even}
\end{cases}.
\end{equation}
In particular,  $S^n_\sigma V_\ell$ is multiplicity-free and we have for $n\ge 4$:
\begin{equation}
\label{eq:dimension symmetric power}
\dim S^n_\sigma V_\ell =\begin{cases}
\frac{(\ell +1)(\ell (n-1)+2)}{2} & \text {if $\ell$ is odd} \\
\binom{\frac{n\ell }{2} +2}{2} & \text {if $\ell$ is even}
\end{cases} .
\end{equation}
\end{conjecture}

%In support of the conjectures we prove the following

%\begin{proposition}
%The braided symmetric fourth power of $V_{3}$  splits

%  $$S_\sigma^{4}V_{3}\cong V_{12}\oplus V_{8}$$

%\end{proposition}

So far we have collected some partial evidence for for Conjecture
\ref{conj:multiplicity-free}, obtained from analyzing the Poisson
closure $\overline {S(\overline V_\ell)}$ of the bracketed
superalgebra  $(S(\overline V_\ell),\{\cdot,\cdot\})$ (see
(\ref{eq:barSV barLambdaV})) which serves as an upper bound for the braided symmetric algebra $S_\sigma(V_\ell)$ - we expect that the upper bound is achieved for all $\ell$.  For instance, using the computer
algebra system SINGULAR (\cite{GPS01}), we have verified that for $\ell=3,4,5,6$ the dimension of each component $\overline {S(\overline V_\ell)}_n$ equals to the right hand side of
(\ref{eq:dimension symmetric power}). Taken  in conjunction with
Proposition \ref{pr:algebra coalgebra duality} and Theorem
\ref{th:flat poisson closure}(a), this  gives a correct upper bound
on $\dim S^n_\sigma V_\ell$.

Another piece of evidence comes from Proposition \ref{pr:n-th
symm-power} which implies that $S_\sigma^n V_\ell$ always contains a
unique copy of the simple $gl_2(\CC)$-module $\overline V_{n\ell}$.
We also verified  that for $\ell\ge 3$ and any $n\ge 4$ the
$gl_2(\CC)$-module $\overline {S(\overline V_\ell)}_n$ contains a
copy of $\overline V_{(n\ell-2,2)}$.

Using \eqref{eq:exterior cube} we now describe all those modules $V_\ell$ whose the  braided symmetric algebra  is Koszul.

\begin{corollary}
\label{cor:non-Koszul}
The braided symmetric algebra $S_\sigma(V_\ell)$ is  Koszul if and only if $\ell\le 2$.
\end{corollary}
\begin{proof}
For $\ell\le 2$ the module $V_\ell$ is flat and the assertion
follows from Proposition \ref{pr:flatness}. Let now $\ell\ge 3$. If
$S_\sigma(V_\ell)$ was Koszul, then the Hilbert series
$h(S_\sigma(V_\ell),t)$ (see Definition \ref{def:Hilbert
series}) would satisfy:
\begin{equation}
\label{eq:would-be Koszulity}
h(S_\sigma(V_\ell),t)=\frac{1}{h(S_\sigma(V_\ell)^!,-t)}\ .
\end{equation}
Since $S_\sigma(V_\ell)^!\cong \Lambda_\sigma ((V_\ell)^*)$ by
Proposition \ref{pr:algebra coalgebra duality}(c), $\Lambda_\sigma
((V_\ell)^*)$ and $\Lambda_\sigma (V_\ell))$ are isomorphic as
graded vector spaces, the Koszulity of $S_\sigma(V_\ell)$ would
imply that all the coefficients in the power series expansion for
the rational function
$\displaystyle{\frac{1}{h(\Lambda_\sigma(V_\ell),-t)}}$  are
non-negative.

According to Corollary \ref{cor:symmetric exterior cube dimension}, one has
$$h( \Lambda_\sigma(V_\ell),-t)=
1-(\ell+1)t+\frac{\ell (\ell+1)}{2}t^2-\delta_\ell\cdot {\frac{\ell(\ell+2)}{8}}t^3 \ .$$
Therefore,
$$\frac{1}{h( \Lambda_\sigma(V_\ell),-t)}=1+mt+\frac{m(m+1)}{2}t^2+(m^2+\delta_{m-1}\cdot {\frac{m^2-1}{8}})t^3+ct^4+O(t^5) \ ,$$
where $m=\ell+1$ and  $c=m-\frac{m(m^2-1)( m-4-\delta_{m-1})}{4}$.
This implies that $c<0$ for all  $\ell\ge 4$. In the remaining case
$\ell=3$  it is easy to see that:
$$\frac{1}{h( \Lambda_\sigma(V_3),-t)}=1+4t+10t^2+16t^3+4t^4-80t^5+O(t^6) \ .$$
Therefore, (\ref{eq:would-be Koszulity}) makes no sense for $\ell\ge 3$ and $S_\sigma(V_\ell)$ is not Koszul.
The corollary is proved.
\end{proof}

We will conclude the section with a yet conjectural
computation  of  the Hilbert series of (non-Koszul) quadratic algebras $S_\sigma(V_\ell)$.

\begin{corollary}[from Conjecture \ref{conj:multiplicity-free}]
\label{cor:character hilbert series}

\noindent(a) If $\ell$ is odd, then
$$ h(S_\sigma(V_\ell),t)=  1+ (\ell+1)\frac{t}{1-t}+\binom{\ell+1}{2}\frac{t^2}{(1-t)^2}\ .$$
\noindent(b) If $\ell$ is even, then
$$ h(S_\sigma(V_\ell),t)=
1+\ell t+\frac{t}{1-t}+\left(\binom{\frac{\ell}{2}}{2}+\ell\right)\left( \frac{t^2}{1-t}+\frac{t^2}{(1-t)^2}\right)+\frac{\ell^2}{4}\frac{t^2}{(1-t)^3} \ .$$

\end{corollary}

\proof Prove (a) first.   According to \eqref{eq:dimension symmetric power}, one has for $n\ge 1$: 
$$\dim (S_\sigma V_\ell)_n=\dim S_\sigma^n V_\ell=\frac{(\ell +1)(\ell (n-1)+2)}{2}=\binom{\ell+1}{2}(n-1)+\ell+1 \ .$$
Therefore the Hilbert series of the co-algebra $S_\sigma V_\ell$ is given by:
$$h(S_\sigma (V_\ell),t)=1+\sum_{n\ge 1} \binom{\ell+1}{2}(n-1)t^n+(\ell+1)t^n  =1+ (\ell+1)\frac{t}{1-t}+\binom{\ell+1}{2}\frac{t^2}{(1-t)^2}\ .  $$
Part (a) is proved. Prove (b) now.   Using  \eqref{eq:dimension symmetric power} again, one has for $n\ge 2$: 
$$\dim S_\sigma^n V_\ell=\binom{\frac{n\ell }{2} +2}{2}=\frac{\ell^2}{4}\binom{n}{2}+ \frac{\ell^2+6\ell}{8}n+1 \ .$$
Therefore,
$$h(S_\sigma (V_\ell),t)=1+((\ell+1)t-\binom{\frac{\ell}{2}+2}{2})+\sum_{n\ge 1} \frac{\ell^2}{4}\binom{n}{2}t^n+ \frac{\ell^2+6\ell}{8}nt^n+t^n$$
$$= 1+\ell t+\frac{t}{1-t}+\left(\binom{\frac{\ell}{2}}{2}+\ell\right)\left( \frac{t^2}{1-t}+\frac{t^2}{(1-t)^2}\right)+\frac{\ell^2}{4}\frac{t^2}{(1-t)^3}  \ .$$
This proves part (b). The corollary is proved.
\endproof

This result implies, in particular, that the associated scheme $X(\overline V_\ell)$ (see Problem \ref{pr:schemes}) is a curve if $\ell$ is odd and a surface if $\ell$ is even. 

\begin{remark} It follows from  Corollary \ref{cor:character hilbert series} that the Hilbert series for $S_\sigma (V_{2k})$ is regular at $\infty$ and  
$\lim\limits_{t\to \infty} h(S_\sigma(V_{2k}),t)=\binom{2k}{2}$. We expect that this phenomenon is caused by vanishing of the cubic component in the  dual algebra $\Lambda_\sigma(V_{2k}^*)=S_\sigma(V_{2k})^!$ (Theorem \ref{th:exterior powers and symmetric cube}). 
\end{remark}

\section{Proofs of main results}

\label{sect:proofs of main results}

\subsection{Bracketed superalgebras and the proof of Theorem \ref{th:symmetric exterior almost poisson}}
We first formulate and prove a number of results about bracketed
superalgebras. The following fact is obvious:

\begin{lemma} $~$
\label{le:category bracketed}
Tensor product of bracketed superalgebras is a bracketed superalgebra and tensor product of Poisson superalgebras is a Poisson superalgebra.
 \end{lemma}

\begin{definition}
\label{def:homogeneous}
We say that a bracketed superalgebra $A=A_{\bar 0}\oplus A_{\bar 1}$ is {\it homogeneous} if:

\noindent (i) $A_{\bar 0}=\oplus_{n\ge 0} A_{\bar 0,n}$,
$A_1=\oplus_{n\ge 0} A_{\bar 1,n}$ is $\ZZ_2\times \ZZ_{\ge
0}$-grading of $A$ as algebra, i.e., $A_{\varepsilon,m}\cdot
A_{\delta,n}\subset A_{\varepsilon+ \delta,m+n}$ for any
$\varepsilon,\delta\in \ZZ_2=\ZZ/(2)=\{\overline 0,\overline 1\}$
and any $m,n\in \ZZ_{\ge 0}$;

\noindent (ii) $A$ is generated by $A_{\bar 0,1}\oplus A_{\bar 1,1}$ as an algebra;

\noindent (iii) $\{A_{\varepsilon,m},A_{\delta,n}\}\subset
A_{\varepsilon+ \delta,m+n}$ for any $\varepsilon,\delta\in \ZZ_2$
and any $m,n\in \ZZ_{\ge 0}$.
\end{definition}

The following fact is obvious.

\begin{lemma}  The tensor product of homogeneous bracketed superalgebras is also a homogeneous bracketed superalgebra.

\end{lemma}

For any homogeneous bracketed superalgebra $A$ define the subspaces $I_n\subset A_n$, $n=3,4,\ldots$ recursively:
$$I_3(A)=  Span\{J(a,b,c)\,|\,a,b,c\in A_{\bar 0,1}+A_{\bar 1,1}\}$$
and $I_{n+1}(A)=\{A_{\bar 0,1}+A_{\bar 1,1},I_n(A)\}$ for $n\ge 3$.
%The following result gives the characterization of the ideal $I(A)$.
Denote by $I(A)$  the (super)ideal generated by all $I_n(A)$. We say
that a (super)ideal $J\subset A$ is {\it Poisson} if it contains
$I_3(A)$ and $\{A,J\}\subset J$. Equivalently, $J$ is Poisson if the
quotient superalgebra $A/J$ is naturally bracketed and Poisson.

\begin{proposition}
\label{pr:universal poisson closure}

\noindent  For any homogeneous bracketed superalgebra one has

\noindent (a) $I(A)$ is Poisson.

\noindent (b) Any Poisson (super)ideal $J$ in $A$ contains $I(A)$.

\end{proposition}

\proof We need the following fact.

\begin{lemma}
Let $A$ be a homogeneous bracketed superalgebra. Then $A$ is Poisson  if and only if $J(u,v,w)=0$ for all $u,v,w\in A_{\bar 0,1}+A_{\bar 1,1}$.

\end{lemma}

\proof
One easily checks that for any $a\in A_{\varepsilon}, b\in A_{\delta},c\in A_{\gamma}$, and $d\in A_{\theta}$ one has
$$J(a,b,c\cdot d)=(-1)^{\delta \gamma}cJ(a,b, d)+(-1)^{\varepsilon \theta} J(a,b,c)d$$
and  $J(c,a,b)=J(a,b,c)$. An obvious induction by the degree
(defined via the direct sum decomposition of $A$ ) completes the
proof.
\endproof

This proves part (a). Part (b) is obvious because by the very definition,  $I(A)$ is the smallest Poisson ideal in $A$.
\endproof

\begin{definition}
\label{def:almost Poisson}
 We say that a homogeneous bracketed superalgebra is {\it almost Poisson}, if  the ideal
 $I(A)$ is generated by its third component $I_3(A)$, i.e. the obstruction for the algebra $A$ to  be Poisson lies entirely in the  $3$rd degree.

 \end{definition}

\begin{lemma}
The tensor product of two almost Poisson superalgebras is almost Poisson.
\end{lemma}

\proof It is easy to see that
$I_3(A\otimes B)=I_3(A)\otimes 1+1\otimes I_3(B)$
for any homogeneous bracketed superalgebras $A$ and $B$. That is,
$$I_3(A\otimes B)\cdot (A \otimes B)=I_3(A)\cdot A\otimes B+A\otimes I_3(B)\cdot B=I(A)\otimes B+A\otimes I(B)$$
because  $I(A)=I_3(A)\cdot A$ and $I(A)=I_3(A)\cdot A$. Therefore, by Proposition \ref{pr:tensor product Poisson},
$$I_3(A\otimes B)\cdot (A \otimes B)=I(A\otimes B)$$

The lemma is proved.
\endproof

For a homogeneous bracketed superalgebra $A$ denote $P(A):=A/I(A)$ and will refer to $P(A)$ as the {\it Poisson closure} of $A$.

 The following result is an immediate corollary of Proposition \ref{pr:universal poisson closure}.

\begin{corollary} The  Poisson closure  $P(A)$ is characterized by the following universal property: Let  $P$  be any Poisson superalgebra.
Then any homomorphism of bracketed superalgebras $f: A\to P$
factors through $P(A)$, i.e., there exists a surjective homomorphism
$\pi:A\twoheadrightarrow P(A)$ of bracketed algebras and a
homomorphism $j:P(A)\to P$ of Poisson algebras such that $f=j\circ
\pi$.
% $$\xymatrix{A \ar[r]^{\pi}\ar[d]& P\\
%                    P(A)\ar[ur] & . }\ .$$

\end{corollary}

Clearly, the correspondence $A\mapsto P(A)$ is a functor from the
category of homogeneous bracketed superalgebras to the category of
Poisson superalgebras.

\begin{proposition}
\label{pr:tensor product Poisson} Let $A$ and $B$ be homogeneous
bracketed  superalgebras. Then one has a canonical isomorphism of
Poisson algebras: $P(A\otimes B)= P(A)\otimes P(B)$. Equivalently,
$I(A\otimes B)= I(A)\otimes B+A\otimes I(B)$.
\end{proposition}

\proof By Lemma \ref{le:category bracketed}, $P(A)\otimes P(B)$ is
Poisson as a product of Poisson superalgebras. One also has a
surjective algebra homomorphism $A\otimes B\to P(A)\otimes P(B)$.
Due to the universality of the Poisson closure $P(A\otimes B)$, the
latter homomorphism  factors through the homomorphism  $P(A\otimes
B)\to P(A)\otimes P(B)$ of Poisson algebras. Therefore,
$$I(A\otimes B)\subseteq  I(A)\otimes B+A\otimes I(B)\ .$$
On the other hand, $I_3(A\otimes B)$ contains both $I_3(A\otimes
1)=I_3(A)\otimes 1$ and $I_3(1\otimes B)=1\otimes I_3(B)$. This
implies that  $I_n(A\otimes B)$ contains both $I_n(A)\otimes 1$ and
$1\otimes I_n(B)$. Therefore, $I(A\otimes B)$ contains both
$I(A)\otimes 1$ and $1\otimes I(B)$ and we obtain the opposite
inclusion:
$$I(A\otimes B)\supseteq  I(A)\otimes B+A\otimes I(B)\ .$$

The obtained double inclusion proves the assertion.
\endproof

%One  easily checks
%$$I_{3}(A\otimes B)=I_{3}(A)\otimes 1+1\otimes I_{3}(B)\ .$$
%Then, one obtains inductively
%$$\{A\otimes B, I_{n}(A\otimes B)\}\subset I_{n+1}(A)\otimes B+A\otimes I_{n+1}(B)\subset I(A)\otimes B+A\otimes I(B)\ . $$

Now we are ready to finish the proof of Theorem \ref{th:symmetric exterior almost poisson}.

Let ${\bf U}$ be a  coalgebra with the co-product $\Delta:{\bf U}\to
{\bf U}\otimes {\bf U}$ and let $A=\bigoplus\limits_{n\ge 0} A_n$
(where $A_n=A_{\overline 0,n}\oplus A_{\overline 1,n}$) be a
homogeneous bracketed superalgebra. We say that ${\bf U}$ is {\it
acting} on $A$ if  $A$ is a ${\bf U}$-module algebra, ${\bf
U}(A_{\varepsilon,n})\subset A_{\varepsilon,n}$ for each
$\varepsilon\in \ZZ_2$, $n\ge 0$,  and the bracket satisfies:
\begin{equation}
\label{eq:homogeneous bracket differential}
\{u,v\}={\bf m}_2({\bf r}(u\otimes v))
\end{equation}
for each $u, v\in A_1$, where ${\bf m}_2:A_1\otimes A_1\to A_2$ is
the multiplication, and ${\bf r}\in {\bf U}\wedge {\bf U}$ is an
element such that $(\Delta\otimes id)({\bf r})={\bf r}_{13}+{\bf
r}_{23}$.

\begin{proposition}
\label{pr:general almost poisson}

Let ${\bf U}$ be a co-algebra acting on a homogeneous bracketed superalgebra $A$. Then:

\noindent (a)
$\{a,b\}={\bf m}({\bf r}(a\otimes b))$ for all $a, b\in A$,
where ${\bf m}:A\otimes A\to A$ is the multiplication.

\noindent (b)  $\{B,C\}\subset B\cdot C$ for any ${\bf U}$-submodules $B$ and $C$ of $A$. In particular, $A$ is almost Poisson.

\end{proposition}

\proof Prove (a). We will proceed by induction on degree. Indeed,
let $a\in A_{\varepsilon,n_1}$, $b\in A_{\delta,n_2}$, $c\in
A_{n_3}$ such that $n_1\le n$, $n_2+n_3=n$, $n_2>0$, $n_3>0$.  The
inductive hypothesis (applied when $n_2<n$) reads:
$$\{a,b\}={\bf r}_{(1)}(a)\cdot {\bf r}_{(2)}(b), ~\{a,c\}={\bf r}_{(1)}(a)\cdot {\bf r}_{(2)}(c)\ ,$$
where we have used the Sweedler notation ${\bf r}={\bf r}_{(1)}\otimes {\bf r}_{(2)}$. Then
$$\{a,bc\}=\{a,b\}c+(-1)^{\varepsilon \delta} b\{a,c\}
= {\bf r}_{(1)}(a)\cdot {\bf r}_{(2)}(b)\cdot c+(-1)^{\varepsilon \delta} b\cdot {\bf r}_{(1)}(a)\cdot {\bf r}_{(2)}(c)$$
$$= {\bf r}_{(1)}(a)({\bf r}_{(2)}(b)\cdot c+b\cdot {\bf r}_{(2)}(c))={\bf m}({\bf r}(a\otimes bc))\ .$$
This proves part (a). To prove (b) note that  $\{B,C\}\subset {\bf
r}_{(1)}(B)\cdot {\bf r}_{(2)}(C) \subset B\cdot C$ because ${\bf
r}_{(1)}(B)\subset B$, ${\bf r}_{(2)}(C)\subset C$. In particular,
taking $B=\overline V$ and $C=J(\overline V,\overline V,\overline
V)$ proves that $A$ is almost Poisson. The proposition is proved.
\endproof

Now we will set ${\bf U}:=U(\gg)$, ${\bf r}=r^-=r-\tau(r)$. Then
setting respectively  $A=S(\overline V)$ and $A=\Lambda(\overline
V)$ with the bracket given by (\ref{eq:homogeneous bracket
differential}) and using Proposition \ref{pr:general almost
poisson}(a) finishes the proof of  Theorem \ref{th:symmetric
exterior almost poisson}(a). Then Proposition \ref{pr:general almost
poisson}(b) guarantees that both $S(\overline V)$ and
$\Lambda(\overline V)$ are almost Poisson in the sense of Definition
\ref{def:almost Poisson}. Therefore, $\overline {S(\overline V)}$
and $\overline {S(\overline V)}$ are Poisson superalgebras.  This
finishes the proof of  Theorem \ref{th:symmetric exterior almost
poisson}(b).

Theorem \ref{th:symmetric exterior almost poisson} is proved. \endproof

%\subsection{Open problems}

%\begin{problem} ??? morphisms are no longer defined?????
% Classify all the pairs $U,V\in Ob(\OO_f)$ such that  $S_\sigma(U\oplus V)\cong S_\sigma(U)\otimes_\RR S_\sigma(V)$ in $\OO_{gr,f}$ and/or %$\Lambda_\sigma(U\oplus V)\cong \Lambda_\sigma(U)\otimes_{-\RR} \Lambda_\sigma(V)$.
%\end{problem}

\subsection{$({\bf k},{\bf A})$-algebras and the proof of Theorem \ref{th:flat poisson closure}}
In this section we   will develop  a general framework of $({\bf
k},{\bf A})$-algebras and will prove Theorem \ref{th:flat poisson closure} along with its 
generalizations in the category of $({\bf k},{\bf A})$-algebras.

Let ${\bf k}$ be a field and ${\bf A}$ be a local subring of ${\bf
k}$. Denote by $\mm$ the only maximal ideal in ${\bf A}$ and by
$\tilde {\bf k}$ the  residue field of ${\bf A}$, i.e., $\tilde {\bf
k}:={\bf A}/\mm$.

%Clearly, any ${\bf A}$-sub-module $L$ of a  ${\bf k}$-vector space $V$ is a torsion-free.   Note, that therefore  any finitely generated %${\bf A}$-sub-module $L$ of   $V$ is a free,   because   ${\bf A}$ is a PID.
We say that an ${\bf A}$-submodule $L$ of a ${\bf k}$-vector space
$V$ is an {\it ${\bf A}$-lattice} of $V$ if $L$ is a free ${\bf
A}$-module and ${\bf k}\otimes_{\bf A} L=V$, i.e., $L$ spans $V$ as
a ${\bf k}$-vector space. Note that for any ${\bf k}$-vector space
$V$ and any ${\bf k}$-linear basis ${\bf B}$ of $V$ the ${\bf
A}$-span $L={\bf A}\cdot {\bf B}$ is an ${\bf A}$-lattice in $V$.
Conversely, if $L$ is an ${\bf A}$-lattice in $V$, then any ${\bf
A}$-linear basis ${\bf B}$ of $L$ is also a ${\bf k}$-linear basis
of $V$.

Denote by $({\bf k},{\bf A})-Mod$ the  category which objects are
pairs $\VV=(V,L)$ where $V$ is a ${\bf k}$-vector space and
$L\subset V$ is an ${\bf A}$-lattice of $V$; an arrow $(V,L)\to
(V',L')$ is any ${\bf k}$-linear map $f:V\to V'$ such that
$f(L)\subset L'$.

Clearly, $({\bf k},{\bf A})-Mod$ is an abelian category. Moreover,
$({\bf k},{\bf A})-Mod$ is ${\bf A}$-linear because each
$Hom(\UU,\VV)$ in $({\bf k},{\bf A})-Mod$ is an ${\bf A}$-module.

\begin{definition}
We say that a functor $F:\mathcal C\to \mathcal D$ is {\it almost equivalence} of $\mathcal C$ and $\mathcal D$ if:

\noindent (a) for any objects $c,c'$ of $\mathcal C$ an isomorphism $F(c)\cong F(c')$ in $\mathcal D$ implies that $c\cong c'$ in  ${\mathcal C}$;

\noindent (b) for any object in $d$ there exists an object $c$ in ${\mathcal C}$ such that $F(c)\cong d$ in ${\mathcal D}$.

\end{definition}

\begin{lemma}
\label{le:trafo of bases} The forgetful functor $({\bf k},{\bf
A})-Mod \to  {\bf k}-Mod$ given by $(V,L)\mapsto V$ is an almost
equivalence of categories.

\end{lemma}

\proof It suffices to show that any two objects of the form $(V,L)$
and $(V,L')$ are isomorphic in $({\bf k},{\bf A})-Mod$. This follows
from the fact that there exist ${\bf k}$-linear bases ${\bf B}$ and
${\bf B}'$ for $V$ such that $L={\bf A}\cdot {\bf B}$ and $L'={\bf
A}\cdot {\bf B}'$. Fix any bijection $\varphi:{\bf B}\widetilde \to
{\bf B}'$. This extends to  a ${\bf k}$-linear automorphism $f:V\to
V$ such that $f(L)=L'$. This $f$ is a desirable isomorphism
$(V,L)\widetilde \to (V,L')$.
\endproof

The following fact is, apparently, well-known.

\begin{lemma}
\label{le:tensor category}
$({\bf k},{\bf A})-Mod$ is a symmetric tensor category under the  operation
$$(U,L)\otimes (V,L')=(U\otimes_{\bf k} V,L\otimes_{\bf A} L') $$
for any objects $(U,L)$ and $(V,L')$ of $({\bf k},{\bf A})-Mod$.
\end{lemma}

\proof Let $(U,L)$ and $(V,L')$ be objects of $({\bf k},{\bf A})-Mod$. Clearly,
$$U\otimes_{\bf k} V= ({\bf k}\otimes_{\bf A} L)\otimes_{\bf k} ({\bf k}\otimes_{\bf A} L') \ .$$
This taken together with the following general identity
$$({\bf k}\otimes_{\bf A} L)\otimes_{\bf k} ({\bf k}\otimes_{\bf A} L')={\bf k}\otimes_{\bf A} (L \otimes_{\bf A} L')$$
implies that $L \otimes_{\bf A} L'$ is an ${\bf A}$-lattice in $U\otimes_{\bf k} V$.

The associativity of the product is obvious. And the unit object is $({\bf k},{\bf A})$.

The lemma is proved.
\endproof

Define a functor $\FF:({\bf k},{\bf A})-Mod \to \tilde {\bf k}-Mod$ as follows. For any object
$(V,L)$ of $({\bf k},{\bf A})-Mod$ we set $\FF(V,L)=L/\mm L$
and for any morphism $f:(V,L)\to (V',L')$ we set
$\FF(f):L/\mm L\to L'/\mm L'$ to be a natural  $\tilde {\bf k}$-linear map.

\begin{lemma}

\label{le:specialization functor} $\FF:({\bf k},{\bf A})-Mod \to \tilde {\bf k}-Mod$ is a tensor functor and  an almost equivalence of categories.
\end{lemma}

\proof First, let us show that $\FF$ is a tensor functor. It follows
from the fact that for any ${\bf A}$-module $M$ and any ${\bf
A}$-module $M'$ we have a canonical isomorphism (even without any
restrictions on a commutative ring ${\bf A}$ and an ideal
$\mm\subset {\bf A}$):
\begin{equation}
\label{eq:tensor product quotient}
(M/\mm  M)\bigotimes_{{\bf A}/\mm}(M'/\mm M' )= (M\bigotimes_{\bf A} M)/(\mm M\bigotimes_{{\bf A}} M') \ .
\end{equation}

The fact that $\FF$ is an almost equivalence of categories follows from the following two observations.

1. $(V,L)\cong (V',L')$ in $({\bf k},{\bf A})-Mod$ if and only if
$V\cong V'$ as ${\bf k}$-vector spaces. This is a direct consequence
of the proof of Lemma \ref{le:trafo of bases}.

2. For any basis $\overline {\bf B}$ of  $\overline V$ in  $\tilde
{\bf k}-Mod$ we have $(Span_{\bf k}[\overline {\bf B}] ,Span_{\bf
A}[\overline {\bf B}])$ is an object of $({\bf k},{\bf A})-Mod$ and
$\FF(Span_{\bf k}[\overline {\bf B}] ,Span_{\bf A}[\overline {\bf
B}])=Span_{\tilde {\bf k}}\{\overline {\bf B}\}   \cong \overline
V$.

The lemma is proved.
\endproof

In what follows we will sometimes abbreviate $\overline V=\FF(V,L)$ and $\overline f=\FF(f)$. Similarly, for an algebra ${\mathcal A}=(A,L_A)$ in the $({\bf k},{\bf A})-Mod$ we will sometimes abbreviate $\overline {\mathcal A}=\FF{\mathcal A}$, the algebra over the residue field $\tilde {\bf k}={\bf A}/\mm$.

\medskip

Let $U$ be a ${\bf k}$-Hopf algebra and let $U_{\bf A}$ be a Hopf
${\bf A}$-subalgebra of $U$. This means that $\Delta(U_{\bf
A})\subset U_{\bf A}\otimes_{\bf A} U_{\bf A}$ (where $U_{\bf
A}\otimes_{\bf A} U_{\bf A}$ is naturally an ${\bf A}$-subalgebra
of $U\otimes_{\bf k} U$), $\varepsilon(U_{\bf A})\subset {\bf A}$,
and $S(U_{\bf A})\subset U_{\bf A}$.
%such that ${\bf k}\otimes_{\bf A} U_{\bf A}=U$, i.e., $U_{\bf A}$ spans $U$.
We will refer to the above pair  $\UU=(U,U_{\bf A})$ as to $({\bf
k},{\bf A})$-{\it Hopf algebra} (please note that $U_{\bf A}$ is not
necessarily a free ${\bf A}$-module, that is, $\UU$ is not
necessarily a $({\bf k},{\bf A})$-module).

Given a $({\bf k},{\bf A})$-Hopf algebra $\UU=(U,U_{\bf A})$, we say
that an object $\VV=(V,L)$ of $({\bf k}, {\bf A})-Mod$  is a
$\UU$-module if $V$ is a $U$-module and $L$ is an $U_{\bf
A}$-module.

Denote by $\UU-Mod$ the category which objects are $\UU$-modules and
arrows are those morphisms of $({\bf k}, {\bf A})$-modules which
commute with the $\UU$-action.

 Clearly, for $({\bf k},{\bf A})$-Hopf algebra $\UU=(U,U_{\bf A})$
 the category $\UU-Mod$ is a tensor (but not necessarily symmetric) category.

 For each $({\bf k},{\bf A})$-Hopf algebra $\UU=(U,U_{\bf A})$ we define $\overline \UU:=U_{\bf A}/\mm U_{\bf A}$.
 Clearly, $\overline \UU$ is a Hopf algebra over $\tilde {\bf k}={\bf A}/\mm$.

The following fact is obvious.

\begin{lemma}
\label{le:Hopf modules}
 In the notation of Lemma \ref{le:specialization functor}, for any   $({\bf k},{\bf A})$-Hopf algebra $\UU$ the functor $\FF$
 naturally extends to a tensor functor
\begin{equation}
\label{eq:dequantization functor}
\UU-Mod\to \overline \UU-Mod \ .
\end{equation}

\end{lemma}

%Since $\VV$ is

%a $\UU$-module, the action of  $\UU$ yields $({\bf k}, {\bf A})$-module homomorphisms
%$$\UU^{\otimes i}\otimes \VV\to \VV$$, which descends to a $\tilde k$-linear maps
%$$\FF ( \UU^{\otimes i}\otimes \VV)\to \overline V\ .$$

%From Lemma \ref{le:tensors} one obtains maps

%$$\left( \FF (\UU)\right)^{\otimes i} \otimes \overline V\to \overline V$$.

%For $\VV=\UU$ this yields the algebra structure on $\overline U$ and for arbitrary $\VV$ the $\overline U$-module structure.
%\endproof

Now let ${\bf k}=\CC(q)$ and ${\bf A}$ is the ring of all those rational functions in $q$ which are defined at $q=1$.
Clearly, ${\bf A}$ is a local PID with maximal ideal $\mm=(q-1){\bf A}$ (and, moreover, each ideal in ${\bf A}$ is of the form $\mm^n=(q-1)^n{\bf A}$).
Therefore,  $\tilde {\bf k}:={\bf A}/\mm=\CC$.

Recall from Section \ref{sect:definitions} that $U_{q}(\gg)$ is the
quantized universal enveloping algebra (of a reductive Lie algebra
Let $\gg$) generated (over ${\bf k}=\CC(q)$) by $E_i,F_i$,
$i=1,2,\ldots,r$, and all $K_\lambda$, $\lambda\in P$. Denote
$h_\lambda=\frac{K_\lambda-1}{q-1}$ and let $U_{{\bf A}}(\gg)$ be
the ${\bf A}$-algebra generated by all  $h_\lambda$, $\lambda\in P$
and all $E_ i, F_i$.

Denote by $\UU_{q}(\gg)$ the pair $(U_{q}(\gg), U_{{\bf A}}(\gg))$.

\begin{lemma}
\label{le:Hopf pair}
The pair $\UU_{q}(\gg)=(U_{q}(\gg), U_{{\bf A}}(\gg))$ is a $({\bf k},{\bf A})$-Hopf algebra.

\end{lemma}

\proof All we have to show is that $U_{{\bf A}}(\gg)$ is a Hopf
${\bf A}$-subalgebra of $U_{q}(\gg)$. That is, we have to show that
$U_{{\bf A}}(\gg)$ is closed under the coproduct, counit and the
antipode.  Clearly,  $\Delta(E_i), \Delta(F_i) \in U_{{\bf
A}}(\gg)\otimes_{\bf A} U_{{\bf A}}(\gg)$, $\varepsilon(E_i),
\varepsilon(F_i)\in {\bf A}$, and $S(E_i), S(F_i) \in U_{{\bf
A}}(\gg)$ for all $\lambda\in P$ and $i=1,2,\dots,r$. Therefore, it
remains to show that  $\Delta(h_\lambda)\in U_{{\bf
A}}(\gg)\otimes_{\bf A} U_{{\bf A}}(\gg)$, $\epsilon(h_\lambda)\in
{\bf A}$, and $S(h_\lambda)\in U_{{\bf A}}(\gg)$. It is easy to see
that
$$\Delta(h_\lambda)= \frac{K_\lambda\otimes K_\lambda-1\otimes 1}{q-1}=h_\lambda \otimes K_\lambda+1\otimes h_\lambda\ .$$
$$\epsilon(h_\lambda)= 0, S(h_\lambda)=-h_\lambda K_{-\lambda}\ .$$

This proves the Lemma because $K_\lambda=(q-1)h_\lambda+1$, e.g., all $K_\lambda\in U_{{\bf A}}(\gg)$.  \endproof

 \begin {lemma}
 \label{le:U(g)-action} We have $\overline {\UU_{q}(\gg)}=U(\gg)$.
 \end{lemma}

\proof The Hopf algebra $\overline {\UU_{q}(\gg)}$ is generated (over $\tilde {\bf k}=\CC$) by
$$\overline h_{\lambda}=h_{\lambda} \mod (q-1)U_{\bf A}(\gg), ~ \overline K_{\lambda}=K_{\lambda} \mod (q-1)U_{\bf A}(\gg)$$
for $\lambda\in  P$ and  by
$$\overline E_i=E_i\mod (q-1)U_{\bf A}(\gg),~\overline F_i=F_i \mod (q-1)U_{\bf A}(\gg)$$
for $i=1,\ldots r$. Note that $\overline K_{\lambda}=1$ for all $\lambda$ because $K_\lambda=(q-1)h_\lambda+1$. Furthermore, the relations
$$h_{\lambda+\mu}=(q-1)h_{\lambda}h_\mu+h_\lambda+h_\mu$$
imply that $\overline h_{\lambda+\mu}=\overline
h_{\lambda}+\overline h_{\mu}$ for all $\lambda,\mu\in P$. That is,
the $\CC$-span of all  $\overline h_{\lambda}$ is the Cartan
subalgebra of $\gg$. Next, the relations in $\UU_{q}(\gg)$:
$$ h_\lambda E_i-E_ih_\lambda =\frac{q^{(\lambda\,|\,\alpha_i)}-1}{q-1}
E_iK_{\lambda},~h_\lambda F_i-F_ih_\lambda =\frac{q^{-(\lambda\,|\,\alpha_i)}-1}{q-1} F_iK_\lambda $$
yield
$$ \overline h_{\lambda}\overline E_i-\overline E_i\overline h_{\lambda}=
 (\lambda\,|\,\alpha_i)\overline E_i\ ,~\overline h_{\lambda}\overline F_i-\overline
 F_i\overline h_{\lambda}=-(\lambda\,|\,\alpha_i)\overline F_i\ .$$

Furthermore, the relation (\ref{eq:upper lower relations}) in $\UU_{\bf A}(\gg)$:
$$ E_iF_j-F_jE_i=\delta_{i,j} \frac{K_{\alpha_i}- K_{-\alpha_i}}{q^{d_i}-q^{-d_i}}=
\delta_{i,j}h_{\alpha_i}\frac{1+K_{-\alpha_i}}{1+q^{-d_i}}\cdot \frac{q-1}{q^{d_i}-1}$$
becomes $\overline E_i\overline F_j-\overline F_j\overline E_i=\delta_{i,j}\frac{1}{d_i}\overline h_{\alpha_i}$.

Also the quantum Serre relations (\ref{eq:quantum Serre relations})
in $\UU_{\bf A}(\gg)$  turn into the classical Serre relations in
$\gg$. Therefore, if we denote $\overline H_i=\frac{1}{d_i}\overline
h_{\alpha_i}$, then the elements $\overline E_i,\overline
F_i,\overline H_i$ become the standard Chevalley generators of the
semisimple part of $\gg$.

This proves that the Hopf algebra $\overline {\UU_{q}(\gg)}$ is naturally isomorphic to the universal enveloping algebra $U(\gg)$.
 The lemma is proved.
\endproof

Therefore, based on Lemmas \ref{le:Hopf modules}, \ref{le:Hopf
pair}, and   \ref{le:U(g)-action} we define the ''dequantization''
functor $\UU_q(\gg)-Mod\to U(\gg)-Mod$.

Denote by $\OO_{f}(\UU_{q}(\gg))$ the full sub-category of
$\UU_{q}(\gg)-Mod$ which objects $(V,L)$ are such that $V$ is an
object of $\OO_f$ and $L$ is an ${\bf A}$-lattice in $V$ compatible
with the weight decomposition $V=\oplus_{\mu\in P} V(\mu)$, i.e.,
each $L(\mu):=V(\mu)\cap L$ is an ${\bf A}$-lattice in $V(\mu)$.

Let  $V_\lambda\in Ob(\OO_f)$ be an irreducible $U_{q}(\gg)$-module
with highest weight $\lambda\in P^+$ and let $v_{\lambda}\in
V_{\lambda}$ be a highest weight vector.  Define $L_{v_{\lambda}}=
U_{{\bf A}} (\gg) \cdot v_{\lambda}$.

\begin{lemma} $(V_{\lambda},L_{v_{\lambda}})\in \OO_{f}(\UU_{q}(\gg))$.

\end{lemma}
 \proof
Clearly, ${\bf k} \otimes_{{\bf A}} L_{v_{\lambda}}= V_{\lambda}$.
It remains to  show that $L_{v_{\lambda}}$ is a free ${\bf
A}$-module. It is easy to see that  $U_{\bf
A}(\gg)(v_\lambda)=U^-_{\bf A}(v_\lambda)$, where $U^-_{\bf A}$ is
the ${\bf A}$-subalgebra of $U_{\bf A}(\gg)$ generated by
$F_1,F_2,\ldots,F_r$. Clearly, $U^-_{\bf A}$ is graded via
$deg(F_i)=1$: $U^-_{\bf A}=\oplus_{n\ge 0} (U^-_{\bf A})_n$, where
each $(U^-_{\bf A})_n$ is a finitely generated ${\bf A}$-module.
Clearly, $(U^-_{\bf A})_n(v_\lambda)=0$ for $n>>0$. Therefore,
$L_{v_{\lambda}}$ is a finitely generated ${\bf A}$-module of finite
rank. Finally, since ${\bf A}$ is a PID, this implies that
$L_{v_{\lambda}}$ is a free ${\bf A}$-module. The Lemma is proved.
 \endproof

The following fact is obvious and, apparently, well-known.

 \begin{lemma}
 \label{le: equiv of OOs}

(a) Each object $(V_{\lambda},L_{v_{\lambda}})$ is irreducible in
$\OO_{f}(\UU_{q}(\gg))$; and each irreducible object of
$\OO_{f}(\UU_{q}(\gg))$ is isomorphic to one of
$(V_{\lambda},L_{v_{\lambda}})$.

\noindent (b)  The category $\OO_{f}(\UU_{q}(\gg))$  is semisimple.

\noindent (c) The forgetful functor  $(V ,L)\mapsto V$ is an almost equivalence of tensor categories $\OO_{f}(\UU_{q}(\gg))\to \OO_{f}$.
 \end{lemma}

The following fact shows that the category $\OO_{f}(\UU_{q}(\gg))$ is naturally braided.

\begin{lemma}
\label{le:kA-braiding}
Let $\VV=(V,L)$ and $\VV'=(V',L')$ be objects of $\OO_{f}(\UU_{q}(\gg))$. Then:

\noindent (a) The braiding operator $\RR_{V,V'}:V\otimes_{\bf k} V'\to V'\otimes_{\bf k} V$ satisfies:
$$\RR_{V,V'}(L\otimes_{\bf A} L')=L'\otimes_{\bf A} L \ .$$

\noindent (b) The operator $\RR_{V,V'}$ extends to an operator:
$$\RR_{\VV,\VV'}:\VV\otimes  \VV'\to \VV'\otimes  \VV$$
and the collection of all  $\RR_{\VV,\VV'}$ constitutes a braiding on $\OO_{f}(\UU_{q}(\gg))$.

\noindent (c) The forgetful functor  $(V ,L)\mapsto V$ is an almost equivalence of braided tensor categories $\OO_{f}(\UU_{q}(\gg))\to \OO_{f}$.
\end{lemma}

\proof Note that the universal $R$-matrix $R\in
U_q(\gg)\widehat\otimes_{\bf k} U_q(\gg)$ belongs, in fact, to
$U_{{\bf A}}(\gg)\widehat\otimes_{\bf A} U_{{\bf A}}(\gg)$. Recall
also that $R$ acts on $V\otimes_{\bf k} V'$. Therefore,
$R(L\otimes_{\bf A} L')=L\otimes_{\bf A} L'$ and
$\RR_{V,V'}(L\otimes_{\bf A} L')=\tau R(L\otimes_{\bf A}
L')=\tau(L\otimes_{\bf A} L')=L'\otimes_{\bf A} L$. This proves (a).

Parts (b) and (c) also follow.  \endproof

Denote by $\overline \OO_f$ the full sub-category of
$U(\gg)-Mod$, which objects $\overline V$ are finite-dimensional
$U(\gg)$-modules having a weight decomposition $V=\oplus_{\mu\in P}
V(\mu)$.

\begin{lemma}
\label{le:dequantization equivalence} (a) The restriction of the
functor $\UU_q(\gg)-Mod\to U(\gg)-Mod$ defined by
(\ref{eq:dequantization functor}) to the sub-category
$\OO_{f}(\UU_{q}(\gg))$ is a tensor functor
\begin{equation}
\label{eq:dequantization q=1}
\OO_{f}(\UU_{q}(\gg))\to \overline \OO_{f} \ .
\end{equation}

\noindent (b) The functor (\ref{eq:dequantization q=1}) is an almost equivalence of  categories.

\end{lemma}

\proof Prove (a). It suffices to show that for each object $(V,L)$
in $\OO_{f}(\UU_{q}(\gg))$ the corresponding $U(\gg)$-module
$\overline V=L/\mm L$ has a weight decomposition. By definition of
$\OO_{f}(\UU_{q}(\gg))$, the ${\bf A}$-lattice $L$ has a weight
decomposition $L=\oplus_{\mu\in P} L(\mu)$ and for any $\mu,\nu\in
P$ and any weight vector $v\in L(\mu)$ one has
$$  h_{\nu} v = \frac{q^{(\mu,\nu)}-1}{q-1}  v_{\lambda}\ .$$
This implies that passing to the quotient by $(q-1)L$, we obtain
$$ \overline h_{\nu} (\overline v)=(\mu,\nu) \overline v \ .$$
That is, the $U(\gg)$-module $\overline V$ has a weight
decomposition $\overline V=\oplus_{\mu\in P} \overline V(\mu)$. This
proves (a). Prove (b) now. Since both categories
$\OO_{f}(\UU_{q}(\gg))$ and $\overline \OO_{f}$ are semisimple, it
suffices to compute the functor (\ref{eq:dequantization q=1})  on
the irreducible objects. In fact, it suffices to prove that
$\overline V_\lambda:=L_{v_\lambda}/(q-1)L_{v_\lambda}$   is
isomorphic to  corresponding irreducible $U(\gg)$-module with the
highest weight $\lambda$.

 Furthermore, if we denote
$\overline v_\lambda=v_\lambda + (q-1)L_{v_{\lambda}}$, and for any
$u\in U_{\bf A}(\gg)$ we denote $\overline u=u+(q-1)U_{\bf A}(\gg)$,
then we obtain
$$\overline u(\overline v_\lambda)=(u+(q-1)U_{\bf A}(\gg))(v_\lambda + (q-1)L_{v_\lambda})=u(v_\lambda) + (q-1)L_{v_\lambda} \ .$$

This proves that $\overline V_\lambda$ is a cyclic
finite-dimensional $U(\gg)$-module with the highest weight
$\lambda$. Therefore, $\overline V_\lambda$ is irreducible. The
lemma is proved. \endproof

Combining Lemmas \ref{le: equiv of OOs} and \ref{le:dequantization equivalence}, we obtain an obvious (and well-known) corollary.

\begin{corollary}
\label{cor:almost equivalent O bar O}
The categories  $\OO_f$ and  $\overline \OO_f$ are almost equivalent.
\end{corollary}

For any object $V$ of  $\OO_f$ written as $V= \bigoplus_{i=1}^{n} V_{\lambda_i}$ we
will refer to the object $\overline V = \bigoplus_{i=1}^{n} \overline V_{\lambda_i}$ in $\overline \OO_f$ as the {\it classical limit} of $V$.
Denote by  $\OO_{gr,f}(\UU_q(\gg))$ the  sub-category of $\UU_q(\gg)-Mod$ whose objects are $\ZZ_{\ge 0}$-graded:
$$V=\bigoplus_{n\in \ZZ_{\ge 0}} V_n$$
where each $V_n$ is an object of $\OO_f$; and  morphisms are those homomorphisms of  $\UU_q(\gg)$-modules which preserve the $\ZZ_{\ge 0}$-grading.
Similarly,  denote by  $\overline \OO_{gr,f}$ the sub-category of $U(\gg)-Mod$ whose objects are $\ZZ_{\ge 0}$-graded:
$$\overline V=\bigoplus_{n\in \ZZ_{\ge 0}} \overline V_n$$
where each $\overline V_n$ is an object of $\OO_f$; and  morphisms are those homomorphisms of $U(\gg)$-modules which preserve the $\ZZ_{\ge 0}$-grading.

 \begin{lemma}
The functor (\ref{eq:dequantization q=1}) extends to a an almost
equivalence of tensor categories $\OO_{gr,f}(\UU_q(\gg))\to
\overline \OO_{gr,f}$ via
$$V=\bigoplus_{n\in \ZZ_{\ge 0}} V_n\mapsto  \overline V=\bigoplus_{n\in \ZZ_{\ge 0}} \overline V_n \ .$$
 \end{lemma}
We will sometimes call $\overline V=\bigoplus_{n\in \ZZ_{\ge 0}}
\overline V_n$ the {\it classical limit} of $V=\bigoplus_{n\in
\ZZ_{\ge 0}} V_n$.

Let $A=A_{\overline 0}\oplus A_{\overline 1}$ be a $\ZZ_2$-graded ${\bf k}$-algebra. Define the {\it super-commutator} $[\cdot,\cdot]$ on $A$ by
$$[a,b]=ab-(-1)^{\varepsilon \delta}ba$$
for any $a\in A_\varepsilon$, $b\in A_\delta$. It is well-known that
the super-commutator  satisfies the super-Leibniz rule and
super-Jacobi identity:
\begin{equation}
\label{eq:super-Leibniz}
[a,bc]=[a,b]c+(-1)^{\varepsilon\cdot \delta}b[a,c]
\end{equation}
\begin{equation}
\label{eq:Jacobian super-commutator}
(-1)^{\varepsilon\gamma} [a,[b,c]]+  (-1)^{\gamma\delta} [c,[a,b]]   +    (-1)^{\delta\varepsilon}[b,[c,a]]=0
\end{equation}
for any $a\in A_\varepsilon, b\in A_\delta, c\in A_\gamma$.

Let now $L$ be an ${\bf A}$-lattice in $A$ such that $L$ is closed
under multiplication and compatible with the $\ZZ_2$-grading (i.e.,
the pair ${\mathcal A}=(A,L)$ is a $\ZZ_2$-graded $({\bf k}, {\bf
A})$-algebra).

\begin{definition}
\label{def:almost commutative}
We say that a $\ZZ_2$-graded $({\bf k}, {\bf A})$-algebra ${\mathcal
A}=(A,L)$ is an {\it almost commutative} superalgebra if
$[L,L]\subset \mm L$.
\end{definition}

Let ${\VV}={\VV}_{0}\oplus {\VV}_{1} $ be a $\ZZ_{2}$-graded  $({\bf
k}, {\bf A})$-module. Let $\mathcal A$ be a $\ZZ_{2}$-graded $({\bf
k}, {\bf A})$-algebra generated by ${\VV}$ (i.e., there is  a
surjective homomorphism $T(\VV)\twoheadrightarrow \mathcal A$ of
$({\bf k},{\bf A})$-algebras, where $T(\VV)$ is the tensor algebra
of $\VV$). For any  map $f:L\times  L\to L\cdot L$ we define the
${\bf A}$-submodule $I_f$ of $L\cdot L$ by
$$I_f=Span_{\bf A}\{[l,l']+\mm f(l,l')\,|\,l,l'\in L\}  $$
and then denote  ${\mathcal I}_f=({\bf k}\otimes_{\bf A} I_f,I_f)$.
The following fact is obvious.
\begin{proposition}
\label{prop: almost super comm}
Let ${\mathcal A}$ be a $\ZZ_2$-graded $({\bf k}, {\bf A})$-algebra  generated by a $\ZZ_2$-graded $({\bf k}, {\bf A})$-module ${\VV}=(V,L)$ and let ${\mathcal J}$ be any $\ZZ_2$-homogeneous ideal containing ${\mathcal I}_f$.
Then the quotient superalgebra ${\mathcal A}/{\mathcal J}$ is almost commutative.
\end{proposition}

\proof  Let $L^k=L\cdot L\cdots  L$, where the product is taken $k$ times. We will proceed by induction in $k$. If $k=1$ then
$$[a,b]\in {\mathcal I}_f+\mm L^2$$
for $a,b\in L^1=L$, which implies that $[\overline a,\overline b\,]
\equiv 0 \mod \mm$ in the quotient algebra ${\mathcal A}/{\mathcal J}$. 
The inductive step follows directly from the super-Leibniz rule
(\ref{eq:super-Leibniz}). This proves the proposition. \endproof

Note that the functor $\FF$ from Lemma \ref{le:specialization functor} takes a $\ZZ_2$-graded $({\bf k}, {\bf A})$-superalgebra ${\mathcal A}=(A,L_A)$ into  the superalgebra $\FF({\mathcal A})=L_A/\mm L_A$ over the residue field $\tilde {\bf k}={\bf k}/\mm$.
The following fact is obvious. 

\begin{lemma} 
\label{le:commutative}
For any almost commutative $({\bf k}, {\bf A})$-algebra ${\mathcal A}=(A,L_A)$ the 
$\tilde {\bf k}$-superalgebra $\FF({\mathcal A})$ is a commutative superalgebra (see Definitions \ref{def:commutative superalgebra} and \ref{def:almost commutative}).

\end{lemma}

In what follows we assume additionally that  ${\bf A}$ is PID, and we fix an element $h\in {\bf A}$ such that $\mm=h{\bf A}$.

 Let
${\mathcal A}=(A,L_A)$ be an almost  commutative   $({\bf k}, {\bf
A})$-superalgebra. Theidentity $[L_A,L_A]  \subset hL_A$ 
allows for defining a bilinear map $\{\cdot,\cdot\}:\FF({\mathcal A})\times \FF({\mathcal A})\to \FF({\mathcal A})$ via:
\begin{equation}
\label{eq:poisson bracket dequantized}
\{l+hL_A,l'+hL_A\}:=\frac{1}{h}[l,l']+hL_A
\end{equation}
for $l,l'\in L_A$. 

\begin{lemma}
\label{le:poisson specialization} For any almost commutative $({\bf k}, {\bf A})$-algebra ${\mathcal A}$ the operation
\eqref{eq:poisson bracket dequantized} is a
Poisson bracket on the commutative  $\tilde {\bf
k}$-superalgebra $\FF({\mathcal A})$. In the other words,
$\FF({\mathcal A})$ is a Poisson superalgebra.
\end{lemma}

\proof It suffices to verify the super-Leibniz rule (\ref{eq: super
Leibniz rule}) and super-Jacobi identity (\ref{eq:Jacobian}). For
any $a\in A_\varepsilon, b\in A_{\delta}, c\in A_\gamma$ we
abbreviate $\overline a= a+\mm L_{A},\overline b= b+\mm L_{A},
\overline c= c+\mm L_{A}$.   Then
$$\{\overline a,\overline b\overline c \}- \{\overline a,\overline b\}
\overline c-(-1)^{\varepsilon\delta}\overline b\{\overline a,\overline c\}= \frac{1}{h}
\left([a,bc]-[a,b]c-(-1)^{\varepsilon\cdot \delta}b[a,c] \right)+\mm L_{A}= \mm L_{A}$$
by (\ref{eq:super-Leibniz}).
$$J(\overline a,\overline b,\overline c)=(-1)^{\varepsilon\gamma}
\{\overline a,\{\overline b,\overline c\}\} +(-1)^{\gamma\delta}
\{\overline c,\{\overline a,\overline b\}\}+(-1)^{\delta\varepsilon}
\{ \overline b,\{\overline c,\overline a\}\}= $$
$$  \frac{1}{h^{2}} \left((-1)^{\varepsilon\gamma} [a,[b,c]]+  (-1)^{\gamma\delta} [c,[a,b]]
+ (-1)^{\delta\varepsilon}[b,[c,a]] \right)+\mm L_{A}= \mm L_{A}$$
by (\ref{eq:Jacobian super-commutator}). The lemma is proved.
\endproof

%%%%%%%%%%%%%%%%%%%%%%%%%%%%%%%
Let ${\VV}=(V,L)$ be an object in $({\bf k},{\bf A})-Mod$ and let
$\sigma:{\VV}\otimes {\VV}\to {\VV}\otimes {\VV}$ be an involution
in this category. Similarly to Definition \ref{def:symmetric power}
above we define the symmetric and exterior squares of ${\VV}=(V,L)$
by
$$S_\sigma ^2{\VV}=Ker (\sigma-id)=Im(\sigma+id),\Lambda^2{\VV}=Ker (\sigma+id)=Im(\sigma-id)\ .$$
Then we define the braided symmetric and exterior algebras.
\begin{equation}
\label{eq:general symmetric and exterior algebras}
S_\sigma ({\VV})=T({\VV})/\left< \Lambda_\sigma^{2}{\VV}\right>,~\Lambda_\sigma({\VV})=T({\VV})/ \left<S_\sigma^{2}{\VV}\right> \ .
\end{equation}

Clearly, both $({\bf k},{\bf A})$-algebras $S_\sigma ({\VV})$ and
$\Lambda_\sigma ({\VV})$ are naturally $\ZZ_2$-graded via:
$$S_\sigma ({\VV})_{\overline 0}=S_\sigma ({\VV}),~S_\sigma ({\VV})_{\overline 1}=\{0\},
~\Lambda_\sigma ({\VV})_{\overline 0}=\bigoplus_{n\ge 0}
\Lambda_\sigma ({\VV})_{2n},~\Lambda_\sigma ({\VV})_{\overline
1}=\bigoplus_{n\ge 0} \Lambda_\sigma ({\VV})_{2n+1} \ .$$ Next, we
will describe all involutions $\sigma$ for which  both algebras
$S_\sigma ({\VV})$ and $\Lambda_\sigma ({\VV})$ almost
commutative.

\begin{definition}
\label{def:almost permutation} Given an object ${\VV}=(V,L)$ in
$({\bf k},{\bf A})-Mod$ and an involution $\sigma:{\VV}\otimes
{\VV}\to {\VV}\otimes {\VV}$. We say that $\sigma$ is {\it almost
permutation} if $(\sigma-\tau)(L\otimes_{\bf A}  L)\subset \mm
L\otimes_{\bf A}  L$, where   $\tau:V\otimes V\to V\otimes V$ is the
permutation of factors.
\end{definition}

\begin{lemma}
\label{le:super-commutative classical limit}
Let ${\VV}$ be an object of $({\bf k},{\bf A})-Mod$  and let $\sigma:{\VV}\otimes {\VV}\to {\VV}\otimes {\VV}$ be an almost permutation. Then:

\noindent (a) The $({\bf k},{\bf A})$-algebras $S_\sigma ({\VV})$ and $\Lambda_\sigma ({\VV})$ are almost commutative.

\noindent (b) The $\tilde {\bf k}$-algebras $\FF(S_\sigma (\VV))$ and $\FF(\Lambda_\sigma(\VV))$ are homogeneous super-Poisson.

\end{lemma}

\proof Prove (a). First, show that $S_\sigma ({\VV})$ is almost commutative. Let us introduce a $\ZZ_2$-grading on the tensor
algebra $T(\VV)$ by setting $T(\VV)_{\overline 0}=T(\VV)$ and
$T(\VV)_{\overline 1}=\{0\}$.  Clearly, the ideal
$\left<\Lambda_\sigma^2\VV\right>$ is $\ZZ_2$-graded. Moreover, this
ideal satisfies the assumptions of Proposition \ref{prop: almost
super comm} with $f(l,l')=\frac{1}{h}(\tau-\sigma)(l\otimes l')$
because
\begin{equation}
\label{eq: f-symmetric}
l\otimes l'-\sigma(l\otimes l') =[l,l']+hf(l,l')
\end{equation}
for each $l,l'\in L$. Therefore, Proposition \ref{prop: almost super
comm} guarantees that the quotient algebra of $\mathcal A=T(\VV)$ by
the ideal $\mathcal J=\left<\Lambda_\sigma^{2} \VV\right>$ is
almost commutative.

Next, we show that $\Lambda_\sigma ({\VV})$ is almost commutative. Let us introduce a $\ZZ_2$-grading on the tensor
algebra $T(\VV)$ by first turning ${\VV}$ into $\ZZ_2$-graded space
via  ${\VV}_{\overline 0}=\{0\}$, ${\VV}_{\overline 1}={\VV}$, and
extending this naturally to each tensor power of ${\VV}$.   Clearly,
the ideal $\left< S_\sigma^{2} \VV \right>$ is $\ZZ_2$-graded.
Moreover, this ideal satisfies the assumptions of Proposition
\ref{prop: almost super comm} with
$f(l,l')=\frac{1}{h}(\sigma-\tau)(l\otimes l')$ because
\begin{equation}
\label{eq: f -exterior}
l\otimes l'+\sigma(l\otimes l')=[l,l']+hf (l,l')
\end{equation}
for each $l,l'\in L$. Therefore, Proposition \ref{prop: almost super
comm} guarantees that the quotient algebra of $\mathcal A=T(\VV)$ by
the ideal $\mathcal J=\left< S_\sigma^{2} \VV\right>$ is almost commutative. Part (a) is proved.

Part (b) follows directly from (a)  and Lemma \ref{le:poisson specialization}.
The lemma is proved.
\endproof

Let ${\VV}=(V,L)$ be an object of $({\bf k},{\bf A})-Mod$ and let
$\sigma:{\VV}\otimes {\VV}\to {\VV}\otimes {\VV}$ be an almost
permutation. Denote $\overline V=\FF(\VV)=L/hL$ and define the map
${\bf r}$ by the formula: ${\bf r}=\FF(\frac{1}{h}(\sigma-\tau))$,
that is,
\begin{equation}
\label{eq:rr}
{\bf r}((l+hL)\otimes (l'+hL'))= \frac{1}{h}(\tau\circ \sigma-id)(l\otimes l')+hL\otimes_{\bf A} L \ .
\end{equation}
for $l,l'\in L$. Clearly, $\tau\circ {\bf r}=-{\bf r}\circ \tau$.
Therefore, we can define homogeneous brackets
$\{\cdot,\cdot\}_+=\{\cdot,\cdot\}_{+,{\bf r}}$ on  $S(\overline V)$
and $\{\cdot,\cdot\}_-=\{\cdot,\cdot\}_{-,{\bf r}}$ on
$\Lambda(\overline V)$ by the formula:
\begin{equation}
\label{eq:positive bracket}
\{\overline v,\overline v'\}_+= \frac{1}{2}(id+\tau){\bf r}(\overline v\otimes \overline v')
=\frac{1}{2}{\bf r}(\overline v\otimes \overline v'-\overline v'\otimes \overline v)
\end{equation}
\begin{equation}
\label{eq:negative bracket} \{\overline v,\overline v'\}_-=
\frac{1}{2}(id-\tau){\bf r}(\overline v\otimes \overline
v')=\frac{1}{2}{\bf r}(\overline v\otimes \overline v'+\overline
v'\otimes \overline v)
\end{equation}
for $v,v'\in \overline V$.

Clearly, $\{\overline l,\overline l'\}_+=-\{\overline l',\overline
l\}_+$ and $\{\overline l,\overline l'\}_-=\{\overline l',\overline
l\}_-$ for any $v,v'\in \overline V$ and the brackets naturally
extend to $S(\overline V)$ on $\Lambda(\overline V)$ by the
(super)Leibniz rule (see e.g., Definition (\ref{def:bracketed
superalgebra})).

%The following fact is obvious.
%\begin{lemma}
%Let ${\VV}$ be an object of $({\bf k},{\bf A})-Mod$  and let $\sigma:{\VV}\otimes {\VV}\to {\VV}\otimes {\VV}$ be an almost permutation. Then $S^{2}(\overline V)$ and $\overline{S^{2}_\sigma V}$, respectively $\Lambda^{2}\overline V$ and  $\overline{\Lambda^{2}_\sigma V}$,  are canonically isomorphic.
%\end{lemma}

%Via the canonical isomorphism, the bracket defined in Equation \ref{eq: poisson bracket} defines a map

%$$\{ \cdot,\cdot\}: S^{1}(\overline V)\times S^{1}(\overline V)\to S^{2}(\overline V)\ ,$$
%$$\{ \cdot,\cdot\}: \Lambda^{1}(\overline V)\times \Lambda^{1}(\overline V)\to \Lambda^{2}(\overline V)\ .$$

Denote by $P_{\bf r}(S(\overline V))$ and $P_{\bf
r}(\Lambda(\overline V))$ the Poisson closures of the bracketed
superalgebras  $(S(\overline V), \{\cdot,\cdot\}_+)$ and $(\Lambda
(\overline V),\{\cdot,\cdot\}_{-}) $ respectively (along the lines
of Proposition \ref{pr:universal poisson closure}).

Recall that the functor $\FF:({\bf k},{\bf A})-Mod\to \tilde {\bf
k}-Mod$ defined above is an almost equivalence of tensor categories and for any $({\bf k},{\bf A})$-algebra ${\mathcal A}$
the image $\FF({\mathcal A})$ is a $\tilde {\bf k}$-algebra.

\begin{theorem}
\label{th:limits projection} For every object $\VV=(V,L)$ of $({\bf k}, {\bf A})-Mod$ and  any almost permutation
$\sigma: \VV\otimes \VV\to  \VV\otimes \VV$ one has a surjective homomorphism of homogeneous super-Poisson $\tilde {\bf k}$-algebras
\begin{equation}
\label{eq:homomorphism Poisson closures}
P_{\bf r}(S(\overline V))\twoheadrightarrow \FF(S_\sigma (\VV)) ,~ P_{\bf r}(\Lambda(\overline V)) \twoheadrightarrow \FF(\Lambda_\sigma(\VV)) \ ,
\end{equation}
where $\overline V=\FF(\VV)=L/hL$.
\end{theorem}

\proof We need the following obvious fact.

\begin{lemma}
\label{le:generated} Let ${\mathcal A}$ be a $({\bf k},{\bf
A})$-algebra generated by $\VV=(V,L)$ (i.e., there is a surjective
homomorphism $T(\VV)\twoheadrightarrow {\mathcal A}$). Then
$\FF({\mathcal A})$ is a $\tilde {\bf k}$-algebra generated by
$\overline V=\FF(\VV)=L/hL$.
\end{lemma}

By Lemma \ref{le:super-commutative classical limit}(a), $S_\sigma
(\VV)$ is an evenly $\ZZ_2$-graded almost commutative $({\bf k},{\bf A})$-algebra generated by
$\VV=(V,L)$. Hence, by Lemmas \ref{le:super-commutative classical
limit}(b) and  \ref{le:generated}, $\FF(S_\sigma (\VV))$ is an evenly
$\ZZ_2$-graded commutative  $\tilde
{\bf k}$-superalgebra generated   by $\overline V$. Therefore, there
exists a surjective homomorphism
$$\varphi_+:S(\overline V)\twoheadrightarrow  \FF(S_\sigma (\VV)) \ .$$
Note that the restriction of $\varphi_+$ to $S^2(\overline V)$ is an isomorphism $S^2(\overline V)\widetilde \to \FF(S_\sigma(\VV))_2$
because
$$\FF(S_\sigma(\VV))_2=\FF((\VV\otimes \VV)/\Lambda^2_\sigma \VV))=(\overline V\otimes \overline V)/\Lambda^2(\overline V) \ .$$
Similarly,  by Lemma \ref{le:super-commutative classical limit}(a),
$\Lambda_\sigma (\VV)$ is a $\ZZ_2$-graded almost commutative
$({\bf k},{\bf A})$-superalgebra generated by $\VV=(V,L)$. Hence, by
Lemmas \ref{le:super-commutative classical limit}(b) and
\ref{le:generated}, $\FF(\Lambda_\sigma (\VV))$ is a $\ZZ_2$-graded
commutative $\tilde {\bf k}$-superalgebra generated   by $\overline V$.
Therefore, there exists a surjective homomorphism
$$\varphi_-:\Lambda(\overline V)\twoheadrightarrow  \FF(\Lambda_\sigma (\VV)) \ .$$
Note that the restriction of $\varphi_+$ to $\Lambda^2(\overline V)$
is an isomorphism $\Lambda^2(\overline V)\widetilde \to
\FF(S_\sigma(\VV))_2$ because
$$\FF(\Lambda_\sigma(\VV))_2=\FF((\VV\otimes \VV)/\Lambda^2_\sigma \VV))=(\overline V\otimes \overline V)/\Lambda^2(\overline V) \ .$$

Clearly, both $\varphi_+$ and $\varphi_-$ preserves the even $\ZZ_2$-grading. We need  the following obvious fact.

\begin{lemma}
\label{le:bracketed homogeneous homomorphisms} Let $A=\oplus_{n\ge
0} A_n$ and $B=\oplus_{n\ge 0} B_n$ be homogeneous bracketed
superalgebras (see Definition \ref{def:homogeneous}). And let
$\varphi:A\to B$ be an algebra homomorphism preserving the $\ZZ_{\ge
0}$-grading and such that
$\varphi(\{a,a'\})=\{\varphi(a),\varphi(a')\}$ for all $a,a'\in
A_1$.

Then $\varphi$ is a homomorphism of bracketed superalgebras.

\end{lemma}

It is easy to see that Poisson bracket on $\FF(S_\sigma(\VV))$ given
by (\ref{eq:poisson bracket dequantized}) satisfies:
$$\{l+hL,l'+hL\}:=\frac{1}{h}(l\cdot l'-l'\cdot l)+hL\cdot L$$
for any $l,l'\in L$. Note that $l\cdot l'=l\otimes
l'+\Lambda^2_\sigma L=\sigma(l\otimes l')+\Lambda^2_\sigma L$ and
$l'\cdot l=l'\otimes l+\Lambda^2_\sigma L=\tau(l\otimes
l')+\Lambda^2_\sigma L$. Therefore,
$$\{l+hL,l'+hL\}={\bf r}((l+hL)\otimes (l'+hL))+ \FF(\Lambda_\sigma ^2 (V,L))$$
Note that $\FF(\Lambda_\sigma ^2 (V,L))=\Lambda^2(\overline V)$
because $\FF(\sigma)=\tau$. On the other hand, 
$$\overline v\otimes \overline v'+\Lambda^2(\overline V)
=\frac{1}{2}\left(id+\tau\right)(\overline v\otimes \overline v')+\Lambda^2(\overline V)$$ 
for any $\overline v,\overline v'\in \overline V$. Therefore,
$$\{\overline v, \overline v'\}+\Lambda^2(\overline V)=\frac{1}{2}\left(id+\tau\right){\bf r}(\overline v\otimes
\overline v')+\Lambda^2(\overline V)=\varphi_+(\{\overline v,\overline v'\}_+) \ ,$$
where $\{\overline v,\overline v'\}_+$ is given by (\ref{eq:positive bracket}).
Similarly, one has in  $\FF(\Lambda_\sigma(\VV))_2$ one has:
$$\{\overline v, \overline v'\}+S^2(\overline V)=\frac{1}{2}\left(id-\tau\right){\bf r}(\overline v\otimes
\overline v')+\Lambda^2(\overline V)=\varphi_-(\{\overline v,\overline v'\}_-) \ ,$$
where $\{\overline v,\overline v'\}_-$ is given by (\ref{eq:negative bracket}).

Hence  $\varphi_+$ and $\varphi_-$ are  homomorphisms of bracketed
superalgebras by Lemma \ref{le:bracketed homogeneous homomorphisms}.
Hence Proposition \ref{pr:universal poisson closure} implies that
the homomorphisms of bracketed algebras $\varphi_{+}$ and
$\varphi_{-}$  factor through the respective Poisson closures
$P_\sigma(S(\overline V))$, $P_\sigma(\Lambda(\overline V))$, and we
obtain (\ref{eq:homomorphism Poisson closures}). Therefore, Theorem
\ref{th:limits projection}  is proved. \endproof

Let now  ${\bf k}=\CC(q)$ and  ${\bf A}$ be the algebra of rational
the functions in $q$ regular at $q=1$ so that $\tilde {\bf  k}=\CC$.
As we argued above, $A$ is a local PID, with maximal ideal
$\mm=h{\bf A}$, where $h=q-1$. Recall that for any objects $U$ and
$V$ of  $\OO_{f}$ we defined in (\ref{eq:sigma}) the normalized
braiding $ \sigma_{U,V}:U\otimes V \to V\otimes U$. In particular,
each $\sigma_{V, V}$ is an involution.

Recall that $\UU_{q}(\gg)=(U_q(\gg),U_{\bf A}(\gg))$ is a Hopf $({\bf k},{\bf A})$-algebra as above.

\begin{lemma}
\label{le:sigma over A} Let $V\in \OO_f$  and let $L$ be a ${\bf
A}$-lattice in $V$ such that $L$ is also a $U_{\bf A}(\gg)$-submodule of
$V$, i.e., $\VV=(V,L)$ is an object of $\OO_f(\UU_{q}(\gg))$. Then

\noindent (a) $\sigma_{V,V}(L\otimes_{\bf A} L)=L\otimes_{\bf A} L$.
Therefore, $\sigma_{V, V}$ extends to an involutive
$\UU_{q}(\gg)$-module homomorphism $\sigma_{\VV,\VV}:\VV\otimes \VV
\to \VV\otimes \VV$.

\noindent (b) $\sigma_{\VV,\VV}$ is an almost permutation (see Definition \ref{def:almost permutation}).

\noindent (c) In the notation (\ref{eq:rr}) we have ${\bf r}=r^-$,
where $r^-=\frac{1}{2}(r-\tau(r))$ is the skew-symmetrized classical $r$-matrix defined in (\ref{eq:rminus}).
\end{lemma}

\proof Prove (a) It suffices to show that $\sigma_{V,V}$ preserves
$L\otimes_{\bf A} L$. Indeed, it follows from the definition
(\ref{eq:D}) that $D\in U_{\bf A}(\gg)\widehat \bigotimes_{\bf A}
U_{\bf A}(\gg)$ and $u\otimes v\mapsto  D(u\otimes v)$ is a
well-defined invertible map $V\otimes_{\bf k} V\to V\otimes_{\bf k}
V$. Therefore, $D(L\otimes_{\bf A} L)=L\otimes_{\bf A} L$. We obtain
by  (\ref{eq:formula for sigma}):
$$\sigma_{V,V}(L\otimes_{\bf A} L)=(\RR_{V,V} D^{-1})(L\otimes_{\bf A} L)=\RR_{V,V}(L\otimes_{\bf A} L)=L\otimes_{\bf A} L$$
by Lemma \ref{le: equiv of OOs}(c). This proves (a).

Prove (b). Indeed, by \eqref{eq:sigma}, $\sigma_{V,V}=\tau RD^{-1}|_{V\otimes V}$. But, taking into account \eqref{eq:R1} and \eqref{eq:hat R}, we see that the operator $RD^{-1}-id$ acts on each $z\in L\otimes_{\bf A} L$ as an element of $(q-1)U_{\bf A}(\gg)$. This implies that $(\sigma_{V,V}-\tau)(L\otimes_{\bf A} L)\subset (q-1)L\otimes_{\bf A} L$. Part (b) is proved.

Prove (c). In order to compute ${\bf r}=\frac{1}{h}\FF(\tau\circ
\sigma-id)$, where $h=q-1$, $\sigma=\sigma_{\VV,\VV}$, recall that
$\sigma=\tau R D^{-1}$ in the notation of Section
\ref{sect:definitions}. Clearly,
$$\FF\left(\frac{1}{h}(R-id)\right)=r, \FF\left(\frac{1}{h}(R^{op}-id)\right)=\tau(r)\ ,$$
$$2\FF\left(\frac{1}{h}(D-id)\right)=\FF\left(\frac{1}{h}(D^2-id)\right)=\FF\left(\frac{1}{h}(R^{op}R-id)\right)=\tau(r)+r \ ,$$
Therefore,
$$\FF\left(\frac{1}{h}(D-id)\right)=\frac{1}{2}\left(r+\tau(r)\right) \ .$$
Finally,
$${\bf r}=\FF\left(\frac{1}{h}(\tau\circ \sigma-id)\right)=\FF\left(\frac{1}{h}(RD^{-1}-id)\right)=r-\frac{1}{2}\left(r+\tau(r)\right)=r^- \ .$$
The lemma is proved.  \endproof

%We need the following fact.
%\begin{lemma}
%The matrix $\widehat R\in U_{q}(\gg)\widehat\otimes U_{q}(\gg) $ defined in  (\ref{eq:hat R})  belongs to $ U_{A}(\gg)\widehat\otimes U_{A}(\gg)$.
%\end{lemma}
%\proof
%Recall that  by (\ref{eq:hat R}) $\sigma_{V,V}(V\otimes V)=\tau \widehat R (V\otimes V)$ for all objects $V$ of $\OO_{f}$. The Lemma follows directly from Lemma \ref{le:sigma over A}. Needs more ??
%\endproof

%Applying the definition  (\ref{eq:hat R}) of $\widehat R$, one obtains
%$$ \widehat R= 1\otimes 1+(q-1)( \tilde R)$$
%for some $ \tilde R\in U_{A}(\gg)\otimes U_{A}(\gg)$. Therefore,
%$$\sigma=\tau \widehat R  $$
% is an almost permutation. The lemma is proved.

Therefore, Theorem \ref{th:limits projection} and Lemma
\ref{le:sigma over A} guarantee that for $\sigma=\sigma_{\VV, \VV}$
both algebras $\FF(S_\sigma(\VV))$ and $\FF(\Lambda_\sigma(\VV)
)$ are super-Poisson and we have surjective homomorphisms
(\ref{eq:homomorphism Poisson closures}). Moreover, Lemma
\ref{le:sigma over A}(c) implies  in  the brackets (\ref{eq:positive
bracket}) and (\ref{eq:negative bracket}) the operator ${\bf r}$ is
now equal to $r^-$. This proves that in the notation of Theorem
\ref{th:flat poisson closure} and \ref{th:limits projection} we see that
$P_{\bf r}(S_\sigma(\overline V))=\overline{S(\overline V)}$,
$P_{\bf r}(\Lambda_\sigma(\overline V))=\overline{\Lambda(\overline
V)}$ and, therefore, obtain the structural surjective homomorphisms of Poisson
superalgebras
\begin{equation}
\label{eq:homomorphism Poisson closures special}
\overline{S(\overline V)}\twoheadrightarrow \FF(S_\sigma (\VV)) ,~ \overline{\Lambda(\overline V)} \twoheadrightarrow \FF(\Lambda_\sigma(\VV)) \ ,
\end{equation}

%%%%%%%%%%%%%%%
This proves Theorem \ref{th:flat poisson closure}. \endproof

\begin{corollary} (from the proof of Theorem \ref{th:flat poisson closure})
\label{cor:flat poisson closure independent} In the notation as
above, for each $\VV=(V,L)$ in $\OO_f(\UU_q(\gg))$ both the algebras
$\FF(S_\sigma (\VV)), \FF(\Lambda_\sigma(\VV))$ and the structure
homomorphisms (\ref{eq:homomorphism Poisson closures special})
depend only on $\overline V$ (but they do not depend on the choice
of $V$ in $\OO_f$ or of the lattice $L$ in $V$).

\end{corollary}

\proof It follows from Lemma \ref{le: equiv of OOs} that both
algebras $\FF(S_\sigma(\VV))$ and $\FF(\Lambda_\sigma(\VV) )$
depend only on $V$, or, more precisely, on the isomorphism class of
$V$ in $\OO_f$. Therefore, by Corollary \ref{cor:almost equivalent O
bar O}, the above algebras depend (up to an isomorphism) only on the
choice of $\overline V$ in $\overline \OO$. By the construction, the
structure homomorphisms (\ref{eq:homomorphism Poisson closures})
also depend only on $\overline V$. Corollary  \ref{cor:flat poisson
closure independent}  is proved.
\endproof

\subsection{Exterior valuations and the proof of Theorem \ref{th:exterior fourth power}} Due to
Proposition \ref{pr:algebra coalgebra duality}(a), the assertion
$\Lambda_\sigma^n V_\ell =0$ for $n\ge 4$ is equivalent to
$\Lambda_\sigma (V_\ell)_n=0$ for $n\ge 4$ which, in its turn,  is
equivalent to  $\Lambda_\sigma(V_\ell)_4 =0$. By Theorem
\ref{th:flat poisson closure},  it suffices to prove that
$\overline{\Lambda (\overline V_\ell)}_4=0$, where $\overline V_\ell$ is the corresponding $\ell+1$-dimensional
$gl_2(\CC)$-module, i.e., $\overline V_\ell$ is the $\ell$-th
symmetric power of the standard $gl_2(\CC)$-module $\overline V_1=\CC^2$. I
n its turn, by definition (\ref{eq:barSV
barLambdaV}) of $\overline{\Lambda (\overline V_\ell)}$, the
latter statement equivalent to following one.

\begin{proposition}
\label{pr:4th exterior power} For each $\ell\ge 0$ we have
$$J_-(\overline V_\ell\times \overline V_\ell\times \overline V_\ell)\wedge \overline V_\ell=\Lambda^4\overline V_\ell\ ,$$
where $J_-:\overline V_\ell\times \overline V_\ell\times \overline V_\ell\to \Lambda^3(\overline V_\ell)$ is
the super-Jacobian map defined in (\ref{eq:Jacobian}).
\end{proposition}

\proof Recall from Example \ref{ex:r-gl2} that the super-Poisson bracket on $\Lambda(\overline V_\ell)$ is given by:
$$\{a,b\}_-=E(a)\wedge F(b)-F(a)\wedge E(b) \ , $$
for any $a,b\in  \Lambda(\overline V_\ell)$.

Recall  that the $gl_2(\CC)$-module $\overline V_\ell$ has a weight basis $\{ v_i , i=0,\ldots,  \ell \}$  such that
 $$E(v_i)=i v_{ i-1},~F(v_i)=(\ell -i)v_{ i+1}\ .$$
This implies that
$$\{v_{ i},v_{ j}\}_-= i(\ell-j)v_{ i-1}\wedge v_{ j+1}-j(\ell-i)v_{ i+1} \wedge v_{j-1}\ .$$

We need the following fact.
\begin{lemma}
\label{le: jacobean explicit}
For any basis vectors $ v_i, v_j, v_{k}\in  V_\ell$ one has
$$J_-( v_i, v_j, v_{k})= i (2j-\ell)(\ell-k) v_{i-1}\wedge v_j\wedge v_{k+1}-i  (\ell-j)(2k-\ell)  v_{i-1}\wedge v_{j+1}\wedge v_{k}$$
$$-(2i-\ell)j(\ell-k)  v_i\wedge v_{j-1}\wedge   v_{k+1}+(2i-\ell)(\ell-j)k  v_i\wedge v_{j+1}\wedge v_{k-1}$$
$$+(\ell-i) j (2k-\ell) v_{i+1}\wedge v_{j-1}\wedge v_{k} -(\ell-i)(2j-\ell)k v_{i+1}\wedge v_j\wedge v_{k-1} \ .$$
\end{lemma}

%The following fact is easy to check and very useful.
%\begin{lemma}
%\label{le:squarefree}
%The map $j_{-}$ is square-free as a bracket on a commutative superalgebra, i.e.
%$$j_{-}(v_i,v_i,v_{k})=j_{-}(v_i,v_j,v_j)=j_{-}(v_i,v_j,v_i)=0\ .$$
%\end{lemma}

For each subset  $I$ of $[0,\ell]=\{0,1,\ldots,\ell\}$ of the form
$I=\{i_1<i_2<\cdots <i_n\}$ denote by $v_I\in \Lambda^n \overline V_\ell$ the monomial $v_I:=v_{i_1}\wedge v_{i_2}\wedge \cdots
\wedge v_{i_n}$. Clearly for each $n$ the monomials $v_I$, $|I|=n$
form a $\CC$-linear  basis  in $\Lambda^n \overline V_\ell$.

Denote by ${\VV}_n$ the set of all $n$-element subsets in
$[0,\ell]$. We define a total ordering on ${\VV}_n$ to be the
restriction of the lexicographic one from the set $[0,\ell]^n$. For
each $n\in [0,\ell]$  define a map $\nu=\nu_n:\Lambda^n \overline V_\ell \setminus \{0\}\to {\VV}_n$ uniquely by the formula:
$$\nu(v_I)=I$$
for any  $0 \le i_1<i_2<\cdots <i_n\le \ell$ and
$\nu(c_1u+c_2v)=\min(\nu(u),\nu(v))$ for any $c_1,c_2\in \CC^\times$
and $u,v\in \Lambda^n \overline V_\ell \setminus \{0\}$ such that
$\nu(u)\ne \nu(v)$ (the minimum is taken in the above ordering of
${\VV}_n$). We will refer to $\nu$ as a {\it valuation} on
$\Lambda^n \overline V_\ell$.  By definition,
${\VV}_n=\nu(\Lambda^n \overline V_\ell \setminus \{0\}$).

The following facts are obvious.

\begin{lemma} If $x\in \Lambda^n  \overline V_\ell$ and $y\in \Lambda^m  \overline V_\ell$ are such that $\nu(x)\cap \nu(y)=\emptyset$,
then $x\wedge y\ne 0$ and $\nu(x\wedge y)=\nu(x) \cup \nu(y)$.

\end{lemma}

\begin{lemma} For any $X\subset  \Lambda^n \overline V_\ell \setminus \{0\}$ such that $\nu(X)={\VV}_n$ the set $X$ spans
$\Lambda^n \overline V_\ell$ as a $\CC$-vector space.

\end{lemma}

Therefore, in order to finish the proof of Proposition \ref{pr:4th
exterior power} it suffices to construct a subset $X$ of
$V_\ell\wedge J_-(\overline V_\ell\times \overline V_\ell\times \overline V_\ell)\setminus \{0\}$ inside
$\Lambda^4 \overline V_\ell \setminus \{0\}$ such that
$\nu(X)={\VV}_4$. We set  $X:=X_{1}\cup X_{y}\cup X_{z}$, where
$X_{1}$, $X_{y}$, and $X_z$ are defined as follows.

\noindent $\bullet$ $X_1=\{x_{i,j,k,m}, 0\le i<j<k\le \ell:\nu(x_{i,j,k,m})=\{i-1,j,k+1,m\}\}$, where
 $$x_{i,j,k,m}:=J_-(v_i,v_j,v_k)\wedge v_m \ .$$

\noindent $\bullet$  $X_{y} =\{ y_{i,j,k}: 0\le i<j<k\le \ell\}$, where
$$y_{i,j,k}=i(2j-2-\ell)  x_{i,j,k,i-1}-   j(2i-\ell)  x_{i,j-1,k,i}\ .$$

\noindent $\bullet$ $X_{z} =\{z_{i,j,k}: 0\le i<j<k\le \ell\}$, where
$$z_{i,j,k}= (2j+2-\ell) (\ell-k) x_{i,j,k,k+1}- (\ell-j)(2k-\ell)  x_{i,j+1,k,k} \ .$$

Now we obtain the desirable fact.
\begin{lemma} In the notation as above we have $\nu(X)=\VV_4$.
\end{lemma}
\proof Clearly, $\nu(X)\supseteq \nu(X_{1})\cup \nu(X_{y})\cup \nu(X_{z})$.
Note that for $i'<j'<k'<m'$ one has:
$$\nu( x_{i'+1,j',m'+1,k'})=\{ i',j',k',m'\}$$
if $i'+1<j'$, $2j\ne \ell$, and
$$\nu( x_{i'+1,k',m'-1,j'} )=\{ i',j',k',m'\} $$
if $k'+1<m'$, $2k\ne \ell$.

Therefore,  $\nu(X_{1})$ contains the following two subsets
$$ \{ \{ i',j',k',m'\}: i'<j'<k'<m',  \text{~$2j'\ne \ell\ $, $i'+1<j'$}\}$$
 and
 $$\{ \{ i',j',k',m'\}: i'<j'<k'<m', \text{~$2k'\ne\ell$, $k'+1< m'$ }\}\ .$$

Furthermore, $ x_{i,j,k,i-1}=  (2i-\ell)j(\ell-k) v_{i-1}\wedge v_i\wedge v_{j-1}\wedge v_{k+1}+lower~terms$
and
$$y_{i,j,k}=\begin{cases}
  (2i-\ell)j  i (\ell-j) (2k-\ell)  v_{i-1}\wedge v_i\wedge v_j\wedge  v_{k}+\text{ lower terms} &\text{ if $j<k-1$}\\
  \delta_{i,j,k} \     v_{i-1}\wedge v_i\wedge v_{k-1}\wedge v_{k}+\text{ lower terms} &\text{ if $j=k-1$}
\end{cases}\ ,
$$
where $  \delta_{i,j,k}=- i(2i-\ell)(\ell-k+1) \left( (2k-4-\ell)  k
+ (k-1) (2k-\ell)\right)$. It is easy to see that $\delta_{i,j,k}=
0$ if and only if $2i=\ell$ because the equation
$$  (2k-4-\ell)  k + (k-1) (2k-\ell)=(2k-1)(2k-\ell-2)-2=0$$
has no integer solutions $k$ for $\ell>0$.
Therefore,
$$\nu (X_{y})\supset \VV''_4= \{\{i',j',k',m'\}: i'+1=j', 2j',2m'\ne \ell \}\ .$$
Finally,
$x_{i,j,k,k+1}=       -i  (\ell-j)(2k-\ell)  v_{i-1}\wedge v_{j+1}\wedge v_{k}\wedge v_{k+1}+lower~terms $
and
$z_{i,j,k}=-(\ell-j)(2k-\ell) (2i-\ell)(j+1)(\ell-k)\ v_i\wedge v_j\wedge v_{k}\wedge v_{k+1}+lower~terms$.
Therefore,  $X_{z}\supset \VV'''_4=\{ \{i',j',k',m'\}:  2k',2i'\ne \ell, m'=k'+1 \}$.

This proves the lemma because  $\VV'_4\cup\VV''_4\cup \VV'''_4=\VV_4$. \endproof

Therefore, Proposition \ref{pr:4th exterior power}  is proved. \endproof

Therefore, Theorem \ref{th:exterior fourth power}  is proved. \endproof

\subsection{Braided triple products and the proof of Theorem \ref{th:exterior powers and symmetric cube}} The proof will be rather long and technical. The strategy is is as follows. First, we will generalize the relevant submodules of $V_\ell\otimes V_\ell\otimes V_\ell$ to {\it braided triple products} inside triple tensor products $V_{\beta_1}\otimes V_{\beta_2}\otimes V_{\beta_3}$ of $U_q(gl_2(\CC))$-modules and  will formulate the appropriate generalization (Theorem \ref{th:exterior and symmetric  products}) of Theorem \ref{th:exterior powers and symmetric cube}. The key ingredient in the proof of Theorem \ref{th:exterior and symmetric  products} is Howe duality (cited below in Proposition \ref{pr: Howe-Duality}) between the highest weight vectors in $V_{\beta_1}\otimes V_{\beta_2}\otimes V_{\beta_3}$ and weight spaces in $U_q(gl_3(\CC))$-modules. Using the Howe duality, we describe the highest weight vectors of braided triple products as certain subspaces inside these $U_q(gl_3(\CC))$-modules  and then, in Theorems \ref{th:intersections} and \ref{th:generic bases in general} establish transversality of these subspaces. Later on, in Proposition \ref{pr:generic bases any beta}, we establish the real reason for this transversality - the presence of the {\it dual canonical basis} in the $U_q(gl_3(\CC))$-modules. The most technical part of the section is the proof of Theorem  \ref{th:intersections}, which  is based on a surprising result in combinatorial optimization (Proposition \ref{pr:convex}) that computes the absolute maximum of a certain linear functional on a discrete convex set.

In order to define {\it braided triple products}, let us write the decomposition \eqref{eq:isotypic decomposition multiplicity free} from Appendix for $U_q(gl_2(\CC))$-modules $V_{\ell_1}$, $V_{\ell_2}$:
$$V_{\ell_1}\otimes V_{\ell_2}=\bigoplus_{0\le m\le \min(\ell_1,\ell_2)} J_{\ell_1,\ell_2}^m \ , $$
where $J_{\ell_1,\ell_2}^m\cong V_{(\ell_1+\ell_2-m,m)}$.
For $\varepsilon\in \{-,+\}$ we define the submodule $V_{\ell_1}\bullet_\varepsilon  V_{\ell_2}$ of $V_{\ell_1}\otimes  V_{\ell_2}$ by:
\begin{equation}
\label{eq:braided product}
V_{\ell_1}\bullet_\varepsilon  V_{\ell_2}:=
\begin{cases} 
\bigoplus\limits_{0\le m\le \min(\ell_1,\ell_2),m\in 2\ZZ} J_{\ell_1,\ell_2}^m & \text{if $\varepsilon=+$} \\
\bigoplus\limits_{0\le m\le \min(\ell_1,\ell_2),m\in 2\ZZ+1} J_{\ell_1,\ell_2}^m & \text{if $\varepsilon=-$} \\
\end{cases}.
\end{equation}
Clearly, $V_\ell\bullet_+  V_\ell=S^2_\sigma V_\ell$ and $V_\ell\bullet_-  V_\ell=\Lambda^2_\sigma V_\ell$ .

Then for any
$\beta_1,\beta_2,\beta_3\in \ZZ_{\ge 0}$ and $\varepsilon_1,\varepsilon_2\in \{-,+\}$ define the {\it braided triple product} $V_{\beta_1}\bullet_{\varepsilon_1}  V_{\beta_2} \bullet_{\varepsilon_2} V_{\beta_3}\subset V_{\beta_1}\otimes  V_{\beta_2} \otimes  V_{\beta_3}$ by 
$$ V_{\beta_1}\bullet_{\varepsilon_1}  V_{\beta_2} \bullet_{\varepsilon_2}  V_{\beta_3}:=V_{\beta_1}\bullet_{\varepsilon_1}  V_{\beta_2} \otimes  V_{\beta_3}\cap V_{\beta_1}\otimes V_{\beta_2} \bullet_{\varepsilon_2}  V_{\beta_3} \ .$$
Clearly, $V_\beta \bullet_+  V_\beta \bullet_+  V_\beta=S^3_\sigma V_\beta$ and $V_\beta \bullet_-  V_\beta \bullet_-  V_\beta=\Lambda^3_\sigma V_\beta$.

\begin{theorem}
\label{th:exterior and symmetric  products}
For any $\beta_1,\beta_2,\beta_3\in \ZZ_{\ge 0}$ and $\varepsilon\in \{-,+\}$ we have an isomorphism of $U_q(gl_2(\CC))$-modules:
$$V_{\beta_1}\bullet_\varepsilon  V_{\beta_2} \bullet_\varepsilon  V_{\beta_3}\cong \oplus_\lambda V_\lambda \ ,$$
where the summation is over all $\lambda=(\lambda_1\ge \lambda_2)$ such that:

\noindent $\bullet$ $\lambda_1+\lambda_2= \beta_1+\beta_2+\beta_3$, $\min(\lambda_2,\beta_2)\ge  (\lambda_2-\beta_1)_++(\lambda_2-\beta_3)_+$,

\noindent $\bullet$ $\min(\lambda_2,\beta_2)$ is even  and $(-1)^{(\lambda_2-\beta_1)_+}=(-1)^{(\lambda_2-\beta_3)_+}=\varepsilon$.

\noindent (with the convention $(a)_+:=\max(a,0)$ for $a\in \ZZ$).
\end{theorem}

\proof  The following fact is well-known.

\begin{lemma} For any $\lambda=(\lambda_1\ge \lambda_2\ge \lambda_3)$
the restriction of each $U_q(gl_3(\CC))$-module $V_\lambda$ to $U_q(gl_2(\CC))_i$, $i=1,2$ is multiplicity-free:
\begin{equation}
\label{eq:Gelfand-Tsetlin decomposition}
V_\lambda|_{U_q(gl_2(\CC))_i}\cong \oplus_\mu  V_\mu^i \ ,
\end{equation}
where each $V_\mu^i$ is the irreducible  $U_q(gl_2(\CC)_i$-submodule
the of the highest weigh $\mu=(\mu_1,\mu_2)$, and the summation is
over all $\mu=(\mu_1,\mu_2)\in \ZZ^2$ such that
\begin{equation}
\label{eq:GT inequality}
\lambda_1\ge \mu_1\ge \lambda_2\ge \mu_2\ge \lambda_3
\end{equation}
\end{lemma}

Therefore, for each $i\in \{1,2\}$ there exists a unique (up to
scalar multiplication)  basis $\BB^i_\lambda$  for $V_\lambda$ such
that   the intersection  $\BB^i_\lambda\cap V_\mu^i$ is the
(canonical) basis for $V_\mu$. Each $b^i\in \BB^i_\lambda$ is
labeled by  a unique {\it Gelfand-Tsetlin pattern} $(\mu,\nu)$ where
$\mu=(\mu_1,\mu_2)$ is such that $b\in V_\mu^i$ and $\nu\in
[\mu_2,\mu_1]$ is such that
\begin{equation}
\label{eq:Gelfand-Tsetlin basis}
b^i_{(\mu,\nu)}=\begin{cases}
F_1^{\mu_1-\nu}(v_\mu^1) & \text{if $i=1$}\\
F_2^{\nu-\mu_2}(v_\mu^2) & \text{if $i=2$}\\
\end{cases},
\end{equation}
where $v_\mu^i$ is the highest weight vector of $V_\mu^i$.

We will refer to each $\BB^i_\lambda$ as a {\it Gelfand-Tsetlin
basis} for $V_\lambda$. For $\varepsilon\in \{-,+\}$ and $i\in
\{1,2\}$ denote by $\BB_\lambda^i(\beta)^\varepsilon$ the set of all
$b^i_{(\mu,\nu)}\in \BB^i_\lambda(\beta)$ such that $(-1)^{\mu_2}=
\varepsilon$.

By the construction, $b_{(\mu,\nu)}^1$ belongs to
$V_\lambda(\beta)$, where
$\beta_1+\beta_2+\beta_3=\lambda_1+\lambda_2+\lambda_3$,
$\mu_1+\mu_2=\beta_1+\beta_2$, $\nu=\beta_1$ and $b_{(\mu,\nu)}^2$
belongs to $V_\lambda(\beta)$, where
$\beta_1+\beta_2+\beta_3=\lambda_1+\lambda_2+\lambda_3$,
$\mu_1+\mu_2=\beta_2+\beta_3$, $\nu=\beta_3$.

%Let $\lambda=(\lambda_{1},\lambda_{2},\lambda_{3})\in \ZZ$ with $\lambda_{1}\ge\lambda_{2}\ge\lambda_{3}$ and $\lambda_{1}+\lambda_{2}+\lambda_{3}=0$.  % We abbreviate $\ell_{1}=\lambda_{1}-\lambda_{2}$ and $\ell_{2}=\lambda_{2}-\lambda_{3}$ and $V_{\lambda}({\bf 0})=V_{\lambda}(0,\ldots,0)$.
%Denote by $^{-}V_{\lambda}({\bf 0})$  the kernel of $(\sigma_1+id):V_{\lambda}({\bf 0})\to V_{\lambda}({\bf 0})$ (resp. of %$(\sigma_2+id):V_{\lambda}({\bf 0})\to V_{\lambda}({\bf 0})$), where we abbreviate   $V_{\lambda}({\bf 0})=V_{\lambda}(\ell,\ldots,\ell)$.  Similarly, %denote by $ {}^{+}V_{\lambda}({\bf 0}):V_{\lambda}({\bf 0})\to V_{\lambda}({\bf 0})$ (resp. by $V_{\lambda}({\bf 0}))^{+}:V_{\lambda}({\bf 0})\to %V_{\lambda}({\bf 0})$) the kernel of $\sigma_1-id$ (resp. of $\sigma_2-id$).

We say that a subset $\BB$ of  an $m$-dimensional vector space $V$
are {\it in general position} if any $m$-element subset of $\BB$ is
linearly independent. Clearly, if $\BB\subset V$ is generic then for
any linearly independent subsets   $S, S'\subset \BB$  we have (with
the convention that $\langle S\rangle$ is the linear span of $S$ in
$V$)
\begin{equation}
\label{eq:generic intersections}
\dim \left(\langle S\rangle\cap \langle S'\rangle\right)=\max(0,|S|+|S'|-m)  \ .
\end{equation}

\begin{theorem}
\label{th:intersections} For any $\lambda,\beta\in \ZZ^3$ we have:
$\BB^1_\lambda(\beta)\cap \BB^2_\lambda(\beta)=\emptyset$ and the
union $\BB^1_\lambda(\beta)\cup \BB^2_\lambda(\beta)$  is a generic
subset of $V_\lambda(\beta)$.

\end{theorem}

%We define for each $\beta=(\beta_{1},\beta_{2},\beta_{3})$ such that $V_{\lambda}(\beta)\ne 0$ the spaces
%$$ ^{1}V_{\lambda}(\beta)=\begin{cases}
 % {\rm Span}\{ u_{d-1}, u_{d-3},\ldots, u_{m } \} & \text{ if $d-m$ is even} \\
%{\rm Span}\{ u_{d-1}, u_{d-3},\ldots, u_{m-1 } \} & \text{ if $d-m$ is odd}
%\end{cases}$$
%$$  V_{\lambda}(\beta)^{1} =\begin{cases}
 % {\rm Span}\{ v_{d-1}, v_{d-3},\ldots, v_{m }  & \text{ if $d-m$ is even} \\
%{\rm Span}\{ v_{d-1}, v_{d-3},\ldots, v_{m-1 }  & \text{ if $d-m$ is odd}
%\end{cases}$$

\proof We will construct a third basis ${\mathcal B}_\lambda$ for
$V_\lambda(\beta)$ with certain convexity properties of the
coefficients of the transition matrix between the bases. First we
will address a general situation.

For each  Laurent polynomial $p\in \CC[q,q^{-1}]$ denote by  $\deg (p)$ the maximal degree of $p$.  We also denote $deg(0)=-\infty$.

\begin{theorem}
\label{th:generic bases in general} Let ${\mathcal
B}^i=\{b_1^i,\ldots,b_m^i\}$, $i=1,2$ be bases in an $m$-dimensional
vector space $V$ over $\CC(q)$ and ${\mathcal B}=\{b_1,\ldots,b_m\}$
be another basis in $V$ such that
\begin{equation}
\label{eq:generic bases in general}
b_\ell^1=\sum_{k=1}^\ell c_{k\ell}^1 b_k, ~b_\ell^2=\sum_{k=1}^\ell c_{k\ell}^2 b_{m+1-k}\ ,
\end{equation}
where all the coefficients $c_{k\ell}^i$ belong to $\CC[q,q^{-1}]\setminus \{0\}$ and satisfy:

\noindent $\bullet$ $\deg (c^i_{k+1,\ell})>\deg (c^i_{k,\ell})$ for all $1\le k<\ell\le m$.

\noindent $\bullet$ $\deg (c^i_{k,\ell})+\deg(c^i_{k+1,\ell+1})>\deg (c^i_{k+1,\ell})+\deg(c^i_{k,\ell+1})$ for all $1\le k<\ell< m$.

Then $\BB^1 \cap \BB^2 =\emptyset$ and  the union $\BB^1 \cup \BB^2 $  is a generic subset of $V$.
\end{theorem}

\proof First, we will prove a general combinatorial result.

Let $\lambda=(0\le\lambda_1\le\lambda_2\le\cdots\le\lambda_m\le
n)\in \ZZ^m$ be a reversed partition. Denote by $S^-=S^-_\lambda$ be
the set of all $(i,j)\in [1,m]\to [1,n]$ such that $j\le \lambda_i$
and let $S^+=S^+_\lambda$ be the complement of $S^-$, i.e.,
$S^+=[1,m]\to [1,n]\setminus S^-$. That is, $S^+$ is the set of all
$(i,j)\in [1,m]\to [1,n]$ such that $j\ge  \lambda_i+1$.

We say that a $m\times n$-matrix $A=(a_{ij})\in Mat_{m\times
n}(\ZZ)$ is {\it $\lambda$-convex} if the matrix coefficients
$a_{i,j}$  satisfy
\begin{equation}
\label{eq:lambda-convex}
a_{i,j}+a_{i',j'}>a_{i',j}+a_{i,j'}
\end{equation}
for any $i<i'$ in $[1,m]$, $j<j'$ in $[1,n]$ such that:  either
$\{i,i'\}\times \{j,j'\}\subset  S^-$, or $\{i\}\times
\{j,j'\}\subset S^-$, $\{i'\}\times \{j,j'\}\subset S^+$, or
$\{i,i'\}\times \{j,j'\}\subset  S^+$.

Clearly, the set of all $\lambda$-convex matrices is not empty and forms a semigroup (but not a monoid) under the matrix addition.

%and consider a natural linear order on the sets $[1,m]$,  $[1,n]$ and $[-n,-1]$.
%Denote
%$S_m^{n,n'}=\{(i,j)\in [1,m]\times([1,n]\sqcup [-n',-1]):either $i\le j$ or $i\ge .

%and equip it  with the product partial ordering.

For any map $\kappa:[1,m]\to [1,n]$ define the {\it $\kappa$-weight} of each $m\times n$-matrix $A$
 by
$$w_\kappa(A)=\sum_{i=1}^m a_{i,\kappa(i)} \ .$$

For each $\kappa:[1,m]\to [1,n]$ define two multiplicity functions $\mu^\varepsilon_\kappa: [1,n]\to \ZZ_{\ge 0}$, $\varepsilon\in \{-,+\}$ by
$$\mu^\varepsilon_\kappa(j)=|\{i:(i,j)\in S^\varepsilon, \kappa(i)=j\}| \ .$$

For any given functions $K^-,K^+:[1,n]\to \ZZ_{\ge 0}$ denote by
$\FF(K^-,K^+)$ the set of all maps $\kappa:[1,n]\to [0,m]$ such that
$\mu^\varepsilon_\kappa=K^\varepsilon$ for $\varepsilon\in \{-,+\}$.

% FALSE Clearly, $\FF(K^-,K^+)\ne \emptyset$ if and only if $\sum\limits_{i=1}^m K^-(i)+K^+(i)=m$ FALSE.

We say that a pair $(i,i')$ is an {\it inversion} of $\kappa:[1,m]\to [1,n]$ if $1\le i<i'\le m$,   $\kappa(i)>\kappa(i')$, and
\begin{equation}
\label{eq:non-preserved multiplicities}
\mu^\varepsilon_{\kappa\circ (i,i')} = \mu^\varepsilon_\kappa
\end{equation}
for all $\varepsilon\in \{-,+\}$, where $(i,i'):[1,m]\to [1,m]$ is the transposition interchanging only $i$ and $i'$.

Denote by $In(\kappa)\subset [1,m]\times [1,m]$ the set of all inversions of $\kappa$.
% and we write $\kappa\preceq \kappa'$ if $In(\kappa)\supseteq In(\kappa')$.

%Let    $K^{+} (\kappa), K^{-}(\kappa)\subset [1,n]\times \ZZ_{\ge 0}$  be the sets
%$K^{+}(\kappa)=\{ (i,j): j=|(\kappa^{-1}(i),i)\cap S^{+}|\}$ and   $K^{-}(\kappa)=\{ (i,j): %j=|(\kappa^{-1}(i),i)\cap S^{-}|\}$.

%Denote for  $K^{+},K^{-}\subset [1,n]$  by $\FF=\FF( K^{+},K^{-})$ the set of maps $\kappa:[1,m]\to [1,n]$ for %which $K^{+}(\kappa)=K^{+}$ and $K^{-}(\kappa)=K^{-}$.

\begin{proposition}
\label{pr:convex}
For any $K^-,K^+:[1,n]\to \ZZ_{\ge 0}$ such that $\FF(K^-,K^+)\ne \emptyset$ we have:

\noindent (a)   There exists a unique map $\kappa_{(K^-,K^+)}\in \FF( K^-,K^+)$ such that $In(\kappa_{(K^-,K^+)})=\emptyset$.

\noindent (b) For any $\kappa\in \FF(K^-,K^+)\setminus \{\kappa_{(K^-,K^+)}\}$  and any  $\lambda$-convex matrix $A$ one has:
\begin{equation}
\label{eq:kappa maximum}
w_{\kappa_{(K^-,K^+)}} (A)>w_\kappa(A) \ .
\end{equation}

\end{proposition}

%\proof
%Fix a $\lambda$-convex matrix $A$. Let   $K^-=\{ k_{1}^{-},\ldots k_\ell^{-}\}$ and
%$K^+=\{ k_{1}^{+},\ldots k_{\ell'}^{+}\}$, where $ k_{1}^{-}<\ldots <k_\ell^{-}$ and $k_{1}^{+}<\ldots< %k_{\ell'}^{+}$.

%\begin{proposition}
% \label{pr:convex}
% For any  functions $K^-,K^+:[1,n]\to [0,m]$  such that $\FF(K^-,K^+)\ne \emptyset$ there exists a unique map  %$\kappa=\kappa_{(K^-,K^+)}\in \FF( K^-,K^+)$ such that
%\begin{equation}
%\label{eq:kappa maximum}
%w_\kappa(A)=\max\limits_{\kappa'\in \FF( K^-,K^+)}  ~w_{\kappa'}(A)
%\end{equation}
% for any $\lambda$-convex $m\times n$-matrix $A$.
% \end{proposition}

\proof Prove (a).
We need the following result.

\begin{lemma}
\label{le: transpositions}
 Let $\kappa:[1,m]\to [1,n]$ and $i,i'\in [1,m]$. Then:

\noindent (a) $\mu^\varepsilon_{\kappa\circ (i,i')}=
\mu^\varepsilon_\kappa$ for all $\varepsilon\in \{-,+\}$ if and only
if: either $\{i,i'\}\times \{\kappa(i),\kappa(i')\}\subset  S^-$, or
$\{i\}\times \{\kappa(i),\kappa(i')\}\subset S^-$, $\{i'\}\times
\{\kappa(i),\kappa(i')\}\subset S^+$, or $\{i,i'\}\times
\{\kappa(i),\kappa(i')\}\subset  S^+$.

\noindent (b) If $(i,i')\in In(\kappa)$ then
\begin{equation}
\label{eq:inversion increases weight}
w_{\kappa\circ (i,i')}(A)>w_{\kappa}(A)
\end{equation}
for any $\lambda$-convex matrix $A$.
\end{lemma}

\proof Part (a) is obvious. Prove (b). Recall that $(i,i')\in
In(\kappa)$ if and only if $\mu^\varepsilon_{\kappa\circ (i,i')}=
\mu^\varepsilon_\kappa$ for all $\varepsilon\in \{-,+\}$,   $i<i'$,
$j<j'$, where $j=\kappa(i')$, $i'=\kappa(i)$. By the assumption,
$i<i'$, $j<j'$. Then
 $$w_{\kappa\circ (i,i')}(A)-w_{\kappa}(A)= a_{i',j'}+a_{i,j}-a_{i',j}-a_{i,j'}>0 $$
 by (a) and the $\lambda$-convexity of $A$.
 The lemma is proved.
\endproof

Now let us prove  the existence of $\kappa_0\in \FF(K^-,K^+)$ with no inversions.

\begin{lemma}
\label{le:extremal kappa0}
If $\kappa_0\in \FF(K^-,K^+)$ satisfies
\begin{equation}
\label{eq:extremal kappa0}
w_{\kappa_0}(A)\ge w_\kappa(A)
\end{equation}
for all $\kappa\in \FF(K^-,K^+)$ and  some $\lambda$-convex matrix
$A$, then $\kappa_0$ has no inversions. Moreover, if $\kappa\in
\FF(K^-,K^+)$ is such that $In(\kappa)\ne \emptyset$ then  then
there exists $\kappa_0\in \FF(K^-,K^+)$ with no inversions such that
$w_{\kappa_0}(A)> w_\kappa(A)$.
\end{lemma}

\proof Indeed, let $\kappa_0\in \FF(K^-,K^+)$ satisfy
(\ref{eq:extremal kappa0}). Then $\kappa_0$ has no inversions,
because if $(i,i')\in In(\kappa_0)$, then  Lemma \ref{le:
transpositions} would imply that $\kappa_0\circ (i,i')\in
\FF(K^-,K^+)$ and $w_{\kappa_0\circ (i,i')}(A)>w_{\kappa_0}(A)$,
which contradicts to (\ref{eq:extremal kappa0}). Finally, if
$\kappa\in \FF(K^-,K^+)$ and $In(\kappa)\ne \emptyset$, then Lemma
\ref{le: transpositions} implies that there exists $\kappa'$ that
$w_{\kappa'}(A)> w_\kappa(A)$ and then choose $\kappa_0\in
\FF(K^-,K^+)$ in such a way that $\kappa_0$ brings the maximum of
the weight $w_{\kappa'}(A)$ (such ''maximal'' $\kappa_0$ always
exists because the set  $\FF(K^-,K^+)$ is finite). And
(\ref{eq:extremal kappa0}) implies that $In(\kappa_0)=\emptyset$ and
$w_{\kappa_0}(A)> w_\kappa(A)$.
\endproof

Now we prove the uniqueness of such $\kappa_0$. We need the following result.

\begin{lemma}
\label{le:kappa m unique} Assume that $\kappa_0 \in\FF (K^-,K^+)$
has no inversions. Then $\kappa_0(m)$ is uniquely determined by
$(K^-,K^+)$, or, more precisely,
$$\kappa_0(m)=
\begin{cases} k^- & \text{if  $\exists$ $i_0\in [1,m-1]$: $\sum\limits_{i=\lambda_{i_0}+1}^m K^-(i)=m-i_0$}\\
k^+ & \text{otherwise}
\end{cases} ,
$$
where $k^\varepsilon=\max\{j:K^\varepsilon(j)\ne 0\}$ for $\varepsilon\in \{+,-\}$.
\end{lemma}

\proof Suppose $\kappa_0\in \FF(K^-,K^+)$ has no inversion and let
us fix any $i^+\in [1,m]$ such that $\kappa_0(i^+)=k^+$ and let us
fix any $i^-\in [1,m]$ such that $\kappa_0(i^-)=k^-$.

Suppose, by contradiction, that $\kappa_{0}(m)\notin \{k^+,k^{-}\}$,
i.e., $i^-<m$ and $i^+<m$. Then we construct an inversion of
$\kappa_0$ as follows.  If $(m,\kappa_0(m))\in S^+$ then, clearly,
$(i^+,m)$ is an inversion of $\kappa_0$ and if $(m,\kappa_0(m))\in
S^-$ then, clearly, $(i^-,m)$ is an inversion of $\kappa_0$. The
obtained contradiction proves that $\kappa_0(m)\in \{k^+,k^{-}\}$.

Denote also $i_{0,\varepsilon}=\max\{i: \kappa_0(i)=k^{\varepsilon}\}$.

Clearly,  $i_{0, \varepsilon}<i_{1,\varepsilon}$  for some
$\varepsilon\in \{+,-\}$, then Lemma \ref{le: transpositions}(a)
implies that the pair $(i_{0,\varepsilon},i_{1,\varepsilon})$ is an
inversion of $\kappa_0$.

Therefore $\kappa_0(m)$ equals either $k^+$ or $k^-$. Now assume
that $\kappa_0(m)=k^-$. Let $i_0<m$ be largest index such that
$\kappa(i_0)=k^+$. Since $\kappa_0$ has no inversions, we
immediately see that $(i,\kappa_0(i))\in S^-$ for each $i\in
[i_0,m]$. In its turn this implies that
$\sum\limits_{i>\lambda_{i_0}}  K^-(i)=m-i_0$.

In order to finish the proof we need to show the opposite
implication. Now assume that $\sum\limits_{i>\lambda_{i_0}}
K^-(i)=m-i_0$ for some $i_0<m$. This implies that
$(i,\kappa_0(i))\in S^-$ for each $i\in [i_0,m]$. In particular,
$(m,\kappa_0(m))\in S^-$, therefore, $\kappa_0(m)\ne k^+$ and, by
the above, one must have $\kappa_0(m)= k^-$. The lemma is proved.
\endproof

We finish the proof of the uniqueness of $\kappa_0$ by  induction in
$m$. According to Lemma \ref{le:kappa m unique}, if $\kappa_0\in
\FF(K^-,K^+)$ has no inversions then $\kappa_0(m)$ is uniquely
determined by $(K^-,K^+)$.  It is easy to see that the restriction
$\kappa_0|_{[1,m-1]}$ has no inversions and belongs to
$\FF(L^-,L^+)$, where
$$ L^\varepsilon(i)=\begin{cases}
         K^\varepsilon(i) & \text{ if $i\ne k^\varepsilon$  }\\
           K^\varepsilon(k^\varepsilon)-1 &  \text{ if $i= k^\varepsilon$, $\kappa_{0}(m)=k^\varepsilon$  }\\
               K^\varepsilon(k^\varepsilon) &  \text{ if $i= k^\varepsilon$, $\kappa_{0}(m)\ne k^\varepsilon$  }\\
        \end{cases} .
$$
Since $\kappa_0|_{[1,m-1]}$ has no inversions and  the pair
$(L^-,L^+)$ is uniquely determined by $(K^-,K^+)$, it follows that
$\kappa_0|_{[1,m-1]}$ is  (by the inductive hypothesis) uniquely
determined by $(K^-,K^+)$ in $\FF(L^-,L^+)$. This proves (a).

Part (b) follows from (\ref{eq:extremal kappa0}) for
$\kappa_0=\kappa_{(K^-,K^+)}$, which is unique according to the
already proved Part (a). Therefore, the inequality (\ref{eq:extremal
kappa0}) is strict and holds for all $\lambda$-convex matrices $A$.

This proves part (b). Therefore,  Proposition \ref{pr:convex} is proved.
\endproof

Now we are ready to finish the proof of Theorem \ref{th:generic
bases in general}. For simplicity, we identify $V$ with $\CC(q)^m$
by means of the basis $\BB$. Let $U$ be the $m\times 2m$ matrix
whose columns are $b^2_m,b^2_{m-1},\ldots,b^2_1;b_1^1,\ldots,b_m^1$.
For each embedding $\tau:[1,m]\hookrightarrow [1,2m]$ denote
$u_\tau=u_{1,\tau(1)}\cdots u_{m,\tau(m)}$; and for each
 subset $J\subset [1,2m]$ denote by $T_J$ of a of all embeddings
 $\tau:[1,m]\to [1,2m]$ such that $u_\tau\ne 0$. Since each minor is an alternating sum of $u_\tau$, $\tau\in T_J$,
 it suffices to prove  that  for each $J\subset [1,2m]$ there exists a unique $\tau_J$ such that $deg(u_{\tau_J})>deg(u_\tau)$ for any
 $\tau\in T_J\setminus \{\tau_J\}$.
Indeed, we will prove it using Proposition \ref{pr:convex}.

First, define a $m\times (m+1)$ integer matrix $A$ by:
$$a_{k,\ell}=\begin{cases}
\deg(u_{k,\ell}) & \text{ if $ k\ge \ell $ }\\
\deg(u_{k,\ell+m-1})    &\text{  if $k<\ell$}
\end{cases} =\begin{cases}
\deg(c_{m+1-k,m+1-\ell}^2) & \text{ if $ k\ge \ell $ }\\
\deg(c_{k,\ell-1}^1)    &\text{  if $k<\ell$}
\end{cases} $$

Clearly,  $a_{k,\ell}+a_{k+1,\ell+1}>a_{k+1,\ell}+a_{k,\ell+1}$ for
all $k,\ell\in [1,m-1]$ such that  $\ell-k\notin \{0,1\}$;
$a_{k,\ell}>a_{k,\ell+1}$ if $k\ge \ell$  and
$a_{k,\ell}<a_{k+1,\ell}$ if $\ell>k+1$.

Therefore, if we set $\lambda=(1,2,\ldots,n)$, then  $A$ is $\lambda$-convex in the notation (\ref{eq:lambda-convex}).

Furthermore, given an $m$-element subset $J\subset [1,2m]$, denote
$K^-=J\cap [1,m]$ and $ K^+=\{j-m:j\in J, j>m\}$. We will view these
subsets of $[1,m]$ as functions $K^-,K^+:[1,m]\to \{0,1\}$. For any
embedding $\tau:[1,m]\to [1,2m]$ such that $u_\tau\ne 0$ we define a
map $\varphi(\tau):[1,m]\to [1,m+1]$ by
$$\varphi(\tau)(k)=\begin{cases}
\tau(k) & \text{ if $\tau(k)\le k$}\\
\tau(k)+1-m &  \text{ if $\tau(k)\ge k+m$}
\end{cases} . $$
Clearly, $\varphi(\tau)$ is well-defined. It is easy to see that if
$\tau\in T_J$ then $\varphi(\tau)\in \FF(K^-,K^+)$ in the notation
of Proposition \ref{pr:convex} and, moreover, $\varphi:T_J\to
\FF(K^-,K^+)$ is a bijection. It is also easy to see that under this
bijection we have $$\deg(u_\tau)=w_{\varphi(\tau)}(A) \ .$$
Therefore, by  Proposition \ref{pr:convex}, there exists a unique
$\tau_J\in T_J$ such that $\varphi(\tau_J)=\kappa_{(K^-,K^+)}$ and
$deg(u_{\tau_J})>deg(u_\tau)$ for any $\tau\in T_J\setminus
\{\tau_J\}$.  Theorem \ref{th:generic bases in general} is proved.
\endproof

Now we will finish the proof of Theorem \ref{th:intersections}.
In what follows we fix $\lambda, \beta\in \ZZ^3$ such that $V_\lambda(\beta)\ne 0$.

Denote by $\BB_\lambda$ the {\it dual canonical basis} for
$V_\lambda$ (see e.g.,  \cite{bz1,bz2}). Recall from \cite{bz1,bz2}
that each element $b\in \BB_\lambda$ is naturally labeled by  a
quadruple ${\bf m}=(m_{1},m_{2},m_{12},m_{21})$ of non-negative
integers such that $ m_{1} m_{2}=0$, such that  $m_i+m_{ij}\le
\ell_i$ for $\{i,j\}=\{1,2\}$ (here we abbreviated
$\ell_1=\lambda_1-\lambda_2$, $\ell_2=\lambda_2-\lambda_3$). For any
$\beta\in \ZZ^3$ the intersection $\BB_{\lambda}
(\beta)=\BB_{\lambda} \cap V_{\lambda}(\beta)$ is a basis in
$V_{\lambda}(\beta)$ and $b_{\bf m}\in \BB_{\lambda} (\beta)$ if and
only if
$m_{1}\alpha_1+m_2\alpha_2+(m_{12}+m_{21})(\alpha_1+\alpha_2)=\lambda-\beta$.

The basis $\BB_\lambda$ determined (up to a scalar multiple) by the
property that for each $d\ge 0$, $\beta\in \ZZ^3$, and $i\in
\{1,2\}$ a part of $\BB_\lambda$ spans the subspace $\ker
F_i^{d+1}\cap V_\lambda(\beta)$  and (another) part of
$\BB_\lambda$ spans $\ker E_i^{d+1}\cap V_\lambda(\beta)$.

\begin{lemma}[\cite{bz1,bz2}]
\label{le:canonical basis}
 The action of  $E_i, F_i \in U_{q}(gl_{3})$, $i\in \{1,2\}$ on  $\BB_\lambda $ is given by:
\begin{equation}
\label{eq: Ei-action}
E_i(b_{\bf m})= ( m_i+m_{ji})_{q}b_{{\bf m}-e_i}+ (m_{ji})_{q} b_{{\bf m}-e_i'} +(m_{ji})_{q} b_{{\bf m}-e_i''}
\end{equation}
for $\{i,j\}=\{1,2\}$, where $e_1=(1,0,0,0)$, $e_1'=(1,0,-1,1)$,$e_1''= (0,-1,0,1)$, and $e_2=(0,1,0,0)$, $e_2'=(0,1,1,-1)$, $e_2''= (-1,0,1,0)$.
\begin{equation}
\label{eq: Fj-action}
F_i(b_{\bf m})= ( m_j+m'_{ij})_{q}b_{{\bf m}+e_i''}+ (m'_{ij})_{q} b_{{\bf m}-f_i'} +(m'_{ij})_{q} b_{{\bf m}+e_i}
\end{equation}
for $\{i,j\}=\{1,2\}$, where $f'_1=(-1,0,0,1)$, $f_2'=(0,-1,1,0)$, $m'_{ij}=\ell_i-m_i-m_{ij}$.

\end{lemma}

For $v\in V_\lambda\setminus \{0\}$ and $i\in \{1,2\}$ denote
$$\ell_i^+(v):=\max(k:E_i^k(v)\ne 0\} ,   \ell_i^-(v):=\max(k:F_i^k(v)\ne 0\} \ .$$
It follows from Lemma \ref{eq: Ei-action} that $\ell_i^+(b_{\bf m})=m_i+m_{ji}$ and $\ell_i^-(b_{\bf m})=m_j+m'_{ij}$.
%By definition (\ref{eq:Gelfand-Tsetlin basis}),
%\begin{equation}
%\label{eq:Gelfand-Tsetlin exponents}
%\ell_i^+(b^i_{(\mu,\nu)})=\begin{cases}
%\mu_1-\nu & \text{if $i=1$}\\
%\nu-\mu_2 & \text{if $i=2$}\\
%\end{cases}, \ell_i^-(b^i_{(\mu,\nu)})=\begin{cases}
%\nu-\mu_2  & \text{if $i=1$}\\
%\mu_1-\nu & \text{if $i=2$}\\
%\end{cases}.
%\end{equation}
Denote also
$$d_i^-=d_i^-(\lambda,\beta)=\min_{v\in V_\lambda(\beta)\setminus \{0\}} \ell_i^-(v),
d_i^+=d_i^+(\lambda,\beta)=\min_{v\in V_\lambda(\beta)\setminus \{0\}} \ell_i^+(v) \ .$$
%$$M_i^-=d_i(\lambda,\beta)^-=\max_{v\in V_\lambda(\beta)\setminus \{0\}} \ell_i^+(v), M_i^+=d_i(\lambda,\beta)^+=\max_{v\in V_\lambda(\beta)\setminus %\{0\}} \ell_i^+(v)$$

The properties of $\BB_\lambda$ established above imply the following simple result.

\begin{lemma}
\label{le:kernel of E,F} We have $d_1^+=d_1^-+\beta_2-\beta_1$, $d_2^+=d_2^-+\beta_3-\beta_2$, and
$$d_1^-=(\beta_1-\lambda_2)_++(\lambda_2-\beta_2)_+, \, d_2^-=(\lambda_2-\beta_3)_++(\beta_2-\lambda_2)_+$$
where $(a)_+:=\max(a,0)$ for $a\in \ZZ$.
\end{lemma}

%Let $d_iEach of the spaces $V_\lambda(\beta)\cap\ker E_i$, $V_\lambda(\beta)\cap\ker F_i$, $i\in \{1,2\}$ is at most $1$-dimensional. Moreover,

%$\bullet$ $V_\lambda(\beta)\cap \ker E_1\ne 0$ if and only if $\beta_1\ge \lambda_2\ge \beta_2$,
%
%$\bullet$ $V_\lambda(\beta) \cap \ker E_2\ne 0$ if and only if $\beta_3\le \lambda_2\le \beta_2$,
%
%$\bullet$ $V_\lambda(\beta) \cap \ker F_1\ne 0$ if and only if $\beta_1\le \lambda_2\le \beta_2$,
%
%$\bullet$  $V_\lambda(\beta)\cap \ker F_2\ne 0$ if and only if $\beta_3\ge \lambda_2\ge \beta_2$.
%\end{lemma}

Furthermore, we will label the  Gelfand-Tsetlin bases
$\BB^i_\lambda(\beta)=\BB_\lambda^i\cap V_\lambda(\beta)$, $i\in
\{1,2\}$ and the dual canonical basis
$\BB_\lambda(\beta):=\BB_\lambda\cap V_\lambda(\beta)$ as follows
(with the convention that $m=m_{\lambda,\beta}=\dim
V_\lambda(\beta)$).
\begin{equation}
\label{eq:labeled bases}
\BB^i_\lambda(\beta)=\{b_1^i,b_2^i,\ldots,b_m^i\},~ \BB_\lambda(\beta)=\{b_1,b_2,\ldots,b_m\} \ ,
\end{equation}
where $b_k^i$ (for $k=1,2,\ldots,m$, $i\in \{1,2\})$ is the only
element of  $\BB^i_\lambda(\beta)$  such that
$\ell_i^-(b_k^1)=d_i^-+k-1$; and  $b_k$ is the only element of
$\BB_\lambda(\beta)$  such that $\ell_1^-(b_k)=d_1^-+k-1$.

\begin{lemma}
\label{le:triangularity}
The above labeling satisfies for all $\ell\in [1,m]$:
$$b^1_\ell\in Span\{b_1,\ldots,b_\ell\},~b^2_\ell\in Span\{b_{m-\ell+1},\ldots,b_m\} \ .$$

\end{lemma}

\begin{proof} Indeed, by definition, $b_\ell^1\in Ker F_1^{d_1^-+\ell}\cap V_\lambda(\beta)$.
But according to the basic property of $\BB_\lambda$, the vectors
$b_1,b_2,\ldots,b_\ell$ form a basis in $Ker F_1^{d_1^-+\ell}\cap
V_\lambda(\beta)$. This proves the first inclusion. To prove the
second, it suffices to show that for $d=\ell_2^-(b_\ell^2)$ the
intersection $Ker F_1^{d+1}\cap V_\lambda(\beta)$ is spanned by
$\{b_{m-\ell+1},\ldots,b_m\}$. That is, we have to show that
$\ell_2^-(b_\ell^1)\ge \ell_2^-(b_k)$ if and only if $k+\ell\ge
m+1$.

It is easy to see that for $b_k\in V_\lambda(\beta)$, we have
\begin{equation}
\label{eq:explicit relabeling}
b_k=b_{(m_1,m_2,m_{12},m_{21})}\ ,
\end{equation}
where $m_1=(\beta_2-\lambda_2)_+$, $m_2=(\lambda_2-\beta_2)_+$,
$m_{12}=\lambda_1-\beta_2-k-d_1^-(\lambda,\beta)+1$,
$m_{21}=\lambda_1-\beta_1-m_1-m_{12}$. Therefore,
$$\ell_2^+(b_k)=m_2+m_{12}=\lambda_1-\beta_2-k-d_1^-(\lambda,\beta)+1+(\lambda_2-\beta_2)_+=d_2^+(\lambda,\beta)+m-k$$
because
$m=1+\lambda_1-(\beta_1-\lambda_2)_+-(\beta_3-\lambda_2)_+-\max(\lambda_2,\beta_2)$.
Therefore, $\ell_2^-(b_k)=d_2^-(\lambda,\beta)+m-k$ and the
inequality $\ell_2^-(b_\ell^2)\ge \ell_2^-(b_k)$ is equivalent to
$d_2^-(\lambda,\beta)+\ell-1\ge d_2^-(\lambda,\beta)+m-k$, i.e.,
$k+\ell\ge m+1$. This proves the second inclusion. The lemma is
proved.
\end{proof}

The above result implies that $b_1^1$ is proportional to $b_1$ and
$b_1^2$ is proportional to $b_m$. In what follows we set
$b_1^1=b_1$ and $b_1^2=b_m$ for each $\beta$ such that
$d_1^-(\lambda,\beta)=0$. It is easy to see that this setting agrees
with the definition (\ref{eq:Gelfand-Tsetlin basis}) of both
$\BB^1_\lambda$ and $\BB^2_\lambda$, and, moreover, determines them
uniquely.

The following result which completely describes the expansion of
Gelfand-Tsetlin bases along the dual canonical basis of
$V_\lambda(\beta)$ will  finish the proof of Theorem
\ref{th:intersections}.

\begin{proposition}
\label{pr:generic bases any beta} The bases
$\BB^1:=\BB^1_\lambda(\beta)$, $\BB^2:=\BB^2_\lambda(\beta)$, and
$\BB:=\BB_\lambda(\beta)$ labeled by (\ref{eq:labeled bases})
satisfy the hypotheses of Theorem \ref{th:generic bases in general}.
More precisely,
\begin{equation}
\label{eq:generic bases in particular}
 b_\ell^1=\sum_{k=1}^\ell c_{k\ell}^1b_k\ , b_\ell^2=\sum_{k=1}^\ell c_{k\ell}^2 b_{m+1-k}
\end{equation}
where all the coefficients $c_{k\ell}^1,c_{k\ell}^2$ belong to $\CC[q,q^{-1}]\setminus \{0\}$ and satisfy the recursion
\begin{equation}
\label{eq:c-differences}
\deg (c_{k,\ell}^i)-\deg (c_{k-1,\ell}^i)= |\beta_2-\lambda_2|+2(\ell-k)+1
 \end{equation}
for all $1<  k\le\ell \le m$, $i\in \{1,2\}$.

\end{proposition}

\begin{proof} First, the triangularity (\ref{eq:generic bases in particular}) follows from Lemma \ref{le:triangularity}.
Furthermore, we will establish some symmetries which take Gelfand-Tsetlin bases into each other and preserve the dual canonical basis.
The following well-known fact is  a particular case of a more general result (see e.g., \cite{bz2}).

\begin{lemma}

\label{le:symmetries}
For each $\lambda=( \lambda_1\ge\lambda_2\ge\lambda_3)$ we have:

\noindent (a)  There is a unique (up to a scalar multiple) linear isomorphism $\psi_\lambda:V_\lambda\to V_{\lambda^*}$, where
$\lambda^*=(-\lambda_3,-\lambda_2,-\lambda_1)$  such that
$$\psi_\lambda(E_i(v))=F_i(\psi_\lambda(v)),\psi_\lambda(F_i(v))=E_i(\psi_\lambda(v)), \psi_\lambda(K_\beta(v))=K_{-\beta} (\psi_\lambda(v))$$
for any $i\in \{1,2\}$, $v\in V_\lambda$, $\beta\in \ZZ^3$.

\noindent (b) there is a unique (up to a scalar multiple) isomorphism $\eta_\lambda:V_\lambda\to V_{\lambda}$ such that
$$\eta_\lambda(E_i(v))=F_{3-i}(\eta_\lambda(v)),\eta_\lambda(F_i(v))=E_{3-i}(\eta_\lambda(v)), \eta_\lambda(K_\beta(v))=K_{\beta^{op}}(\eta_\lambda(v))$$
for any $i\in \{1,2\}$, $v\in V_\lambda$, $\beta\in \ZZ^3$, where $(\beta_1,\beta_2,\beta_3)^{op}=(\beta_3,\beta_2,\beta_1)$.

\noindent (c) $\psi_\lambda\circ \eta_\lambda= \eta_{\lambda^*}\circ \psi_{\lambda^*}$.
\end{lemma}

Clearly,  in terms of the labeling (\ref{eq:labeled bases}), $ \psi_\lambda(b_k^i)=b_k^i$ for $i\in \{1,2\}$  for all $k\in [1,m]$ and
$\eta_\lambda(b_k^1)=b_k^2, ~\eta_\lambda(b_k^2)=b_k^1$ for all $k\in [1,m]$.

The results of \cite{bz2} imply that the symmetries $\psi_\lambda$ and $\eta_\lambda$ preserve the dual canonical basis:
$\psi_\lambda(\BB_\lambda(\beta))=\BB_{\lambda^*}(-\beta)$ and $\eta_\lambda(\BB_\lambda(\beta))=\BB_\lambda(\beta^{op})$.
We can refine this fact as follows.

\begin{lemma} We have for $k=1,2,\ldots,m$:

\noindent (a) $\psi (b_k)=b_k$ for $k=1,2,\ldots,m$.
%where $m=\dim V_\lambda(\beta)=\dim V_{\lambda^*}(\beta^*)=\dim V_\lambda(\beta^{op})$.

\noindent (b) $\eta_\lambda(b_k)=b_{m+1-k}$.
\end{lemma}

\begin{proof}
%Since $\varphi_\lambda$ is an involution which takes $\beta_2-\lambda_2$ to $\lambda_2-\beta_2$, it suffices to consider the case when $\beta_2\ge %\lambda_2$. Then, u
Prove (a).  Using the fact that
$\ell_1^-(\psi_\lambda(v))=\ell_1^+(v)$ and
$d_1^+(\lambda,\beta)=d_1^-(\lambda^*,-\beta)$, we see that
$\ell_1^-(\psi_\lambda(b_k))=d_1^-(\lambda^*,-\beta)+k-1=\ell_1^-(b_k)$.
Therefore, $\psi_\lambda(b_k)=b_k$.

Part (b) follows follows directly from Lemma \ref{le:triangularity}.
\end{proof}

Using the symmetry $\psi_\lambda$ we see that it suffices to  prove
Proposition  \ref{pr:generic bases any beta} only in the case when
$\lambda_2\le \beta_2$. On the other hand,  using the symmetry
$\eta_\lambda$ and the property (\ref{eq:c-differences relabeled})
for $i=1$ we see that we see  that the  coefficients $c_{k\ell}^2$
also satisfy (\ref{eq:c-differences relabeled}).

Therefore, all we have to prove is that in the case when $\beta_2\ge
\lambda_2$ all the coefficients of the first expansion
(\ref{eq:generic bases in particular}) satisfy the property
(\ref{eq:c-differences relabeled}) with $i=1$. In order to achieve
this goal, it is convenient to re-label the bases
$\BB^1_\lambda(\beta)$ and  $\BB_\lambda(\beta)$  as follows.
\begin{equation}
\label{eq:relabeling}
\tilde b_\ell^1:=b_{\ell-d_1^-+1}^1,~\tilde b_\ell:=b_{\ell-d_1^-+1}
\end{equation}
for $\ell=d_1^-,d_1^-+1,\ldots, D_1^-$, where
$d_1^-=d_1^-(\lambda,\beta)=(\beta_1-\lambda_2)_+$ and
$D_1^-=D_1^-(\lambda,\beta)=\lambda_1-\beta_2-(\beta_3-\lambda_2)_+$;
denote also  $\tilde c_{k,\ell}:=c_{k-d_1^-+1,\ell-d_1^-+1}^1$.

Then in terms of the this re-labeling  it suffices to prove that if
$\beta_2\ge \lambda_2$, then:
\begin{equation}
\label{eq:generic bases in particular relabeled}
\tilde b_\ell^1=\sum_{k=d_1^-}^\ell \tilde c_{k\ell} \tilde b_k\ ,
\end{equation}
for all $\ell\in [d^-_1,D^-_1]$, where each $\tilde c_{k\ell}=\tilde c_{k\ell}^\beta\in \CC[q,q^{-1}]\setminus \{0\}$ and
\begin{equation}
\label{eq:c-differences relabeled}
\deg (\tilde c_{k,\ell}^{\beta})-\deg (\tilde c_{k-1,\ell}^{\beta})= \beta_2-\lambda_2+2(\ell-k)+1 \ .
 \end{equation}
for all $d_1^-(\lambda,\beta)< k\le\ell \le D_1^-(\lambda,\beta)$.

Furthermore, we need the following  reformulation of Lemma \ref{le:canonical basis}.

\begin{lemma} Let $\beta'\in \ZZ^3$ be such that  $\lambda_2\le \beta_2'-1$. Then the action
of the operator $E_1:V_\lambda(\beta')\to
V_\lambda(\beta'+\alpha_1)$ relative to the  bases
$\BB_\lambda(\beta')$ and $\BB_\lambda(\beta'+\alpha_1)$  is given
by
\begin{equation}
\label{eq: E1-action}
E_1(\tilde b_k)= ( \beta'_2-\beta'_1+k)_q \tilde  b_{k+1}+ (\lambda_2-\beta'_1+k)_q \tilde  b_k
 \end{equation}
for $k\in [d_1^-(\lambda,\beta'),D_1^-(\lambda,\beta')]$ (where $\tilde b_k=0$ in $V_\lambda(\beta'+\alpha_1)$ if $k=\beta_1'-\lambda_2)$.
\end{lemma}

\begin{proof} First of all (in the notation of Lemma \ref{le:canonical basis})
we have $m_2=0$ and $m_1=\beta'_2-\lambda_2$ because $m_1-m_2=\beta'_2-\lambda_2> 0$ and $m_1m_2=0$.

Therefore, we can rewrite (\ref{eq: Ei-action}) as follows:
$$E_1(b_{(m_1,0,m_{12},m_{21}})= ( m_1+m_{21})_q b_{(m_1-1,0,m_{12},m_{21})}+ (m_{21})_q b_{(m_1-1,0,m_{12}+1,m_{21}-1)}$$

Then, after ''translating''  the equation $b_{\bf m}=\tilde b_k$ in  $V_\lambda(\beta')$ (based on (\ref{eq:explicit relabeling})) into
$$b_{(\beta_2'-\lambda_2,0,\lambda_1-\beta_2'-k,\lambda_2-\beta_1'+k)}=\tilde b_k\ ,$$
we immediately obtain (\ref{eq: E1-action}). The lemma is proved.
\end{proof}

We proceed by induction in $\beta$ with respect to the following
partial order on $\ZZ^3$: $\beta\ge \beta'$ if
$\beta=\beta'+n\alpha_1$ for some $n\ge 0$. The base of induction
consists of all those $\beta$ which satisfy: $V_\lambda(\beta)\ne 0$
but $V_\lambda(\beta-\alpha_1)=0$. It is easy to see that for such
$\beta$ one has $\dim V_\lambda(\beta)=1$ and $\beta_2\ge \lambda_2$
and, therefore,  we have nothing to prove. Now we assume that both
weight spaces  $V_\lambda(\beta)$ and $V_\lambda(\beta-\alpha_1)$
are non-zero and $\beta_2\ge \lambda_2$.

Note that if $b_{\ell-1}^1\in V_\lambda(\beta-\alpha_1)$ and
$\ell>(\beta_1-\lambda_2)_+$, then $b_\ell^1=E_1(b_{\ell-1}^1)$.
Therefore, applying (\ref{eq: E1-action}) to  (\ref{eq:generic bases
in particular relabeled}) with  $\beta-\alpha_1$ and $\ell-1$, we
obtain an expansion in $V_\lambda(\beta)$:
$$ \tilde b_\ell^1=\sum_{k< \ell} \tilde c_{k,\ell-1}^{\beta'} E_1(\tilde b_k)=\sum_{k<\ell}
\tilde c_{k,\ell-1}^{\beta'} \left(( \beta'_2-\beta'_1+k)_q \tilde  b_{k+1}+ (\lambda_2-\beta'_1+k)_q \tilde  b_k\right) ,$$
 where $\beta'=(\beta_1-1,\beta_2+1,\beta_3)=\beta-\alpha_1$. This implies that
$$ \tilde c_{k,\ell}^{\beta}= \tilde c_{k-1,\ell-1}^{\beta'}(\beta_2-\beta_1+k+1)_q+\tilde c_{k,\ell-1}^{\beta'}(\lambda_2-\beta_1+k+1)_q $$
with the convention that $c_{k,\ell}^{\beta'}=0$ for $k=\ell+1$. In
particular, this proves that each $c_{k,\ell}^{\beta}$ is a Laurent
polynomial in $q$.  Then the inductive hypothesis
(\ref{eq:c-differences relabeled}) taken for $\beta'$ and $\ell-1$
implies that for $k< \ell$ the  the degree of the second summand
minus the degree of the first summand equals $2(\ell-k)-1>0$.
Therefore,
$$ \deg( \tilde c_{k,\ell}^{\beta})=
\begin{cases}
\deg(\tilde c_{k,\ell-1}^{\beta'})+ \lambda_2-\beta_1+k +1  & \text{if $k< \ell$}\\
\deg(\tilde c_{\ell-1,\ell-1}^{\beta'}) +\beta_2-\beta_1+\ell+1 &\text{if $k= \ell$}\\
\end{cases}\ .$$
This implies that for $k<\ell$ we have
$$   \deg(\tilde c_{k,\ell}^{\beta})-\deg(\tilde c_{k-1,\ell}^{\beta} )=\deg(\tilde c_{k,\ell-1}^{\beta'})+
\lambda_2-\beta_1+k +1-(\deg(\tilde c_{k-1,\ell-1}^{\beta'})+ \lambda_2-\beta_1+k ) $$
$$=1+\deg(\tilde c_{k,\ell-1}^{\beta'})-\deg(\tilde c_{k-1,\ell-1}^{\beta'})= \beta_2-\lambda_2+2(\ell-k)+1  $$
by the same inductive hypothesis.
Similarly, for  $k= \ell$ we have
$$ \deg(\tilde c_{\ell,\ell}^{\beta})-\deg(\tilde c_{\ell-1,\ell}^{\beta} )=\deg(\tilde c_{\ell-1,\ell-1}^{\beta'}) +\beta_2-\beta_1+\ell+1-
(\deg(\tilde c_{\ell-1,\ell-1}^{\beta'})+ \lambda_2-\beta_1+\ell )
 $$
 $$=\lambda_2-\beta_2+1 \ .$$

This finishes the induction.  The recursion  (\ref{eq:c-differences
relabeled})  is proved. Therefore, the recursion
(\ref{eq:c-differences}) with $i=1$ and $\beta_2\ge \lambda_2$ is
proved. Proposition \ref {pr:generic bases any beta} is proved.
\end{proof}

Therefore,  Theorem \ref{th:intersections} is proved. \endproof

Now we will finish the proof of Theorem \ref{th:exterior and symmetric  products}.

\begin{lemma}
\label{le:double intersections} Under the identification
$(V_{\beta_1}\otimes  V_{\beta_2} \otimes
V_{\beta_3})^\lambda=V_\lambda(\beta)$  (in the notation of
(\ref{eq:notation highest vectors})) we have for  $\lambda, \beta\in
\ZZ_{\ge 0}^3$, and $\varepsilon\in \{-,+\}$:
$$((V_{\beta_1}\bullet_\varepsilon  V_{\beta_2}) \otimes  V_{\beta_3})^\lambda=Span\{\BB_\lambda^1(\beta)^\varepsilon\},
~(V_{\beta_1} \otimes (V_{\beta_2} \bullet_\varepsilon    V_{\beta_3}))^\lambda=Span\{\BB_\lambda^2(\beta)^\varepsilon\} \ .$$

\end{lemma}

\proof In the notation (\ref{eq:isotypic decomposition multiplicity
free}) each subspace $(J_{\beta_1,\beta_2}^m\otimes V_{\beta_3
})^\lambda$ (resp. $( V_{\beta_1 }\otimes
J_{\beta_2,\beta_3}^m)^\lambda$)  is spanned by all those
$b^1_{(\mu,\nu)}\in V_\lambda(\beta)$ (resp. $b^2_{(\mu,\nu)}\in
V_\lambda(\beta)$) which satisfy $\mu_2=m$.

Using the definition \eqref{eq:braided product}:
$$V_{\beta_1}\bullet_+  V_{\beta_2}  =\bigoplus_{m\in 2\ZZ} J_{\beta_1,\beta_2}^m, V_{\beta_1}\bullet_-
V_{\beta_2}=\bigoplus_{m\in 2\ZZ+1} J_{\beta_1,\beta_2}^m \ $$
and Lemma \ref{le:isotypic gl2} from Appendix, we  finish the proof of Lemma \ref{le:double intersections}.   \endproof

Furthermore, using Theorem \ref{th:intersections}, Lemma
\ref{le:double intersections}, and (\ref{eq:generic intersections})
with $S=\BB^1_\lambda(\beta)^\varepsilon$,
$S'=\BB^2_\lambda(\beta)^\varepsilon$, we obtain
$$\dim (V_{\beta_1}\bullet_\varepsilon  V_{\beta_2} \bullet_\varepsilon  V_{\beta_3})^\lambda
=\delta_{\lambda_3,0}\cdot \max(0,|\BB^1_\lambda(\beta)^\varepsilon|+|\BB^2_\lambda(\beta)^{\varepsilon}|-\dim V_\lambda(\beta)) $$
$$=\delta_{\lambda_3,0}\cdot \max(0,|\BB^1_\lambda(\beta)^\varepsilon|-|\BB^2_\lambda(\beta)^{-\varepsilon}|)$$

Now define a bijection $S_\lambda:\BB^1_\lambda\to \BB^2_\lambda$ by
$b^1_{(\mu,\nu)}\mapsto b^1_{(\mu',\nu')}$, where
$$\mu'_1=\lambda_1-\min(\mu_1,\nu-\mu_2),
\mu'_2=\min(\lambda_2,\mu_1+\mu_2-\nu)-\mu_2,
\nu'=\lambda_1+\lambda_2+\lambda_3-\mu_1-\mu_2-\mu_3$$ Clearly,  the
inverse of $S_\lambda$ is given by the same formula and
$S_\lambda(\BB^1_\lambda(\beta)^\varepsilon)=\BB^2_\lambda(\beta)^\delta$
, where $\delta=\varepsilon\cdot (-1)^{\min(\lambda_2,\beta_2)}$.
Therefore, if $\min(\lambda_2,\beta_2)$ is odd then
$|\BB^1_\lambda(\beta)^\varepsilon|=|\BB^2_\lambda(\beta)^{-\varepsilon}|$
and $(V_{\beta_1}\bullet_\varepsilon  V_{\beta_2}
\bullet_\varepsilon  V_{\beta_3})^\lambda=0$ for $\varepsilon \in
\{-,+\}$. Now assume that $\min(\lambda_2,\beta_2)$ is even, in
particular,
$|\BB^1_\lambda(\beta)^\varepsilon|=|\BB^2_\lambda(\beta)^{\varepsilon}|$
for $\varepsilon \in \{-,+\}$, and $\lambda_3=0$. Therefore,
$$\dim (V_{\beta_1}\bullet_\varepsilon  V_{\beta_2} \bullet_\varepsilon  V_{\beta_3})^\lambda=\max(0,|\BB^1_\lambda(\beta)^\varepsilon|-
|\BB^1_\lambda(\beta)^{-\varepsilon}|) \ .$$

It is easy to see that $|\BB^1_\lambda(\beta)^\varepsilon|=|\{m\in \ZZ:c(\lambda,\beta)\le m\le d(\lambda,\beta), (-1)^m=\varepsilon\}|$,
where
$$c(\lambda,\beta)=\max(0,\lambda_2-\beta_3),~d(\lambda,\beta)=\min(\lambda_2,\beta_2)-\max(0,\lambda_2-\beta_1) \ .$$
Therefore, $\dim (V_{\beta_1}\bullet_\varepsilon  V_{\beta_2}
\bullet_\varepsilon  V_{\beta_3})^\lambda=1$ if and only if
$(-1)^{\max(0,\lambda_2-\beta_3)}=(-1)^{\max(0,\lambda_2-\beta_3)}=\varepsilon$,
and $\dim (V_{\beta_1}\bullet_\varepsilon  V_{\beta_2}
\bullet_\varepsilon  V_{\beta_3})^\lambda=0$ otherwise.  Theorem
\ref{th:exterior and symmetric  products} is proved.   \endproof

To finish the proof of Theorem \ref{th:exterior powers and symmetric
cube}  note that in the notation of Theorem \ref{th:exterior and
symmetric  products}:
$$V_\ell\bullet_\varepsilon  V_\ell \bullet_\varepsilon  V_\ell \cong \bigoplus_\lambda V_\lambda \ ,$$
where the summation is over all $\lambda=(\lambda_1\ge \lambda_2\ge
0)$ such that: $\lambda_1+\lambda_2= 3\ell$ and
$(-1)^{\lambda_2}=(-1)^{(\lambda_2-\ell)_+}=\varepsilon$. Therefore,
Theorem \ref{th:exterior powers and symmetric cube} is proved.
\endproof

\section{Appendix: Quantum matrices}
\label{sect:Appendix}

%%%%%%%%%%%%%%%%%

The algebra $\CC_q[M_{2\times 2}]$ of {\it quantum} $2\times 2$-matrices
is the $Q(q)$-algebra generated by $a,b,c,d$  subject to the relations
$$ca=qac,~ ba=qab,~ dc= qcd, ~db=qbd,~cb=bc,~da-ad=(q-q^{-1})bc \ . $$

Let $\CC_q[M_{d\times k}]$ is the algebra of {\it quantum}
$m\times n$-matrices,
i.e., the $Q(q)$-algebra generated by $x_{ij}$, $1\le i\le d$, $1\le j\le k$ with the
following  defining relations: for each $1\le i<i'\le d,1\le j<j'\le k$ the subalgebra of
$\CC_q[M_{m\times n}]$ generated by the four elements $a=x_{i,j}, b=x_{i,j'}, c=x_{i',j}, d=x_{i',j'}$
is isomorphic to $\CC_q[M_{2\times 2}]$.

The algebra $\CC_q[M_{d\times k}]$ is graded by $\ZZ_{\ge 0}^d\times \ZZ_{\ge 0}^k$
via ${\rm deg}(x_{ij})=e'_i+e_j$,  where $e'_1,\ldots,e'_d$ (resp. $e_1,\ldots,e_k$) is the standard basis in
$\ZZ^d$ (resp. in $\ZZ^k$). That is,
$$\CC_q[M_{d\times k}]=\bigoplus\limits_{(\gamma;\delta)\in \ZZ_{\ge 0}^d\times \ZZ_{\ge 0}^k} \CC_q[M_{d\times k}]_{\gamma;\delta} \ ,$$
where
$\CC_q[M_{d\times k}]_{\gamma;\delta}$ the set of all $x\in \CC_q[M_{d\times k}]$ such that
$\deg(x)=(\gamma;\delta)$.

There is a natural action of $U_{q }(gl_d(\CC)\times gl_k(\CC))$ on $\CC_q[M_{d\times k}]$ via
$$E'_{i'}(x_{i,j})=\delta_{i,i'}x_{i+1,j},~F'_{i'}(x_{i,j})=\delta_{i-1,i'}x_{i-1,j}, K'_\lambda(x_{i,j})=q^{-\lambda_i}\cdot x_{ij}$$
for $i'=1,2,\ldots,d-1$, $\lambda\in \ZZ_{\ge 0}^d$;
$$E_{j'}(x_{ij})=\delta_{j',j-1}x_{i,j-1},~F_{j'}(x_{ij})=\delta_{j,j'}x_{i,j+1},~ K_{\mu}(x_{i,j})=q^{\mu_j}x_{i,j} $$
for $i,i'=1,2,\ldots,d-1$, $j,j'=1,2,\ldots,k-1$, $\lambda\in
\ZZ_{\ge 0}^d$; and $u(xy)={\bf m}(\Delta(u)(x\otimes y))$ for any
$u\in U_{q }(gl_d(\CC)\times gl_k(\CC))$ and $x,y\in
\CC_q[M_{d\times k}]$, where ${\bf m}$ is the multiplication
$\CC_q[M_{d\times k}]\otimes \CC_q[M_{d\times k}]\to
\CC_q[M_{d\times k}]$.

This definition implies that $K'_\lambda K_\mu(x)=q^{(\mu,\delta)-(\lambda,\gamma)}\cdot x$
for $x\in \CC_q[M_{d\times k}]_{\gamma;\delta}$ and
$$E_j(xy)=E_j(x) K_{-\alpha_j}(y)+ xE_j(y), F_j(xy)=F_j(x) y+ K_{\alpha_j}xF_j(y)$$
for $j=1,2,\ldots,k-1$.

Simple objects of $\OO_f$ for $U_{q }(gl_d(\CC)\times gl_k(\CC))$
are labeled by pairs $(\lambda;\mu)$, where $\lambda$ is a dominant
integral $U_{q }(gl_d(\CC))$-weight and $\mu$ is a dominant integral
$U_{q }(gl_k(\CC))$-weight. Each irreducible   $U_q(gl_d(\CC)\times
gl_k(\CC))$-module of the highest weight  $(\lambda;\mu)$ is
isomorphic to $V'_{\lambda}\otimes V_\mu$, where $V'_\lambda$ is an
irreducible $U_q(gl_d(\CC)\times gl_k(\CC))$-module trivially
extended from that for $U_q(gl_d(\CC))$, and $V_\mu$ is an
irreducible $U_q(gl_d(\CC)\times gl_k(\CC))$-module trivially
extended from that for $U_q(gl_k(\CC))$.

  The decomposition  of $\CC_q[M_{d\times k}]$ into simple submodules is given by the following

\begin{proposition}\cite[Theorem 1.1]{Zh}
\label{pr: Howe-Duality}
 If $d\le k$, then for any $n\ge 0$ one has  an  isomorphism of $U_q(gl_d(\CC)\times gl_k(\CC))$-modules:
\begin{equation}
\label{eq:d-k decomposition}
\CC_q[M_{d\times k}]_n\cong\oplus_{\lambda} V'_{\lambda}\otimes V_{\tilde \lambda}
\end{equation}
where the summation is over all $d$-tuples of integers
$\lambda=(\lambda_1\ge \lambda_2\ge \cdots \ge \lambda_d\ge 0)$ such
that $\lambda_1+ \lambda_2+ \cdots + \lambda_d=n$; and  $\tilde\lambda\in\ZZ_{\ge 0}^k$ stands for the $\lambda$ extended with $k-d$ zeros.
\end{proposition}

%The algebra $\CC_{q}[M_{d\times k}]$  has the structure of a $U_{q}(\gl_{d})$-module as well as of a %$U_{q}(\gl_{d}\times \gl_{k})$-module. We have the following results.

%Morover, the above decomposition is compatible with the algebra structure.

Recall that the {\it braided} tensor product $A\otimes_{\RR} B$ of algebras $A$ and $B$ in $\OO_{gr,f}$ is defined as follows:

As a $U_q(\gg)$-module,  $A\otimes_{\RR} B\cong A\otimes B$, and the multiplication is twisted by $\RR$:
$$(a\otimes b)\cdot (a'\otimes b')=a\RR(b\otimes a')b$$
for any $a,a'\in A$, $b,b'\in B$. Associativity follows from that $\RR$ is a braiding.

The following result is an extension of \cite[Theorem 1.1]{Zh}.

\begin{proposition}
\label{pr:sum of products}
 As a $U_q(gl_d(\CC))$-module algebra, $\CC_q[M_{d\times k}]$ is isomorphic to the  braided  tensor $k$-th power
 of the $q$-polynomial algebra $S_q[x_1,\ldots,x_d]$ in the braided monoidal category $\OO_{gr,f}$ for $U_q(gl_d(\CC))$. More precisely,
the correspondence:
\begin{equation}
\label{eq:xij correspondence}
x_{ij}\mapsto 1\otimes 1\otimes\cdots  \otimes x_i\otimes \cdots \otimes 1 \otimes 1\ ,
\end{equation}
 where $x_i$ occurs only in the $j$-th place, $i=1,\ldots,d$, $j=1,\ldots,k$,
defines an isomorphism of algebras in $\OO_{gr,f}$ for $U_q(gl_d(\CC))$:
\begin{equation}
\label{eq:braided power}
\CC_q[M_{d\times k}]\widetilde \to \CC_q[x_1,\ldots,x_d]\otimes_\RR \CC_q[x_1,\ldots,x_d]\otimes _\RR \cdots \otimes _\RR \CC_q[x_1,\ldots,x_d]\ .
\end{equation}
\end{proposition}

\proof Since both algebras are quadratic, it suffices to show that
they have the same quadratic relations. For each $x\in
\CC_q[x_1,\ldots,x_d]$ denote $[x]_j$ the $k$-tensor of the form
$1^{\otimes j-1}\otimes x\otimes 1^{\otimes k-j}$. In particular,
the right hand side of (\ref{eq:xij correspondence}) is $[x_i]_j$.
Also for $x,y\in  \CC_q[x_1,\ldots,x_d]$ and for $1\le i<j\le k$
denote by $[x\otimes y]_{ij}$ the tensor $1^{\otimes i-1}\otimes
x\otimes 1^{\otimes j-i-1}\otimes y\otimes 1^{\otimes k-j}$ and
extend it by linearity: $[\sum_k x^{(k)}\otimes y^{(k)}]_{ij}=\sum_k
[x^{(k)}\otimes y^{(k)}]_{ij}$.

Clearly, $[x_j]_k[x_i]_k=[x_jx_i]_k=[qx_ix_j]_k=q[x_i]_k[x_j]_k$ for
any $i<j$ and  $[x_i]_k[x_j]_\ell=[\RR(x_i\otimes x_j)]_{k,\ell}$ if
$k>\ell$. Taking into account a  well-known computation:
\begin{equation}
\label{eq:braiding gl_d}
\RR(x_i\otimes x_j)=\begin{cases}
x_j\otimes x_i & \text{if $i<j$}\\
qx_i\otimes x_i & \text{if $i=j$}\\
x_j\otimes x_i+(q-q^{-1})x_i\otimes x_j & \text{if $i>j$} \\
\end{cases} \ ,
\end{equation}
we obtain for $k>\ell$
$$[x_i]_k[x_j]_\ell=\begin{cases}
[x_j\otimes x_i]_{\ell,k}=[x_j]_\ell[x_i]_k & \text{if $i<j$}\\
[qx_i\otimes x_i]_{\ell,k} =q[x_i]_k [x_i]_\ell & \text{if $i=j$}\\
[x_j\otimes x_i+(q-q^{-1})x_i\otimes x_j]_{\ell,k}=[x_j]_\ell[x_i]_k+(q-q^{-1})[x_i]_\ell[x_j]_k & \text{if $i>j$} \\
\end{cases}
$$
Therefore, the association $x_{ij}\mapsto [x_i]_j$ is a surjective
homomorphism of algebras. The fact that it is injective follows from
that both algebras are flat, i.e., both algebras are isomorphic to
$S(\CC(q)^{dk})$ as graded vector spaces.
\endproof

%%%%%%%%%%%%%%%%%%%%%%
According to  Proposition \ref{pr: Howe-Duality}, we have
\begin{equation}
\label{eq: k tensor  h/w}
(V_{\beta_1\omega_1}\otimes \cdots \otimes V_{\beta_k\omega_1})^\lambda=  V_{\lambda}(\beta)\ .
\end{equation}
for any $\beta=(\beta_1,\ldots,\beta_k)\in \ZZ_{\ge 0}^k$, where
$V_{\lambda}(\beta)$ is the weight space of the weight $(\beta)$ in
the irreducible $U_q(gl_k(\CC))$-module $V_\lambda$ in $\OO_f$, and
each $V_{n\omega_1}\cong S^n_\sigma V_{\omega_1}$ is a simple
$U_q(gl_d(\CC))$-module.

By definition, $V_{\ell_1\omega_1}\otimes  V_{\ell_2\omega_1}$ decomposes into the sum of its isotypic components:
\begin{equation}
\label{eq:isotypic decomposition multiplicity free}
V_{\ell_1\omega_1}\otimes  V_{\ell_2\omega_1}=\bigoplus_{m=0}^{\min(\ell_1,\ell_2)} J_{\ell_1,\ell_2}^m \ ,
\end{equation}
where $J_{\ell_1,\ell_2}^m\cong V_{(\ell_1+\ell_2-m,m,0,\ldots,0)}$
as a $U_q(gl_k(\CC))$-module (with the convention that
$J_{\ell_1,\ell_2}^m$ if $m\notin [0,\min(\ell_1,\ell_2)]$).

\begin{lemma}
\label{le:isotypic gl2} Under the identification (\ref{eq: k tensor
h/w}) we have for each $i\in [1,n-1]$,
 $\beta=(\beta_1,\ldots,\beta_k)\in \ZZ_{\ge 0}^k$, and $m\ge 0$:
\begin{equation}
\label{eq:isotypic gl2} (V_{\beta_1\omega_1}\otimes  \cdots \otimes
V_{\beta_{i-1}\omega_1}\otimes J_{\beta_i,\beta_{i+1}}^m\otimes
\cdots \otimes V_{\beta_k\omega_1})^\lambda=
E_i^m(V_{\lambda}(\beta-m\alpha_i)\cap Ker F_i) \ .
\end{equation}

\end{lemma}

\proof  We will prove (\ref{eq:isotypic gl2}) by induction in
$\beta_i\ge 0$. If $\beta_i=0$, the statement is obvious because
$V_0\otimes  V_{\beta_{i+1}\omega_1}=J^0_{0,\beta_{i+1}}\cong
V_{\beta_{i+1}\omega_1}$ and $V_{\lambda}(\beta)\subset Ker F_i$.

It is easy to see that under the identification $\bigoplus\limits_{\ell_1,\ell_2\ge 0}
V_{\ell_1\omega_1}\otimes  V_{\ell_2\omega_1}=\CC_q[M_{d\times 2}]$, we have $E_1(J_{\ell_1,\ell_2}^m)=J_{\ell_1+1,\ell_2-1}^m$
for all $\ell_1,\ell_2,m\ge 0$. Therefore, applying $F_i$ to (\ref{eq:isotypic gl2}) with  $(\beta_1,\ldots,\beta_i-1,\beta_{i+1}+1,\ldots,\beta_k)$,
we obtain (\ref{eq:isotypic gl2}) with $(\beta_1,\ldots,\beta_i,\beta_{i+1},\ldots, \beta_k)$. \endproof

%%%%%%%%%%%%%%%%%%%%%%

%For a $U_{q}(\gg)$-module $V$ in $\OO_{f,gr}$ and for a dominant weight $\lambda\in P_+$ we denote  $V^\lambda=\{v\in V(\lambda):E_i(v)=0 ~\forall %i\}$,  the space of highest weight vectors   of weight $\lambda$.

\end{document}